\documentclass[twocolumn]{autart_arxiv}

\usepackage{soul} % for strikethrough
\usepackage{caption}

\usepackage{pgfplots}
\pgfplotsset{compat=newest}
\tikzset{every picture/.style={font issue=\scriptsize},
         font issue/.style={execute at begin picture={#1\selectfont}}
        } 

\usepackage[tight,normalsize,sf,SF]{subfigure}

\usepackage{csquotes}
\usepackage{amssymb}
\usepackage{amsmath}
\usepackage{graphicx}
\graphicspath{ {fig/} }
\usepackage{epstopdf}

\begin{document}

\begin{frontmatter}

\title{On Transitive Consistency for Linear Invertible Transformations 
between Euclidean Coordinate Systems\thanksref{footnoteinfo}}

\thanks[footnoteinfo]{The authors gratefully acknowledge the financial support from the Fonds
National de la Recherche, Luxembourg (6538106, 8864515).}

\author[Johan]{Johan Thunberg}
\ead{johan.thunberg@uni.lu},
\author[Johan,Florian1]{Florian Bernard}\ead{bernard.florian@chl.lu},
\author[Johan,Cambrdige]{Jorge Goncalves}
\ead{jorge.goncalves@uni.lu}

\address[Johan]{Luxembourg Centre for Systems Biomedicine, University of Luxembourg, Esch-sur-Alzette, LUXEMBOURG}

\address[Florian1]{
Centre Hospitalier de Luxembourg, Luxembourg City, LUXEMBOURG}

\address[Cambrdige]{Control Group, Department of Engineering, University of Cambridge, Cambridge, UNITED KINGDOM}

\begin{keyword}
Distributed optimization, transformation synchronization, Procrustes problem, consensus algorithms, graph theory.
\end{keyword}

\begin{abstract}
Transitive consistency is an intrinsic property
for collections of linear invertible transformations
between Euclidean coordinate frames. In practice,
when the transformations are estimated from data, this
property is lacking.
This work addresses the problem of synchronizing transformations that are not transitively consistent. Once 
the transformations have been synchronized, they satisfy the transitive consistency condition -- a transformation from frame $A$ to frame $C$ is equal to the 
composite transformation of first transforming 
$A$ to  $B$ and then transforming $B$ to $C$. 
The coordinate frames correspond to nodes in 
a graph and the transformations correspond to  edges in the same graph.
Two direct or centralized synchronization methods are presented
for different graph topologies; the first one for quasi-strongly connected graphs, and the second one for connected graphs. As an extension of the second method, an iterative Gauss-Newton method is presented, which is later adapted to the case of affine and Euclidean transformations.
Two distributed synchronization methods are also presented for orthogonal matrices, which can be seen as distributed versions of the two direct or centralized methods; they are similar in nature to standard
consensus protocols used for distributed averaging.
When the transformations are orthogonal matrices,
a bound on the optimality gap can be computed. Simulations
show that the gap is almost tight, even for noise
large in magnitude. 
This work also contributes on a theoretical level by providing 
linear algebraic relationships for 
transitively consistent transformations. 
One of the benefits of the proposed 
methods is their simplicity -- basic linear algebraic methods are used, e.g., the Singular Value Decomposition (SVD). For a wide range of parameter settings, the methods are numerically validated.
\end{abstract}

\end{frontmatter}

\section{Introduction}\label{sec:introduction}
Collections of linear invertible transformations between Euclidean 
coordinate systems must be transitively consistent. 
In practice however, when the transformations are estimated
from data, this condition does not hold. 
This issue is present in the 3D localization problem, where transformations are rigid and estimated from e.g.,
camera measurements; in the multiple images registration
problem where the transformations are affine (or linear 
by using homogeneous coordinates); in the generalized Procrustes problem where scales, rotations and translations are calculated from multiple  point clouds. In order to resolve the issue,
the estimated transformations need to be synchronized in the sense of finding transitively consistent transformations close to the estimated ones.

\subsection{Problem}
This work addresses the problem of
synchronizing linear invertible transformations or matrices between 
Euclidean coordinate systems or frames. More precisely,
given a collection $\{G_{ij}\}$ of matrices in $GL(d, \mathbb{R})$, another collection $\{G_{ij}^*\}$ of
matrices in $GL(d, \mathbb{R})$ is constructed such that 
\begin{equation}\label{eq:8}
G^*_{ij}G^*_{jk} = G^*_{ik}, \text{ for all } i,j,k,
\end{equation}
where $G_{ij}^*$ is \enquote{close} to $G_{ij}$ for all $i,j$. 
By satisfying \eqref{eq:8}, the collection
$\{G_{ij}^*\}$ is said to be {transitively consistent}.

\subsection{Background}
There are many applications for the proposed methods.  
One such application is the 3D localization problem in camera networks~\cite{tron2014distributed} where a network of cameras are observing a scene and epipolar geometry is used to calculate/measure $G_{ij}$ transformations between $(i,j)$-pairs of cameras. If the cameras are fully calibrated these transformations are Euclidean, otherwise they could be e.g., affine (or linear by using homogeneous coordinates). Since the transformations are calculated
from measurements, they do not satisfy \eqref{eq:8} in general. Hence our proposed methods can be used
to synchronize the matrices. 
For the 3D localization problem,
we do not have to limit ourselves to the case of cameras
and epipolar geometry. The transformations could be 
calculated in a setting where the geometry of the scene is known. In  the case of known point features, the perspective-n-point problem
can be solved in order to get estimates of the relative transformations~\cite{lepetit2009epnp}.

Another important problem is image registration, which has attracted much attention in the 
medical imaging community. The number of applications
is vast, ranging from surgery planning to longitudinal studies. To register a (moving) image with another (fixed)
image is to transform the former into the the latter in such 
way that they fit in the \enquote{best} way. For that, optimization methods are used 
to calculate a transformation which minimizes a suitable 
objective. Registration of multiple
images poses a greater challenge. There are several approaches
in the literature. For example: Finding a path of pairwise transformations, which contains all images \cite{vskrinjar2008symmetric}; aligning images
with a reference frame \cite{joshi2004unbiased}; image congealing, where variability along known axes of variation is removed in an 
iterative manner \cite{joshi2004unbiased}; considering a minimum description length (MDL) approach of a statistical shape model built from the correspondences given due to groupwise image registration \cite{cootes2004groupwise}; Bayesian formulations and 
Expected Maximization (EM) \cite{allassonniere2007towards}.
 
Another way to solve
the (affine) multiple images registration problem is to use the 
transitive consistency criterion \eqref{eq:8} \cite{bernard_sami}. Let 
the $G_{ij}$ correspond to the affine transformations calculated 
from pairwise registrations, then our method can be 
used to create transitively consistent 
$G_{ij}^*$ transformations. Registration methods using transitive consistency have also been proposed for deformable transformations \cite{gass2014detection,geng2007transitive}.

A related problem to the one posed in this paper is the problem of calculating 
the \enquote{best} translations, rotations and scales between 
pairs of point clouds. If only one pair is considered the 
problem is referred to as the Procrustes Problem~\cite{Gower:2004uu}. This problem can be solved by means of
singular value decomposition or eigenvalue decomposition~\cite{Arun:1987uu,Schonemann:1966ch,Horn:1988bq}, or in the 
case case of 3D transformations, by a quaternion-based approach \cite{Horn:1987hf,Walker:1991kt}. The problem
restricted to 3D is referred to as the absolute orientation problem~\cite{Horn:1988bq,Horn:1987hf}.
In the general setting, when $n$ point clouds are considered,
the problem is referred to as the Generalized Procrustes 
Problem~\cite{Gower:2004uu}. In order to solve this problem,
iterative methods are often used; when the dimension is 
two or three, direct methods have recently been proposed~\cite{Pizarro:2011ta}. 
Our previous work in \cite{bernard_cvpr} has tackled the Generalized Procrustes Problem 
using an approach based on transitive consistency. The present paper will extend and generalize these ideas as well
as describe many theoretical properties of the generalizations.

Our methods can be used for solving the Generalized Procrustes Problem in the following way: Between each 
pair of point clouds a $G_{ij}$ transformation is calculated
using any standard technique~\cite{Arun:1987uu,Schonemann:1966ch,Horn:1988bq}, then our methods are used to improve
the pairwise transformations by calculating transitively consistent transformations.

In the special case when the $G_{ij}$
are orthogonal matrices, 
Singer \emph{et al.}  have presented methods for the optimization
of transitive consistency~\cite{Singer:2011ba,Hadani:2011hb,Hadani:2011tw,Chaudhury:2013un}. These works were later adapted 
by Pachauri \emph{et al.} to the special case when the $G_{ij}$ are permutation 
matrices~\cite{Pachauri:2013wx}. In the latter work, a 
relaxation of the original problem is considered -- in the original problem the transformations 
shall be orthogonal matrices --
and then 
permutation matrices are obtained by means of projection from the solution of the relaxed problem.
The method presented by Singer \emph{et al.}  is said 
to be a synchronization method for minimization of
transitive consistency errors -- a formalism adopted in this work.

\subsection{Methods and results}
The approach in this work share similarities with the approaches of Singer \emph{et al.} and Pachauri \emph{et al.}; it continues along
the lines of the the recently proposed methods in \cite{bernard_cvpr,thunberg_synbchronization}. 

In \cite{bernard_cvpr,thunberg_synbchronization} a so called
$Z$-matrix is constructed from the $G_{ij}$ matrices. If
the index set for the (available) transformations has a certain property, transitively
consistent transformations can be obtained by a method 
where the Singular Value Decomposition (SVD) is calculated 
for $Z$. The property that must be fulfilled for the index set
$\{(i,j)\} = \mathcal{E}$, is that it is the edge set of a quasi-strongly connected (QSC)
directed graph (see Definition~\ref{qscGraph}). In the $Z$-matrix approach, 
a set of linear algebraic equations are formulated -- equations which 
shall be satisfied for the case of transitively consistent
transformations. When the transformations are not transitively 
consistent, the problem is solved in the sense of least
squares minimization.

As we will show in this work, the $Z$-matrix appears in the construction of 
a Hessian matrix $H$ for a quadratic convex function
of the $G_{ij}$,
and under certain conditions it holds that $H = Z + Z^T$. 
From the SVD of the Hessian matrix $H$, transitively consistent
transformations can be calculated in the same manner as 
for the $Z$-matrix. The justification for using the $H$-matrix  stems from the fact that
it is the Hessian matrix of the objective function in a relevant optimization problem. The justification of using 
the $Z$-matrix stems purely from the linear algebraic constraints
that should be satisfied for transitively consistent transformations.

The $Z$-matrix method and the $H$-matrix method are both direct methods,
i.e., the solution is found at once. As an extension we also 
propose an iterative Gauss-Newton method, which uses the solution 
from the $H$-matrix method as initialization. For orthogonal 
matrices one can prove that this iterative scheme cannot decrease
the objective function at all. The Gauss-Newton method is also adapted to the cases of affine and Euclidean transformations. In this case -- as opposed to the result for orthogonal matrices -- significant 
improvement over the $Z$-matrix method and the $H$-matrix method can be seen in numerical simulations.

Many properties of the $Z$-matrix and the $H$-matrix are 
proved in this work. For example it is shown that transitive consistency 
in the case of connected graphs is equivalent to the condition 
that the nullspace of $H$ has dimension $d$. Furthermore,
the transitively consistent transformations can be obtained as the $d \times d$ blocks in a
matrix, the columns of which span the nullspace of $H$. 
For the $Z$-matrix only a weaker condition is  formulated;
if the graph is QSC and the transformations are 
transitively consistent, the transformations can be obtained as the $d \times d$ blocks in a
matrix, the columns of which span the nullspace of $Z$.

Now, in most aspects the $H$-matrix approach seems to be superior to the
$Z$-matrix approach. However, one large benefit of using the $Z$-matrix
over the $H$-matrix is that it can be used in a distributed 
algorithm when the communication graph is directed. 

In a later part of the paper, two distributed methods are introduced for the case of orthogonal $G_{ij}$ transformations. The first method is using the $Z$-matrix under the assumption 
that the
communication graph is directed and QSC. The other method is using the $H$-matrix under the assumption that the communication
graph is symmetric. The performance of the two methods are almost the same in numerical experiments. The distributed 
methods are similar in structure to linear consensus protocols~\cite{mesbahi2010graph,jad03,sab04,olfati2007consensus,olfati2004consensus}. Key differences to those approaches is that the states
here are matrices instead of vectors, and the states combined converge
to a $d$-dimensional linear subspace instead of the 
consensus set.   

The distributed iterative methods are introduced 
mainly with communication between agents in mind, e.g., in networks
of robots with limited communication range, where the robots only communicate with their neighbors (directly or indirectly). However, a further scenario of the distributed methods is parallelisation in order to better deal with the computational burden in the case of very large problem instances. 

When it is known that the transitively consistent transformations
are orthogonal matrices, i.e., elements of $O(d) = \{R: R \in \mathbb{R}^{d \times d}, R^TR = I\}$, 
a method is provided for calculating an upper
bound on the optimality gap. In the case when the $G_{ij}$ are
also orthogonal, simulations show that this gap is almost
tight. As an example, for $n = 100$ coordinate systems,  dimension $d = 3$, and randomly generated $G_{ij}$ matrices
in $O(3)$, the gap is smaller than a tenth of a percent in average.   
There
are (and will be even more in the future) applications where 
large networks of cameras, robots, satellites or unmanned vehicles, need to synchronize their pairwise relative rotations. 
In such applications methods that are near optimal and run almost in real time are of utmost importance 
to have.

\subsection{Outline}
The paper proceeds as follows. In Section~\ref{sec:2}, graphs
and properties thereof are introduced, followed by the introduction of the $G_{ij}$ transformations and their connections to the graphs.
We have chosen to incorporate 
graphs in the very definition of transitive consistency. 
Section~\ref{sec:33} addresses linear invertible transformations. In Section~\ref{sec:zMatrix}, the $Z$-matrix is introduced,
followed by a collection of results and a least squares 
method. In Section~\ref{sec:hMatrix}, the $H$-matrix is introduced;
in the same manner as in Section~\ref{sec:zMatrix}, a collection
of results is provided in conjunction with an algorithm. In Section~\ref{sec:gn} a Gauss-Newton method is presented, where the matrices obtained from the $H$-matrix method are used as initialization. 
Section~\ref{sec:44} consider the special case of orthogonal matrices, i.e., elements of 
$O(d)$.
The section starts with some bounds on the optimality
gap, and continues in Section~\ref{sec:5} with the introduction of distributed algorithms. The reader interested in the 
distributed methods can go directly to this section and consult the earlier sections only for reference. Section~\ref{sec:6} is a small detour, where 
a gradient flow method is presented for orthogonal matrices. This method
is employed as a baseline method, used for comparison
in some of the simulations in Section~\ref{sec:results} -- the
section where the proposed methods are thoroughly numerically evaluated.

\section{Preliminaries}\label{sec:2}
\subsection{Directed Graphs}
Let $\mathcal{G} = (\mathcal{V}, \mathcal{E})$ be a directed graph,
where $\mathcal{V} = \{1, 2, \ldots, n\}$ is the node set and 
$\mathcal{E} \subset \mathcal{V} \times \mathcal{V}$ is the 
edge set. The set $\mathcal{N}_i$ is defined by 
$$\mathcal{N}_i = \{j: (i,j) \in \mathcal{E}\}.$$

The adjacency matrix $A = [A_{ij}]$ for the graph 
$\mathcal{G}$ is defined by
$${A}_{ij}= \begin{cases}
1  &  \text{ if } (i,j) \in \mathcal{E},\\
0 & \text{ else.}
\end{cases}$$

The graph Laplacian matrix is defined by 
$$L = \text{diag}(A 1_n) - A,$$
where $1_n \in \mathbb{R}^n$ is a vector 
with all entries equal to $1$. In order to emphasize that 
the adjacency matrix $A$, the Laplacian matrix $L$ 
and the $\mathcal{N}_i$ sets depend on the graph $\mathcal{G}$, 
we may write $A(\mathcal{G})$, $L(\mathcal{G})$ and $\mathcal{N}_i(\mathcal{G})$ respectively. For simplicity however, we mostly omit this notation
and simply write $A$, $L$, and $\mathcal{N}_i$.

\begin{defn} \emph{(connected graph, {undirected} path)}\\
The directed graph $\mathcal{G}$ is connected 
if there is an {undirected} path from any node
in the graph to any other node. An {undirected} path is defined as a (finite) sequence of unique nodes such that for any pair $(i,j)$ of consecutive nodes in the sequence it holds that
$$((i,j) \in \mathcal{E}) \text{ or } ((j,i) \in \mathcal{E}).$$
\end{defn}

\begin{defn} \emph{(quasi-strongly connected graph, center, directed path)}\label{qscGraph}\\
The directed graph $\mathcal{G}$ is quasi-strongly connected (QSC)
if it contains a center. A center is
a node in the graph to which there is a directed path from any other node
in the graph. A directed path is defined as a (finite) sequence of unique nodes such that any pair of consecutive nodes in the sequence comprises an edge in 
$\mathcal{E}$.
\end{defn}

\begin{defn}\emph{(symmetric graph)}\label{def:3}\\
The directed graph $\mathcal{G} = (\mathcal{V}, \mathcal{E})$ is {symmetric} if $$((i,j) \in \mathcal{E}) \Rightarrow ((j,i) \in \mathcal{E}) \text{ for all } (i,j) \in \mathcal{V} \times \mathcal{V}. $$
\end{defn}

Given a graph $\mathcal{G} = (\mathcal{V}, \mathcal{E})$, the graph $\bar{\mathcal{G}} = (\mathcal{V}, \bar{\mathcal{E}})$ is the graph constructed by reversing the direction
of the edges in $\mathcal{E}$, i.e., $(i,j) \in \bar{\mathcal{E}}$ if and only if
$(j,i) \in {\mathcal{E}}$. 
It is easy to see that 
$$A(\bar{\mathcal{G}}) = (A({\mathcal{G}}))^T \text{ and } L(\bar{\mathcal{G}}) =  \text{diag}((A({\mathcal{G}}))^T 1_d) - A({\mathcal{G}})^T.$$

\subsection{Transformations}
Given a directed graph $\mathcal{G} = (\mathcal{V}, \mathcal{E})$, let there be a collection of matrices $\{{G}_{ij}\}_{(i,j) \in \mathcal{E}}$ where 
${G}_{ij} \in GL(d, \mathbb{R})$ for all $(i,j) {\in \mathcal{E}}$.
 Let $n = |\mathcal{V}|$. The ${G}_{ij}$ are 
not necessarily transitively consistent in that 
$${G}_{ik} \neq {G}_{ij}{G}_{jk}$$
may hold if $(i,j), (j,k)$ and $(i,k)$ are elements of $\mathcal{E}$.

In the methods to be defined, the goal is to find a transitively consistent collection
$\{{G}^*_{ij}\}_{(i,j) \in \mathcal{E}}$ of matrices in $GL(d, \mathbb{R})$, such that for all $(i,j) \in \mathcal{E}$, $G^*_{ij}$ is 
close to $G_{ij}$ in some appropriate sense.
Notation-wise, $G_{ij}^*$ is simply (a name of) a matrix. This notation
should not be mixed up with the conjugate transpose -- in this paper, all matrices considered 
are real and the conjugate transpose will not be used.

\begin{defn} \emph{(transitive consistency)}
\begin{enumerate}
\item The matrices in the collection $\{{G}^*_{ij}\}_{(i,j) \in \mathcal{V} \times \mathcal{V}}$ of matrices in $GL(d, \mathbb{R})$ are \emph{transitively 
consistent for the complete graph} if 
$${G}^*_{ik} = {G}^*_{ij}{G}^*_{jk}$$
for all $i,j$ and $k$. \\
\item Given a graph $\mathcal{G} = (\mathcal{V}, \mathcal{E})$, the matrices in the collection 
$\{{G}^*_{ij}\}_{(i,j) \in \mathcal{E}}$ of matrices in $GL(d, \mathbb{R})$ are  \emph{transitively 
consistent for $\mathcal{G}$} if there is a collection
$\{{G}^*_{ij}\}_{(i,j) \in \mathcal{V} \times \mathcal{V}} \supset \{{G}^*_{ij}\}_{(i,j) \in \mathcal{E}}$
such that $\{{G}^*_{ij}\}_{(i,j) \in \mathcal{V} \times \mathcal{V}}$ is transitively consistent for the complete graph.
\end{enumerate}

\end{defn}
If it is apparent by the context, sometimes we will be less strict and omit to mention which graph a collection of transformations is transitively consistent for.
 A sufficient
condition for transitive consistency of the $G^*_{ij}$ matrices for
any graph is that
there is a collection $\{G^*_i\}_{i \in \mathcal{V}}$ of matrices in $GL(d, \mathbb{R})$ such that 
$$G^*_{ij} = G^{*-1}_{i}G^*_{j}$$
for all $i,j$. Lemma~\ref{lem:1} below and the proof thereof provides additional important information. 
The result is similar to that in \cite{tron2014distributed}.
For the statement of the lemma, the following 
definition is needed. 

\begin{defn}
Two collections $\{G^*_i\}_{i \in \mathcal{V}}$ and 
$\{G^{**}_{i}\}_{i \in \mathcal{V}}$ of
matrices in $GL(d, \mathbb{R})$ are equal up to transformation 
from the left, if there is $Q \in GL(d, \mathbb{R})$ such that 
$$QG^*_i = G^{**}_i \text{ for all  } i.$$
\end{defn}

\begin{lem}\label{lem:1}
For any graph $\mathcal{G} = (\mathcal{V}, \mathcal{E})$ and collection $\{{G}^*_{ij}\}_{(i,j) \in \mathcal{E}}$ of matrices
in $GL(d, \mathbb{R})$ that are transitively consistent for $\mathcal{G}$,
\begin{enumerate}
\item there is a collection $\{G^{*}_{i}\}_{i \in \mathcal{V}}$ of matrices in $GL(d, \mathbb{R})$ such that  
\begin{equation}\label{eq:1}
G^*_{ij} = G^{*-1}_{i}G^*_{j} \text{ for all } (i,j) \in \mathcal{E},
\end{equation}
\item all collections $\{ G^*_i\}_{i \in \mathcal{V}}$ satisfying \eqref{eq:1} are equal up
to transformation from the left
if and only if $\mathcal{G}$ is connected,
\item \vspace{1mm} there is a unique collection $\{{G}^*_{ij}\}_{(i,j) \in \mathcal{V} \times \mathcal{V}} \supset \{{G}^*_{ij}\}_{(i,j) \in \mathcal{E}}$ of
transitively consistent matrices for the 
complete graph, if and only if all collections $\{ G^*_i\}_{i \in \mathcal{V}}$ satisfying \eqref{eq:1} are equal up
to transformation from the left.

\end{enumerate}
\end{lem}

\noindent\emph{Proof: }
All matrices appearing in this proof, if the contrary is not explicitly stated,
are assumed to be elements of $GL(d, \mathbb{R})$.\vspace{2mm}\\ 
\noindent
{(1)} Since the matrices in $\{{G}^*_{ij}\}_{(i,j) \in \mathcal{E}}$ are transitively consistent for $\mathcal{G}$, there
is $\{{G}^*_{ij}\}_{(i,j) \in \mathcal{V} \times \mathcal{V}} \supset \{{G}^*_{ij}\}_{(i,j) \in \mathcal{E}}$ in
which the matrices are transitively consistent for the complete 
graph.
Let the $G^*_i$ in a collection $\{G^*_i\}_{i \in \mathcal{V}}$ be defined by
 $$G^*_i = G^*_{1i}.$$
We  shall prove that $$G^{*-1}_iG^*_j = G^{*-1}_{1i}G^*_{1j} = G^*_{ij} \text{ for all } i,j.$$
Using the fact that $G^*_{11}$ is invertible and the fact that 
$G^{*}_{11} = G^{*2}_{11}$, one can show that $G^{*2}_{11}=I$. Now, $G^*_{1i}G^*_{i1} = G_{11} = I$; thus  $G^{*-1}_{1i} = G^{*}_{i1}$.
But then 
$$G^{*-1}_{1i}G^*_{1j} = G^{*}_{i1}G^*_{1j} = G^{*}_{ij}.$$
\vspace{2mm}

\noindent
{(2)} We know that transitive consistency of 
$\{{G}^*_{ij}\}_{(i,j) \in \mathcal{E}}$ for $\mathcal{G}$ is equivalent 
to 
the statement that there is a collection $\{G^*_i\}_{i \in \mathcal{V}}$ of matrices in $GL(d, \mathbb{R})$ such that  
$$G^*_{ij} = G^{*-1}_{i}G^*_{j} \text{ for all } (i,j) \in \mathcal{E}.$$
Let $G_{ij}^* = G_i^{-1}G_i$ for all 
$i,j \in \mathcal{V}$.
For any other collection 
$\{G^{**}_i\}_{i \in \mathcal{V}}$ of matrices in $GL(d, \mathbb{R})$ such that  
$$G^*_{ij} = G^{**-1}_{i}G^{**}_{j} \text{ for all } (i,j) \in \mathcal{E},$$
it holds that 
\begin{align*}
& \begin{bmatrix}
G^{**{T}}_{1} & G^{**{T}}_{2} & \ldots & G^{**{T}}_{n}
\end{bmatrix}^{T} \\
=~& \text{diag}\left(
Q_{1}, Q_{2}, \ldots,  Q_{n} \right) \cdot 
\begin{bmatrix}
G^{*{T}}_{1} & G^{*{T}}_{2} & \ldots & G^{*{T}}_{n}
\end{bmatrix}^{T},
\end{align*}
where the $Q_i$ matrices are elements of $GL(d, \mathbb{R})$. 

{\textbf{If}:} Now, if the graph is connected and at least two of the $Q_i$ are not equal, 
there is $(j,k) \in \mathcal{E}$ such that $Q_j \neq Q_k$. 
We know $$G^{*}_{j}G^*_{jk}G^{*-1}_{k} = I_d,$$
but since $Q_j \neq Q_k$ we can calculate this entity to 
$$G^{*}_{j}G^*_{jk}G^{*-1}_{k} = Q_j^{-1}Q_k \neq I_d,$$
which is a contradiction. $I_d \in \mathbb{R}^{n \times n}$ is the identity matrix.

{\textbf{Only if}:} On the other hand, if the graph is not connected there are 
two disjoint sets $\mathcal{V}_1$ and $\mathcal{V}_2$ such that 
$\mathcal{V}_1 \cup \mathcal{V}_2 = \mathcal{V}$, for which there is 
no pair $(i,j) \in \mathcal{E}$ such that ($i \in \mathcal{V}_1$ and 
$j \in \mathcal{V}_2$) or ($j \in \mathcal{V}_1$ and 
$i \in \mathcal{V}_2$). Thus, the nodes in $\mathcal{V}_1$ and the 
corresponding edges, respective the nodes in $\mathcal{V}_2$ and the 
corresponding edges, can be seen as two different disconnected {(sub)}graphs, {each of them being connected};
the $G^*_i$ matrices in the first graph can be multiplied 
with a matrix $Q_1$ from the left and the  $G^*_i$ matrices in the 
second graph can be multiplied with a matrix $Q_2$ from the left, where $Q_1 \neq Q_2$, generating a collection of matrices 
$\{G^{**}_{i}\}_{i \in \mathcal{V}}$  not equal to $\{G^{*}_{i}\}_{i \in \mathcal{V}}$ up to transformation from the left.
\vspace{2mm}

\noindent
{(3)} \textbf{If:} Any other collection 
$\{G^{**}_{i \in \mathcal{V}}\}$ of matrices in $GL(d, \mathbb{R})$ such that  
$$G^*_{ij} = G^{**-1}_{i}G^{**}_{j} \text{ for all } (i,j) \in \mathcal{E},$$
is equal to $\{G^*_i\}$ up to transformation from the left. 
Now, for any $(i,j) \in (\mathcal{V}\times \mathcal{V}) - \mathcal{E}$ it holds that 
$$G^{**-1}_{i}G^{**}_{j} = G^{*-1}_{i}Q^{-1}QG^{*}_{j} = G^{*-1}_{i}G^{*}_{j} = G^{*}_{ij}$$ for some matrix $Q \in GL(d, \mathbb{R})$. 

\vspace{2mm}
\noindent
 \textbf{Only if:} The approach here is similar
 to that in 2) above. Suppose for $\{G^{*}_{i \in \mathcal{V}}\}$ satisfying \eqref{eq:1}, there is another collection 
$\{G^{**}_{i \in \mathcal{V}}\}$ of matrices in $GL(d, \mathbb{R})$ also satisfying \eqref{eq:1}, but the matrices in the two collections are 
not equal up to transformation from the left. Then it holds that 
\begin{align*}
& \begin{bmatrix}
G^{**{T}}_{1} & G^{**{T}}_{2} & \ldots & G^{**{T}}_{n}
\end{bmatrix}^{T} \\
=~& \text{diag}\left(
Q_{1}, Q_{2}, \ldots,  Q_{n} \right) \cdot 
\begin{bmatrix}
G^{*{T}}_{1} & G^{*{T}}_{2} & \ldots & G^{*{T}}_{n}
\end{bmatrix}^{T},
\end{align*}
where the $Q_i$ matrices are elements of $GL(d, \mathbb{R})$ and there is a pair $(k,l)$ for which $Q_k \neq Q_l$.

Now 
$$G^{**-1}_kG^{**}_l = G^{*-1}_kQ_k^{-1}Q_lG^{*}_l \neq  G^{*-1}_kG^{*}_l.$$
\hfill $\blacksquare$\newline

Lemma~\ref{lem:1} states that connectivity is a necessary
property to determine a unique (up to transformation 
from the left) collection $\{G_i^*\}$ satisfying 
\eqref{eq:1}. As it turns out, a stronger type
of connectivity -- quasi-strong connectivity -- is
useful in order to develop linear algebraic 
methods for solving our synchronization problem. The
first method we present is based on the so called $Z$-matrix.

\section{Linear invertible transformations}\label{sec:33}
\subsection{The $Z$-matrix}\label{sec:zMatrix} 
In this section a certain matrix is defined -- referred to as $Z$. 
It is used as a building block in a matrix $H$,
corresponding to the Hessian of a convex quadratic function, see Section~\ref{sec:hMatrix}.
After its definition, its properties are investigated. Amongst other things,
it is shown that if the $G_{ij}$ transformations are orthogonal, i.e., $G_{ij}^TG_{ij} = I$, the matrix $(-Z)$ is (critically) stable in the linear dynamical systems sense (cf. Lemma \ref{lem:10}). This means that, for directed graphs, the matrix $Z$ can be used in a linear distributed
algorithm for synchronizing orthogonal matrices (Section~\ref{sec:5}).

Define the matrix 
\begin{align*}
W(\mathcal{G}, \{{G}^*_{ij}\}_{(i,j) \in \mathcal{E}}) = [W_{ij}(G_{ij}^*)],
\end{align*}
where $$W_{ij}(G_{ij}^*) = \begin{cases}
G^*_{ij} & \text{ if } j \in \mathcal{N}_i, \\
0 & \text{ else, }
\end{cases}$$
and the matrix 
\begin{align*}
& Z(\mathcal{G}, \{{G}^*_{ij}\}_{(i,j) \in \mathcal{E}})  \\
 = & ~\text{diag}(A(\mathcal{G})1) \otimes I_d - W(\mathcal{G}, \{{G}^*_{ij}\}_{(i,j) \in \mathcal{E}}).
\end{align*}
The symbol $\otimes$ denotes the Kronecker product.

\begin{rem}
A more general way of constructing the $W$-matrix and 
the $Z$-matrix
with positive weights is as follows. Replace the 
$W_{ij}$ in the definition of $W$ with $a_{ij}W_{ij}$,
and replace $\text{diag}(A(\mathcal{G})1) \otimes I_d$ 
in the definition of $Z$ with 
$$\text{diag}(\sum_{j \in \mathcal{N}_1}a_{1j}, \sum_{j \in \mathcal{N}_2}a_{2j}, \ldots, \sum_{j \in \mathcal{N}_n}a_{nj}) \otimes I_d.$$
The $a_{ij}$ are positive for all $i,j$. Equivalent results to all the results
obtained for the $Z$-matrix in this section can also be formulated for the
alternative $Z$-matrix with 
positive weights. The alternative $Z$-matrix can be 
used in a distributed algorithm, equivalent to the 
one that will be presented in Section~\ref{sec:dist:orthogonal}.
\end{rem}

For the collection $\{G^*_i\}_{i \in \mathcal{V}}$ of matrices in $GL(d, \mathbb{R})$,
let 
\begin{align*}
U_1(\{G^*_i\}_{i \in \mathcal{V}}) & = \begin{bmatrix}
G^{*-T}_1 & G^{*-T}_2 & \ldots & G^{*-T}_n 
\end{bmatrix}^T, \\ 
U_2(\{G^*_i\}_{i \in \mathcal{V}}) & = \begin{bmatrix}
G^{*}_1 & G^{*}_2 & \ldots & G^{*}_n 
\end{bmatrix}. 
\end{align*}

\begin{lem}\label{lem:2}
For any \textbf{(QSC)} graph $\mathcal{G} = (\mathcal{V}, \mathcal{E})$,
collection $\{{G}^*_{ij}\}_{(i,j) \in \mathcal{E}}$ of matrices
in $GL(d, \mathbb{R})$ -- transitively consistent for $\mathcal{G}$ -- and collection $\{G^*_i\}_{i \in \mathcal{V}}$ of matrices in $GL(d, \mathbb{R})$ 
it holds that
\begin{equation*}
G^*_{ij} = G^{*-1}_{i}G^*_{j} \text{ for all } (i,j) \in \mathcal{E}
\end{equation*}
{\textbf{(if and)} only if}
\begin{equation}\label{eq:olle}
\text{im}(\text{diag}\left(
G^*_{1}, G^*_{2}, \ldots,  G^*_{n} \right)V) = \text{ker}(L \otimes I_d),
\end{equation}
for any matrix $V$, where the columns thereof form a basis for $\text{ker}(Z(\mathcal{G}, \{{G}^*_{ij}\}_{(i,j) \in \mathcal{E}}))$.
In particular, if  $\mathcal{G}$ is QSC,  \eqref{eq:olle}
can be stated as 
$$\text{im}(U_1(\{G^*_i\}_{i \in \mathcal{V}})) = \text{ker}(Z(\mathcal{G}, \{{G}^*_{ij}\}_{(i,j) \in \mathcal{E}})).$$
\end{lem}

\noindent\emph{Proof: } \\
\textbf{Only if:} Suppose it holds that  
$$G^*_{ij} = G^{*-1}_{i}G^*_{j} \text{ for all } (i,j) \in \mathcal{E}.$$
Then
\begin{align}
\label{eq:Z_laplacian}
& Z(\mathcal{G}, \{{G}^*_{ij}\}_{(i,j) \in \mathcal{E}}) \\
\nonumber
=~ &\text{diag}(G_1^{*-1}, G_2^{*-1}, \ldots, G_n^{*-1})(L \otimes I)\cdot \\
\nonumber
& \text{diag}(G^*_1, G^*_2, \ldots, G^*_n).
\end{align}
Now,
\begin{align*}
Z(\mathcal{G}, \{{G}^*_{ij}\}_{(i,j) \in \mathcal{E}})V & = 0 \Leftrightarrow \\
(L \otimes I) \text{diag}(G^*_1, G^*_2, \ldots, G^*_n)V & = 0 \Leftrightarrow \\
\text{im}(\text{diag}\left(
G^*_{1}, G^*_{2}, \ldots,  G^*_{n} \right)V) & = \text{ker}(L \otimes I_d).
\end{align*}

\noindent 
\textbf{If:} This part is only proven for the case when the 
graph $\mathcal{G}$ is QSC. 

Since $\{{G}^*_{ij}\}_{(i,j) \in \mathcal{E}}$ is transitively consistent 
for $\mathcal{G}$, there is $\{G^{**}_i\}_{i \in \mathcal{V}}$ of matrices in $GL(d, \mathbb{R})$ such that 
\begin{align*}
& Z(\mathcal{G}, \{{G}^*_{ij}\}_{(i,j) \in \mathcal{E}}) \\
\nonumber
=~ &\text{diag}(G_1^{**-1}, G_2^{**-1}, \ldots, G_n^{**-1})(L \otimes I)\cdot \\
\nonumber
& \text{diag}(G^{**}_1, G^{**}_2, \ldots, G^{**}_n).
\end{align*}
Thus, the null-space of $Z(\mathcal{G}, \{{G}^*_{ij}\}_{(i,j) \in \mathcal{E}})$
is given by
\begin{align*}
&~ \text{ker}(Z(\mathcal{G}, \{{G}^*_{ij}\}_{(i,j) \in \mathcal{E}})) = \text{im}(V),
\end{align*}
where $$V = \text{diag}(G^{-1**}_1, G^{-1**}_2, \ldots, G^{-1**}_n)([1, 1, \ldots, 1]^T \otimes I_d).$$

Now, suppose \eqref{eq:olle} holds. Then 
\begin{align*}
\text{diag}\left(
G^*_{1}, G^*_{2}, \ldots,  G^*_{n} \right)V = ([1, 1, \ldots, 1]^T \otimes I_d)Q,
\end{align*}
where $Q$ is some matrix in $GL(d, \mathbb{R})$. This means that 
\begin{align*}
\text{diag}\left(G^*_{1}G^{-1**}_1, G^*_{2}G^{-1**}_2, \ldots,  G^*_{n}G^{-1**}_n \right) = I_n \otimes Q,
\end{align*}
which implies that $\{G^{**}_{i} \}_{i \in \mathcal{V}}$ and $\{G^{*}_{i}\}_{i \in \mathcal{V}}$ are equal up to transformation from the left. By using Lemma~\ref{lem:1} we can conclude that 
$$G^*_{ij} = G^{*-1}_{i}G^*_{j} \text{ for all } (i,j) \in \mathcal{E}.$$
\hfill $\blacksquare$\newline

\begin{rem}\label{rem:1}
In Lemma~\ref{lem:2},
the relation $$\text{im}(\text{diag}\left(
G^*_{1}, G^*_{2}, \ldots,  G^*_{n} \right)V) = \text{ker}(L \otimes I_d)$$
holds if and only if
for any matrix $V_2$, where the columns thereof comprise a basis for $\text{ker}(L \otimes I_d)$, there is a matrix 
$Q$
such that 
$$\text{diag}\left(
G^*_{1}, G^*_{2}, \ldots,  G^*_{n} \right)V = V_2Q.$$
\end{rem}

\begin{rem}\label{rem:2}
In Lemma~\ref{lem:2}, if $\mathcal{G}$ is connected but not 
QSC,  it can hold that $A^T$ is the adjacency matrix of a 
QSC graph $\mathcal{G}' = (\mathcal{V}', \mathcal{E}')$. Then 
it holds that 
$$\text{im}(U_1(\{G^{*-T}_i\}_{i \in \mathcal{V}}) = \text{ker}(Z(\mathcal{G}', \{{G}^{*T}_{ij}\}_{(i,j) \in \mathcal{E}}).$$
\end{rem}

Lemma~\ref{lem:2} is important as it provides a way of finding matrices
$\{{G}^*_{i}\}_{i \in \mathcal{V}}$ fulfilling \eqref{eq:1}. In the 
following subsection, this lemma is used to provide a least squares method.

\subsection{A least squares method}
Suppose the graph $\mathcal{G} = (\mathcal{V}, \mathcal{E})$ is QSC,
and the collection $\{{G}_{ij}\}_{(i,j) \in \mathcal{E}}$ of matrices
in $GL(d, \mathbb{R})$ are not transitively consistent for $\mathcal{G}$, but close to being transitively consistent (closeness 
is in the sense of some matrix norm in $\mathbb{R}^{d \times d}$).
Then,
motivated by Lemma~\ref{lem:2}, the collection $\{G^*_i\}_{i \in \mathcal{V}}$ of matrices in $GL(d, \mathbb{R})$
such that \eqref{eq:1} holds can be found by using the following 
approach. 
\vspace{2mm}
 
\subsection*{Algorithm 1}
\begin{enumerate}
\item 
Solve the problem
$$\underset{{V}}{\min} \|Z(\mathcal{G}, \{{G}_{ij}\}_{(i,j) \in \mathcal{E}})V\|_F^2,$$
where 
$V \in \mathbb{R}^{nd \times d}$, $V^TV = I_d$. This is done by means of the 
Singular Value Decomposition of $Z(\mathcal{G}, \{{G}_{ij}\}_{(i,j) \in \mathcal{E}})$. Let $V_1$ be the optimal solution.\\
 \item Identify the $G_i^*$ in the collection $\{G^*_i\}_{i \in \mathcal{V}}$ by 
$$V_1^T = [G_1^{*-T}, G_2^{*-T}, \ldots, G_n^{*-T}].$$
\end{enumerate}
The algorithm is motivated by Lemma~\ref{lem:2}.
Note that the method is applicable if and only if the graph $\mathcal{G}$ is QSC, (Lemma~\ref{lem:2}). 
In the special case when the transformations are known to be Euclidean (or belong to some other desirable subset of $GL(d, \mathbb{R})$), the collection 
$\{G_i^{**}\}$ can be obtained 
by projecting the {$G_i^{*}$} onto the set of Euclidean transformations (or any other desirable subset of $GL(d, \mathbb{R})$).

If $\{{G}_{ij}\}_{(i,j) \in \mathcal{E}}$ is close
to being transitively consistent, the $d \times d$ block matrices in $V_1$ are 
invertible and can be identified with the $G_1^{*-T}$. This is guaranteed 
by the following lemma \cite{thunberg_synbchronization}.

\begin{lem}\label{lem:11}
In this lemma $Z$ or $Z(\mathcal{G},\{{G}_{ij}\}_{(i,j) \in \mathcal{E}})$ is fixed, whereas 
the matrix $\tilde{Z}$ 
is regarded as a variable in 
$\mathbb{R}^{nd \times nd}$.
Let
\begin{align*}
\mathcal{S}_1 & = \{U \in \mathbb{R}^{nd \times d}: U^TU = I\}, \\
\mathcal{S}_2(\tilde{Z}) & = \arg\min_{U \in \mathcal{S}_1}\textnormal{trace}(U^T\tilde{Z}^T\tilde{Z}U).
\end{align*}
For $\epsilon > 0$, there is $\delta(\epsilon) > 0$ such
that if
$$\|\tilde{Z} - Z \|_{\textnormal{F}} < \delta,$$
it holds that for all $U \in \mathcal{S}_2(\tilde{Z})$,
$$\|U\|_{\mathcal{S}_2({Z})} < \epsilon,$$
where 
$$\|U\|_{\mathcal{S}_2({Z})} = \inf_{V \in \mathcal{S}_2({Z})}\|U - V\|_{\text{F}}.$$
\end{lem}

\subsection{Further results}
Loops in the graph $\mathcal{G}$ are essential for the performance
of Algorithm 1 -- if the graph is QSC and has no
loops, improvement is not possible, see the following lemma.

\begin{lem}\label{lem:3}
If the QSC graph $\mathcal{G}$ is a {spanning tree (containing a center)}, any collection $\{{G}_{ij}\}_{(i,j) \in \mathcal{E}}$ of matrices
in $GL(d, \mathbb{R})$ is transitively consistent for $\mathcal{G}$.
\end{lem}

\noindent\emph{Proof: } \\
\begin{align*}
& Z(\mathcal{G}, \{{G}_{ij}\}_{(i,j) \in \mathcal{E}}) \\
=~& \begin{bmatrix}
I & * & * & \ldots & * & *\\
0 & I & * & \ldots & * & *\\
\vdots & \vdots & \vdots & \ddots & \vdots & \vdots \\
0 & 0 & 0 & \ldots & I & * \\
0 & 0 & 0 & \ldots & 0 & 0
\end{bmatrix},
\end{align*}
where all but one of the $*$ at each (block) row is nonzero and 
an invertible matrix. Due 
to this structure, there is a collection
$\{G^*_i\}_{i \in \mathcal{V}}$ of matrices in $GL(d, \mathbb{R})$ 
such that \eqref{eq:Z_laplacian} holds, which in turn
means that $G_{ij} = G_i^{*-1}G_j$ for all $(i,j) \in \mathcal{E}$. Now since $\mathcal{G}$ is QSC, this means that $G_{ij} = G_i^{*-1}G_j^*$ for all $(i,j) \in \mathcal{V} \times \mathcal{V}$.
\hfill $\blacksquare$\newline

Due to Lemma~\ref{lem:3}, if $\mathcal{G}$ is QSC and a spanning tree and if $\{{G}_{ij}\}_{(i,j) \in \mathcal{E}}$
corresponds to \enquote{disturbed} versions of $\{{G}^*_{ij}\}_{(i,j) \in \mathcal{E}}$,
the solution to Algorithm 1 will only provide the $\{{G}_{ij}\}_{(i,j) \in \mathcal{E}}$
once again.

Lemma \ref{lem:11} provides us with the positive result that the solution to Algorithm 1
depends continuously on the $G_{ij}$ transformations. 
A somewhat negative result is provided by Lemma~\ref{lem:4} below. Unfortunately it is not 
true that \eqref{eq:olle} implies transitive consistency. 

\begin{lem}\label{lem:4}
Let $\mathcal{G} = (\mathcal{V}, \mathcal{E})$ be any QSC graph
satisfying that at least one element in the vector $A(\mathcal{G})[1, 1, \ldots, 1]^T$ is greater
or equal to 2. Let
 $\{{G}^*_{ij}\}_{(i,j) \in \mathcal{E}}$ be a collection of matrices
in $GL(d, \mathbb{R})$, transitively consistent for $\mathcal{G}$. 
Let $\{G^*_i\}_{i \in \mathcal{V}}$ be a collection of matrices in $GL(d, \mathbb{R})$ 
for which it holds that
\begin{equation*}
G^*_{ij} = G^{*-1}_{i}G^*_{j} \text{ for all } (i,j) \in \mathcal{E}.
\end{equation*} 

Now, for any $\epsilon > 0$, there is a collection $\{{G}_{ij}\}_{(i,j) \in \mathcal{E}}$ of matrices in $GL(d, \mathbb{R})$ that are not transitively consistent for $\mathcal{G}$ such that 
\begin{equation}\label{eq:2}
\sum_{(i,j) \in {\mathcal{E}}}\|{G}_{ij} - G^*_{ij}\|_F \leq \epsilon,
\end{equation}
and \eqref{eq:olle} holds for $\{G_{ij}\}_{(i,j) \in \mathcal{E}}$ and a collection $\{G_i\}_{i \in \mathcal{V}}$
of matrices in $GL(d, \mathbb{R})$.
\end{lem}

\noindent\emph{Proof: } \\
Suppose the $k$th element of the vector $A(\mathcal{G})[1, 1, \ldots, 1]^T$ is larger or equal to $2$. Then there is $l,m$ such that
$l \neq k$, $m \neq k$, $G^*_{kl}, G^*_{km} \in GL(d, \mathbb{R})$.
For $0 < \alpha < 1$ let $G_{kl} = (1+\alpha)G^*_{kl}$ and
$G_{km} = (1 - \alpha)G^*_{km}$. Furthermore, let 
$G_{ij} = G^*_{ij}$ for all $(i,j) \not\in \{(k,l),(k,m)\}$.
It is easy to see that the left-hand side of \eqref{eq:2}
is less than or equal to $\alpha(\|G^*_{kl}\|_F + \|G^*_{km}\|_F)$.
Now we choose $$\alpha < \frac{\epsilon}{\|G^*_{kl}\|_F + \|G^*_{km}\|_F}$$ and \eqref{eq:2} is satisfied. By construction, all the $G_{ij}$ are elements of $GL(d, \mathbb{R})$.

Let $G_i = G_i^*$ for all $i$. It holds that 
\begin{align*}
& Z(\mathcal{G}, \{{G}_{ij}\}_{(i,j) \in \mathcal{E}}) \\
\nonumber
=~ &\text{diag}(G_1^{-1}, G_2^{-1}, \ldots, G_n^{-1})((L + Q) \otimes I)\cdot \\
\nonumber
& \text{diag}(G_1, G_2, \ldots, G_n),
\end{align*}
where $Q = [Q_{ij}]$, $Q_{kl} = \alpha$, $Q_{km} = -\alpha$ and $Q_{ij} = 0$ for 
all $(i,j) \not\in \{(k,l),(k,m)\}$. Since $\text{ker}((L + Q) \otimes I)
\supset \text{ker}(L \otimes I)$, \eqref{eq:olle} holds for the $G_i$.
According to Lemma~\ref{lem:2}, if the $G_{ij}$ are transitively consistent and $\mathcal{G}$ is QSC, \eqref{eq:olle} is a condition to guarantee \eqref{eq:1}. But \eqref{eq:1} is not fulfilled 
since $G_kG_{kl}G_l^{-1} = (1 + \alpha)I$. Thus, the $G_{ij}$ are
not transitively consistent.
\hfill $\blacksquare$\newline

After the introduction of Lemma~\ref{lem:4}, one 
might be lead to believe that Algorithm 1 does not 
work well in practice. However, as will be seen in Section~\ref{sec:results}, this is definitely not the case. 

Now, to recap: Transitive consistency is equivalent to \eqref{eq:1}. Lemma~\ref{lem:2} states that when $\mathcal{G}$ is QSC and transitive consistency holds,
\eqref{eq:1} and \eqref{eq:olle} are equivalent. However, Lemma~\ref{lem:4} states that \eqref{eq:olle} is not equivalent to transitive consistency for QSC graphs.

Now we show a stability property of $-Z$. If the $\mathcal{G}_{ij}$
are transitively consistent, it is easy to see (from \eqref{eq:Z_laplacian}) that
$-Z(\mathcal{G}, \{{G}_{ij}\}_{(i,j) \in \mathcal{E}})$ is critically stable, see
definition in Lemma~\ref{lem:10} below. 
However, the following result shows that if the ${G}_{ij}$ transformations
are elements in $O(d)$, i.e., ${G}_{ij}^T{G}_{ij} = I$ for all $i,j$, the matrix
$-Z(\mathcal{G}, \{{G}_{ij}\}_{(i,j) \in \mathcal{E}})$ is critically stable 
regardless if transitive consistency is fulfilled or not.

\begin{lem}\label{lem:10}
For any graph $\mathcal{G} = (\mathcal{V}, \mathcal{E})$ and collection $\{G_{ij}\}_{(i,j) \in \mathcal{E}}$ where $G_{ij} \in O(d)$ for all $(i,j) \in \mathcal{E}$,
the matrix $-Z(\mathcal{G}, \{{G}_{ij}\}_{(i,j) \in \mathcal{E}})$ is critically 
stable, i.e., for any $\epsilon > 0$, there is $\delta(\epsilon)$ such that 
for $x(0) = x_0 \in \mathbb{R}^{nd}$, $\|x(0)\| < \delta$ it holds that 
$$\|x(t)\| < \epsilon,$$
when $$\dot{x}(t) = -Z(\mathcal{G}, \{{G}_{ij}\}_{(i,j) \in \mathcal{E}})x(t).$$
Furthermore, if there are eigenvalues exactly on the imaginary axis,
those eigenvalues are equal to zero.
\end{lem}

\noindent\emph{Proof: } \\
Let 
\begin{equation}\label{eq15}
\dot{x}(t) = -{Z}x(t), \quad x(t) \in \mathbb{R}^{nd},
\end{equation}
where $x(0)$ is the initial state.
We can write $x(t)$ as $x(t) = [x_1^T(t), x_2^T(t), \ldots, x_n^T(t)]^T$, where $x_i(t) \in \mathbb{R}^{d}$ for all $i$. 
Define the function 
$$V(x) = \max_i(x_i^Tx_i).$$
If there is some eigenvalue  of ${Z}$
with negative real part or if there is a Jordan 
block of dimension larger than one  corresponding to an eigenvalue on the imaginary axis, there is $x_0$ such that for the state $x(t)$ with initial
state $x_0$, $V(x(t)) \rightarrow \infty$ as $t \rightarrow \infty$. We want to show that this is not possible. Let us first define the set
$$\mathcal{I}_{\max}(t) = \{i:V(x(t)) = x_i^T(t)x_i(t)\}.$$
Now,
\begin{align}
\label{eq30}
D^+(V(x(t))) & = \max_{i \in \mathcal{I}_{\max}(t)}\frac{d}{dt}x_i^T(t)x_i(t) \\
\nonumber
& =  \max_{i \in \mathcal{I}_{\max}(t)}x_i^T(t)\left(\sum_{j \in \mathcal{N}_i}({G}_{ij}x_j(t) - x_i(t)\right ) \\
\nonumber
& \leq 0,
\end{align}
where $D^+$ is the upper Dini-derivative. A proof 
of the first equality~\eqref{eq30} can be found in~\cite{thunberg2014consensus}
using the results in \cite{yoshizawa1966stability} and \cite{clarke1975generalized}. The result appears frequently in the 
literature~\cite{shi2009,lin2007state}.  Now we can use the 
Comparison Lemma~\cite{khalil2002nonlinear} to show that $V(x(t))$ is decreasing independently 
of the choice of $x_0$. The inequality in~\eqref{eq30}
is a consequence of the fact that the ${G}_{ij}$
are orthogonal matrices.

Now we show that there are no non-zero eigenvalues on the 
imaginary axis. Suppose there are non-zero eigenvalues on the 
imaginary axis, then there must be a nontrivial periodic solution
$\bar{x}(t) = [\bar{x}_1^T, \bar{x}_1^T, \ldots, \bar{x}_n^T]^T$ to  \eqref{eq15}, i.e., $\bar{x}(t)$ is periodic and $\bar{x}(t_1) \neq \bar{x}(t_2)$ for some $t_1 \neq t_2$. It can be shown that $D^+(V(\bar{x}(t))) = 0$ for all $t$ and it can also be shown that a necessary condition for this to hold is that $\bar{x}_i(t) = \bar{x}_i(t)$ for all $t$ and $G_{ij} = I$ for all $(i,j)$.
The procedure to show the latter is a bit intricate and is based on an induction argument
hinging on the fact that $\mathcal{G}$ is QSC. Now, if the $G_{ij} \neq I$, the 
necessary condition is not fulfilled, hence we have a contradiction. In the case when the $G_{ij} = I$ it holds that $Z(\mathcal{G}) = L(\mathcal{G}) \otimes I_d$ and the latter matrix 
does not have any non-zero eigenvalues on the imaginary axis. 
\hfill $\blacksquare$\newline

\subsection{Optimization problems and the $H$-matrix}\label{sec:hMatrix}
In this subsection a matrix $H$ is defined as the Hessian of a quadratic convex function. In the previous subsection the approach was to define
a set of linear constraints, which are fulfilled for 
transitively consistent transformations, and 
then use these constraints to formulate a least squares optimization problem.
In this section the approach is different. Optimization 
problems are formulated directly, without taking a detour 
via algebraic constraints. An assumption throughout this section
is that $\mathcal{G}$ is connected. 

Given the graph $\mathcal{G} = (\mathcal{V}, \mathcal{E})$ and the collection
$\{{G}_{ij}\}_{(i,j) \in \mathcal{E}}$ of matrices in $GL(d, \mathbb{R})$, we formulate three optimization problems, where the the first, (P1), corresponds to the 
exact problem we want to solve. The objective function is non-convex and the 
constraint set is non-compact. {The second problem (P2) is a restriction
of the first problem having a compact constraint set (with a non-convex objective function).
In contrast, the third problem (P3) has
a quadratic convex objective function of 
the $G_i^{-1}$ as well as a compact constraint set.}

\begin{align*}
& \text{(P1)} \quad 
\begin{cases}
\min\limits_{\{G_{i}\}_{i \in \mathcal{V}}}\sum\limits_{(i,j) \in \mathcal{E}}\frac{1}{2}\|G_{ij} - G^{-1}_iG_j \|_F^2, \\
 \quad \text{s.t.}~\quad \: \: G_i \in GL(d,\mathbb{R}).
\end{cases} \\
& \\
& \text{(P2)} \quad 
\begin{cases}
\min\limits_{\{G_{i}\}_{i \in \mathcal{V}}}\sum\limits_{(i,j) \in \mathcal{E}}\frac{1}{2}\|G_{ij} - G^{-1}_iG_j \|_F^2, \\
 \quad \text{s.t.}~\quad \: \: G_i \in O(d). 
 \end{cases} \\
 & \\
& \text{(P3)} \quad 
\begin{cases}
\min\limits_{\{G_{i}\}_{i \in \mathcal{V}}}\sum\limits_{(i,j) \in \mathcal{E}}\frac{1}{2}\|G_{ij}G_j^{-1} - G_i^{-1} \|_F^2, \\
 \quad \text{s.t.}~\quad \: \: U_1(\{G_i\}_{i \in \mathcal{V}})^TU_1(\{G_i\}_{i \in \mathcal{V}}) = Q \succ 0.
\end{cases} 
\end{align*}
Define the two functions
\begin{align}
f(\{G^{-1}_i\}_{i \in \mathcal{V}}) & = \sum\limits_{(i,j) \in \mathcal{E}}\frac{1}{2}\|G_{ij}G_j^{-1} - G_i^{-1} \|_F^2, \nonumber\\
g(\{G_i\}_{i \in \mathcal{V}}) & = \sum\limits_{(i,j) \in \mathcal{E}}\frac{1}{2}\|G_{ij} - G^{-1}_iG_j \|_F^2. \label{gFcn}
\end{align}
The matrix $Q$ is symmetric and positive definite.
We implicitly assume that $g$ and $f$ are parameterized by $\{G_{ij}\}_{(i,j) \in \mathcal{E}}$. There is a 
similar problem to (P3), defined by left-multiplication by the $G_i$ instead 
of right-multiplication by the $G_i^{-1}$:
\begin{align*}
& \text{(P4)} \quad 
\begin{cases}
\min\limits_{\{G_{i}\}_{i \in \mathcal{V}}}\sum\limits_{(i,j) \in \mathcal{E}}\frac{1}{2}\|G_iG_{ij} - G_j \|_F^2, \\
 \quad \text{s.t.}~\quad \: \: U_2(\{G_i\}_{i \in \mathcal{V}})U_2(\{G_i\}_{i \in \mathcal{V}})^T = Q \succ 0,
\end{cases} 
\end{align*}
The two problems are equivalent. We choose to study (P3) instead of (P4) in order to more easily see the connection 
between the Hessian (the $H$-matrix) in the problem and the matrix $Z$ (cf. Section~\ref{sec:problemP2}).

The problem (P2) and variants thereof has received attention lately \cite{tron2014distributed}. 
Exact solutions do not exist in general and local gradient descent methods 
are used. One of the more important contributions of this work is that we 
provide a lower bound for the global solution of this problem as well as a method
for which the bound is almost tight in numerical experiments.

\subsection{Problem (P3) and its connection to problem (P1) -- definition of the $H$-matrix}\label{sec:problemP2}
Let $X = U_1(\{G_i\}_{i \in \mathcal{V}})$. By a slight
abuse of notation, let 
$$f(X) = f(\{G^{-1}_i\}_{i \in \mathcal{V}}) = \sum\limits_{(i,j) \in \mathcal{E}} {\frac{1}{2}}\|G_{ij}^TG_i^{-1} - G_j^{-1} \|_F^2.$$
Now $$\nabla f(X) = X^TH({\mathcal{G}},\{{G}_{ij}\}_{(i,j) \in \mathcal{E}}),$$
where
\begin{align*}
&~  H({\mathcal{G}},\{{G}_{ij}\}_{(i,j) \in \mathcal{E}}) \\
= &~ Z({\mathcal{G}}, \{{G}_{ij}\}_{(i,j) \in \mathcal{E}}) + Z_2(\bar{\mathcal{G}}, \{\bar{{G}}_{ij}\}_{(i,j) \in \bar{\mathcal{E}}});
\end{align*}
and
\begin{align*}
&~ {Z_2(\bar{\mathcal{G}}, \{\bar{{G}}_{ij}\}_{(i,j) \in \bar{\mathcal{E}}})} \\
= ~& \text{diag}(W(\bar{\mathcal{G}}, \{\bar{G}_{ij}\}_{(i,j) \in \bar{\mathcal{E}}})W(\bar{\mathcal{G}}, \{\bar{G}_{ij}\}_{(i,j) \in \bar{\mathcal{E}}})^T) \\
&  - W(\bar{\mathcal{G}}, \{\bar{G}_{ij}\}_{(i,j) \in {\bar{\mathcal{E}}}});
\end{align*}
{$\bar{G}_{ij} = G^T_{ji}$ for all $i,j$ and $\bar{\mathcal{G}} = (\mathcal{V}, \bar{\mathcal{E}})$ is the graph constructed reversing the direction of the edges in $\mathcal{G}$. The operator $\text{diag}(\cdot)$ is here understood in the block-matrix sense, i.e., for a matrix $B \in \mathbb{R}^{nd\times nd}$, $\text{diag}(B) = (I_n \otimes 1_d 1_d^T ) \odot B$, where $\odot$ denotes element-wise multiplication, $I_n$ is the $n$-dimensional identity matrix and $1_d$ is the $d$-dimensional vector containing ones.}

\begin{rem}
A more general formulation of the objective 
functions $f$ and $g$ with positive weights 
is
\begin{align*}
\tilde{f}(\{G^{-1}_i\}_{i \in \mathcal{V}}) & = \sum\limits_{(i,j) \in \mathcal{E}}\frac{1}{2}a_{ij}\|G_{ij}G_j^{-1} - G_i^{-1} \|_F^2, \\
\tilde{g}(\{G_i\}_{i \in \mathcal{V}}) & = \sum\limits_{(i,j) \in \mathcal{E}}\frac{1}{2}a_{ij}\|G_{ij} - G^{-1}_iG_j \|_F^2.
\end{align*}
The $a_{ij}$ are positive for all $i,j$. This way of defining the objective functions lead to a slight modification
of the $H$-matrix. Equivalent results to all the results obtained in
this section for the 
$H$-matrix can also be formulated for this alternative definition
of the $H$-matrix with weights. The alternative $H$-matrix can also be used in an equivalent distributed
algorithm to the one presented in Section~\ref{sec:5.2}.
\end{rem}

\begin{lem}\label{lem:5}
In the special case when all the $G_{ij}$ are elements of $O(d)$, i.e.,
orthogonal matrices, 
$$Z_2(\bar{\mathcal{G}}, \{\bar{{G}}_{ij}\}_{(i,j) \in \bar{\mathcal{E}}}) = Z(\bar{\mathcal{G}}, \{\bar{{G}}_{ij}\}_{(i,j) \in \bar{\mathcal{E}}}).$$
Furthermore, if the graph $\mathcal{G}$ is also symmetric, 
$$Z_2(\bar{\mathcal{G}}, \{\bar{{G}}_{ij}\}_{(i,j) \in \bar{\mathcal{E}}}) = 
Z({\mathcal{G}}, \{{{G}}_{ij}\}_{(i,j) \in {\mathcal{E}}})^T.$$
\end{lem}

\begin{lem}\label{lem:5}
For any connected graph $\mathcal{G} = (\mathcal{V}, \mathcal{E})$,
and collection $\{{G}^*_{ij}\}_{(i,j) \in \mathcal{E}}$ of matrices
in $GL(d, \mathbb{R})$, the collection $\{{G}^*_{ij}\}_{(i,j) \in \mathcal{E}}$ is transitively consistent
for $\mathcal{G}$
if and only if there is a collection $\{{G}^*_{i}\}_{i \in \mathcal{V}}$ of 
matrices in $GL(d, \mathbb{R})$ such that 
$$\text{im}(U_1(\{G^*_i\}_{i \in \mathcal{V}})) \subset \text{ker}(H(\mathcal{G}, \{{G}^*_{ij}\}_{(i,j) \in \mathcal{E}})).$$

The collection $\{{G}^*_{i}\}_{i \in \mathcal{V}}$ satisfies \eqref{eq:1}.
\end{lem}

\noindent\emph{Proof: }
Suppose $\{{G}_{ij}\}_{(i,j) \in \mathcal{E}}$ is transitively consistent, then, according to Lemma~\ref{lem:1}, there is $\{{G}^*_{i}\}_{i \in \mathcal{V}}$ such that \eqref{eq:1} holds, which in turn can be used to show that 
\begin{align}
\nonumber
&f(\{G^{*-1}_i\}_{i \in \mathcal{V}}) \\ 
\nonumber
=~ &U_1(\{G^*_i\}_{i \in \mathcal{V}})^TH(\mathcal{G}, \{{G}^*_{ij}\}_{(i,j) \in \mathcal{E}})U_1(\{G_i\}_{i \in \mathcal{V}}) \\
\label{eq:3}
=~& 0.
\end{align}
On the other hand, if $\{{G}^*_{ij}\}_{(i,j) \in \mathcal{E}}$ is not transitively 
consistent, there is no $\{{G}^*_{i}\}_{i \in \mathcal{V}}$ such that \eqref{eq:1} holds. It can now be shown that \eqref{eq:3} does not hold for any collection $\{{G}^*_{i}\}_{i \in \mathcal{V}}$ of matrices in $GL(d, \mathbb{R})$.
\hfill $\blacksquare$\newline

\begin{lem}\label{lem:6}
For any connected graph $\mathcal{G} = (\mathcal{V}, \mathcal{E})$
and  collection $\{{G}^*_{ij}\}_{(i,j) \in \mathcal{E}}$ of matrices
in $GL(d, \mathbb{R})$ -- transitively consistent for $\mathcal{G}$ --
it holds that 
$$\text{dim}(\text{ker}(H(\mathcal{G}, \{{G}^*_{ij}\}_{(i,j) \in \mathcal{E}}))) = d.$$
\end{lem}

\noindent\emph{Proof: }
Due to Lemma~\ref{lem:5}, we know that 
\begin{equation}\label{eq:4}
\text{dim}(\text{ker}(H(\mathcal{G}, \{{G}^*_{ij}\}_{(i,j) \in \mathcal{E}}))) \geq d.\end{equation}
Thus, we need to show that the inequality in \eqref{eq:4} cannot be strict.
Since $\{{G}^*_{ij}\}_{(i,j) \in \mathcal{E}}$ is transitively consistent, there is $\{{G}^*_{i}\}_{i \in \mathcal{V}}$ fulfilling \eqref{eq:1}.

Suppose the inequality is strict for $\{{G}^*_{ij}\}_{(i,j) \in \mathcal{E}}$. 
We know there is $\{G_i^*\}_{i \in \mathcal{V}}$ where
$G_i^* \in GL(d, \mathbb{R})$ for all $i$, such
that 
$$\text{im}(U_1(\{G^*_i\}_{i \in \mathcal{V}})) \subset \text{ker}(H(\mathcal{G}, \{{G}^*_{ij}\}_{(i,j) \in \mathcal{E}})).$$
Now there must be a vector $y = [y_1^T, y_2^T, \ldots, y_n^T]^T \in \mathbb{R}^{nd}$, where the $y_i$ are in $\mathbb{R}^d$, such that 
 $$y \in \text{ker}(H(\mathcal{G}, \{{G}^*_{ij}\}_{(i,j) \in \mathcal{E}})),$$ $y \neq 0$, and $y^TU_1(\{G^*_i\}_{i \in \mathcal{V}}) = 0$. 
 There must be $k$ and $l$ such that the $l$th element of $y_k$ is nonzero. 
 The set of transformations $\{G^{*-1}_kG^*_i\}_{i \in \mathcal{V}}$ satisfy \eqref{eq:1} (Lemma \ref{lem:1}) and $f(\{(G^{*-1}_kG^*_i)^{-1}\}_{i \in \mathcal{V}}) = 0$. Now, 
 let  $$\bar{X} = [\bar{x}_1, \bar{x}_2, \ldots, \bar{x}_d] = U_1(\{G^{*-1}_kG^*_i\}_{i \in \mathcal{V}}),$$
 and 
  $$\bar{Y} = [\bar{x}_1, \bar{x}_2, \ldots, \bar{x}_{l-1}, y, \bar{x}_{l+1}, \bar{x}_d] = U_1(\{G^{*-1}_kG^*_i\}_{i \in \mathcal{V}}).$$
We know that 
$H({\mathcal{G}},\{{G}_{ij}\}_{(i,j) \in \mathcal{E}})\bar{Y} = 0$. For all $i$, let $\bar{Y}_i$
be the $i$th $d \times d$ block matrix in $\bar{Y}$. We know by construction that
$\bar{Y}_{k} \in GL(d, \mathbb{R})$. Now, for any $j \in \mathcal{N}_k$
it holds that 
$$\|G^*_{kj}\bar{Y}_j - \bar{Y}_k\|_F = 0,$$
which implies that $\bar{Y}_j \in GL(d, \mathbb{R})$.
Also, for any $i$ such that $k \in \mathcal{N}_i$, it holds that 
$$\|G^*_{ik}\bar{Y}_k - \bar{Y}_i\|_F = 0,$$
which implies that $\bar{Y}_i \in GL(d, \mathbb{R})$. Now, due to the fact that $\mathcal{G}$ is connected, an induction argument
can be used to show that all the $\bar{Y}_i$ are elements in $GL(d, \mathbb{R})$.

The collection $\{\bar{Y}_i\}_{i \in \mathcal{V}}$ satisfies 
$$G^*_{ij} = \bar{Y}_{i}\bar{Y}_{j}^{-1} \text{ for all } (i,j) \in \mathcal{E}.$$
It is easy to see that the two collections $\{\bar{Y}^{-1}_i\}_{i \in \mathcal{V}}$ 
and $\{\bar{G}^*_i\}_{i \in \mathcal{V}}$ are not equal up to transformation from the left. But, since the graph is connected, the two must be equal up to 
transformation from the left (Lemma \ref{lem:1}). This is a contradiction.  
Hence it is a false assumption that the inequality in \eqref{eq:4} is strict.
\hfill $\blacksquare$\newline

We summarize the results of Lemma~\ref{lem:5} and Lemma~\ref{lem:6}
in the following proposition.

\begin{prop}\label{prop:1}
The collection $\{{G}^*_{ij}\}_{(i,j) \in \mathcal{E}}$ is transitively consistent
for the connected graph $\mathcal{G}$ if and only if there is a collection $\{{G}^*_{i}\}_{i \in \mathcal{V}}$ of 
matrices in $GL(d, \mathbb{R})$ such that 
$$\text{im}(U_1(\{G^*_i\}_{i \in \mathcal{V}})) = \text{ker}(H(\mathcal{G}, \{{G}^*_{ij}\}_{(i,j) \in \mathcal{E}})).$$
The $G_i^*$ satisfy \eqref{eq:1}.
\end{prop}

\noindent\emph{Proof: }
Direct application of Lemma~\ref{lem:5} and Lemma~\ref{lem:6}.
\hfill $\blacksquare$\newline

The following proposition provides a similar, but somewhat stronger result.

\begin{prop}\label{prop:2}
The collection $\{{G}^*_{ij}\}_{(i,j) \in \mathcal{E}}$ of matrices in $GL(d, \mathbb{R})$ is transitively consistent
for the connected graph $\mathcal{G}$ if and only if 
$${\text{dim}(\text{ker}(H(\mathcal{G}, \{{G}^*_{ij}\}_{(i,j) \in \mathcal{E}})))} = d.$$
\end{prop}

\noindent\emph{Proof: } \\
$\textbf{If:}$
Let $\bar{Y} = [y_1, y_2, \ldots, y_{nd}]^T \in \mathbb{R}^{nd \times d}$ be any full rank matrix such that 
$$H(\mathcal{G}, \{{G}^*_{ij}\}_{(i,j) \in \mathcal{E}})\bar{Y} = 0.$$
All the $y_i \in \mathbb{R}^{d}$.
Let $\bar{Y}_i$
be the $i$th $d \times d$ block matrix in $\bar{Y}$. 
Since $\bar{Y}$ is full rank, there is a (finite) sequence  
$\{y_{i_j}\}_{j =1}^d$ such that $[y_{i_1}, y_{i_2}, \ldots
y_{i_d}] \in GL(d, \mathbb{R})$. 

Now, for $k \in \mathcal{V}$ we know that for any $j \in \mathcal{N}_k$
it holds that 
$$\|G^*_{kj}\bar{Y}_j - \bar{Y}_k\|_F = 0,$$
which implies that $\text{im}(\bar{Y}_j^T) = \text{im}(\bar{Y}_k^T)$.
Also, for any $i$ such that $k \in \mathcal{N}_i$, it holds that 
$$\|G^*_{ik}\bar{Y}_k - \bar{Y}_i\|_F = 0,$$
which implies that $\text{im}(\bar{Y}_i^T) = \text{im}(\bar{Y}_k^T)$. Now, due to the fact that $\mathcal{G}$ is connected, an induction argument
can be used to show that $\text{im}(\bar{Y}_j^T) = 
\text{im}(\bar{Y}_i^T)$ for all $i,j$. But then 
$$\text{im}([y_{i_1}, y_{i_2}, \ldots
y_{i_d}]) \subset \text{im}(\bar{Y}_j^T) \text{ for all }j,$$
which together with the fact that $[y_{i_1}, y_{i_2}, \ldots
y_{i_d}] \in GL(d, \mathbb{R})$ is full rank, can be used to show that $\bar{Y}_i \in GL(d, \mathbb{R})$ for all $i$.
Thus, 
$$\text{im}(U_1(\{\bar{Y}^{-1}_i\}_{i \in \mathcal{V}})) = \text{ker}(H(\mathcal{G}, \{{G}^*_{ij}\}_{(i,j) \in \mathcal{E}})),$$
and the desired result follows from Proposition~\ref{prop:1}
where the $G_i^*$ are replaced by the $\bar{Y}_i^{-1}$.

\vspace{2mm}
\noindent \textbf{Only if:}
Direct application of Lemma~\ref{lem:6}.
\hfill $\blacksquare$\newline

\begin{lem}\label{lem:8}
The optimal solution to (P3) is
$$X^* = VP,$$
where $P \in \mathbb{R}^{d \times d}$, and $V \in \mathbb{R}^{nd\times d}$.
The matrix $V$ is given by the solution to the 
problem
\begin{align*}
& \text{(P5)} \quad 
\begin{cases}
~\min\limits_W \text{trace}(W^TH({\mathcal{G}},\{{G}_{ij}\}_{(i,j) \in \mathcal{E}})W), \\
~W \in \mathbb{R}^{nd \times d}, W^TW = I.
 \end{cases}
\end{align*}
and {$P$} is given by the solution to the problem
\begin{align*}
& \text{(P6)} \quad 
\begin{cases}
~\min\limits_{\tilde{P}} \text{trace}(\tilde{P}^T(V^TH({\mathcal{G}},\{{G}_{ij}\}_{(i,j) \in \mathcal{E}})V)\tilde{P}), \\
~\tilde{P} \in \mathbb{R}^{d \times d}, \tilde{P}^T\tilde{P} = Q.
 \end{cases}
\end{align*}
\end{lem}

\noindent\emph{Proof: }
In the new notation, problem (P3) is written as 
\begin{align*}
& \quad 
\begin{cases}
~\min\limits_{W,\tilde{P}} \text{trace}(\tilde{P}^T(W^TH({\mathcal{G}},\{{G}_{ij}\}_{(i,j) \in \mathcal{E}})W)\tilde{P}), \\
~\tilde{P} \in \mathbb{R}^{d \times d}, \tilde{P}^T\tilde{P} = Q, W \in \mathbb{R}^{nd \times d}, W^TW = I.
 \end{cases}
\end{align*}
In the following derivations it is assumed that $\tilde{P}$ and 
$W$ belong to the constraint sets defined in the problem above.
\begin{align*}
& \min\limits_{W,\tilde{P}}\text{trace}(\tilde{P}^T(W^TH({\mathcal{G}},\{{G}_{ij}\}_{(i,j) \in \mathcal{E}})W)\tilde{P}) \\
=~&  \min\limits_{W,\tilde{P}}\text{trace}(\tilde{P}^T(Q^T(W)D(W)Q(W))\tilde{P}) \\
=~&  \min\limits_{W,\tilde{P}}\text{trace}(\tilde{P}^TD(W)\tilde{P}) \\
=~& \min\limits_{\tilde{P}}\text{trace}(\tilde{P}^T(V^TH({\mathcal{G}},\{{G}_{ij}\}_{(i,j) \in \mathcal{E}})V)\tilde{P}),
\end{align*}
where $Q^T(W)D(W)Q(W)$ is the spectral factorization of $$W^TH({\mathcal{G}},\{{G}_{ij}\}_{(i,j) \in \mathcal{E}})W.$$
\hfill $\blacksquare$\newline

\begin{prop}\label{lem:9}
For any $Q_1 \succ 0$ and $Q_2 \succ 0$ let $\{G^*_i\}_{i \in \mathcal{V}}$ and $\{G^{**}_i\}_{i \in \mathcal{V}}$ be the transformations obtained from the optimal solutions of problem (P3) with $Q$ equal
to $Q_1$ and $Q$ equal to $Q_2$ respectively. It holds that 
$$g(\{G^*_i\}_{i \in \mathcal{V}}) = g(\{G^{**}_i\}_{i \in \mathcal{V}}),$$
i.e., the value of $g$ is independent of $Q$.
\end{prop}

\noindent\emph{Proof: }
According to Lemma~\ref{lem:8} the transformations are equal up to 
transformation from the left.
\hfill $\blacksquare$\newline

\begin{rem}
It is implicitly assumed in Proposition~\ref{lem:9}
that the $G_i^*$ and the $G_i^{**}$ are in $GL(d, \mathbb{R})$. This is guaranteed if the $G_{ij}$ are 
sufficiently close to be transitively consistent. 
The result to guarantee this is omitted but analogous to the statement in Lemma~\ref{lem:11} for the $Z$-matrix.
\end{rem}

\begin{lem}\label{lem:20}
For any graph $\mathcal{G} = (\mathcal{V}, \mathcal{E})$ and collection $\{G_{ij}\}_{(i,j) \in \mathcal{E}}$ where $G_{ij} \in GL(d,\mathbb{R})$ for all $(i,j) \in \mathcal{E}$,
the matrix $-H(\mathcal{G}, \{{G}^*_{ij}\}_{(i,j) \in \mathcal{E}})$ is critically 
stable, i.e., for any $\epsilon > 0$, there is $\delta(\epsilon)$ such that 
for any $x(0) = x_0 \in \mathbb{R}^{nd}$, $\|x(0)\| < \delta$ it holds that 
$$\|x(t)\| < \epsilon,$$
when $$\dot{x}(t) = -H(\mathcal{G}, \{{G}^*_{ij}\}_{(i,j) \in \mathcal{E}})x(t).$$
\end{lem}

\noindent\emph{Proof: } \\
The matrix $H(\mathcal{G}, \{{G}^*_{ij}\}_{(i,j) \in \mathcal{E}})$ is the
Hessian matrix and hence positive semi-definite.
\hfill $\blacksquare$\newline

\subsection{A least squares method}
Proposition~\ref{lem:9} is important, it states that we can without loss of
generality assume that $Q = I$, since the choice of $Q$ does not 
affect the value of $g$, i.e., the cost function we want to minimize. The value of $f$ changes with $Q$, but this is of less importance. Motivated by these results
we introduce a least squares method along the lines of Algorithm 1.

\subsection*{Algorithm 2}
\begin{enumerate}
\item 
Let $V_1$ be the optimal solution to problem (P5).\\
 \item Identify the $G_i^*$ in the collection $\{G^*_i\}_{i \in \mathcal{V}}$ by 
$$V_1^T = [G_1^{*-T}, G_2^{*-T}, \ldots, G_n^{*-T}].$$
\end{enumerate}
The algorithm is motivated by Proposition~\ref{prop:1} and
Proposition~\ref{lem:9}. If the collection $\{{G}_{ij}\}_{(i,j) \in \mathcal{E}}$
is close enough to be transitively consistent, step 2) can be executed,
i.e., each $d \times d$ sub-block of the matrix $V_1$ is invertible. 
The result that guarantees this is analogous to the statement in Lemma~\ref{lem:11}.

\subsection{A Gauss-Newton method}\label{sec:gn}
In this section a Gauss-Newton method is presented.
The solution obtained in Algorithm 2 is used as the
initialization for the algorithm.

The Fr\'echet derivatives of the identity map and 
the inverse map at the point $G_i$ in the direction 
$E_i$ are given by 
\begin{align*}
L_{\text{id}}(G_i,E_i) & = E_i, \\
L_{\text{inv}}(G_i,E_i) & = -G_i^{-1}E_iG_i^{-1},
\end{align*}
respectively. Higham~\cite{higham2008functions}
provides a good 
introduction to Fr\'echet derivatives for matrix 
functions. Let $\{E_i\}_{i \in \mathcal{V}}$ be a
collection of matrices in $\mathbb{R}^{n \times n}$. It holds that
\begin{align}
\nonumber
& g(\{G_i+ E_i\}_{i \in \mathcal{V}}) \\
\nonumber
 = &  \sum\limits_{(i,j) \in \mathcal{E}}\frac{1}{2}\|G_{ij} - G^{-1}_iG_j 
  -  G_i^{-1}L_{\text{id}}(G_j,E_j) \\
\nonumber
  & - L_{\text{inv}}(G_i,E_i)G_j  + o(\|[E_i,E_j]\|_{{F}}^2) \|_F^2.
\end{align}
Let 
\begin{align*}
& \bar{g}(\{G_i\}_{i \in \mathcal{V}},\{E_i\}_{i \in \mathcal{V}}) \\
 = & \sum\limits_{(i,j) \in \mathcal{E}}\frac{1}{2}\|G_{ij} - G^{-1}_iG_j 
  -  G_i^{-1}L_{\text{id}}(G_j,E_j) \\
\nonumber
  & - L_{\text{inv}}(G_i,E_i)G_j  \|_F^2 \\
 = &  \sum\limits_{(i,j) \in \mathcal{E}}\frac{1}{2}\|G_{ij} - G^{-1}_iG_j 
  -  G_i^{-1}E_j \\
\nonumber
  & + G_i^{-1}E_iG_i^{-1}G_j \|_F^2.
\end{align*}
Consider the following problem 
\begin{align*}
& \text{(P7)} \quad 
\begin{cases}
~\min\limits_{\{E_i\}_{i \in \mathcal{V}}} \bar{g}(\{G_i\}_{i \in \mathcal{V}},\{E_i\}_{i \in \mathcal{V}}).
 \end{cases}
\end{align*}
Problem (P7) is solved in each Gauss-Newton step 
of the method we present below (Algorithm 3). Its solution is given
by the collection $\{E_i^{*}\}_{i \in \mathcal{V}}$
obtained by
\begin{equation}\label{eq:gn1}
\text{vec}(U_2(\{{E}^{*}_i\}_{i \in \mathcal{V}})) = x,
\end{equation}
where $x$ is obtained by the solution of
\begin{align}
\nonumber
& H_{\text{GN}}(\{G_i\}_{i \in \mathcal{V}},\{G_{ij}\}_{(i,j) \in \mathcal{E}})x \\
\label{eq:gn2}
=~ & -c_{\text{GN}}(\{G_i\}_{i \in \mathcal{V}},,\{G_{ij}\}_{(i,j) \in \mathcal{E}});
\end{align}
$\text{vec}(\cdot)$ is the vectorization operator,
i.e., it returns a vector with the stacked columns (in consecutive order) of
its matrix-argument. The matrix $H_{\text{GN}} \in \mathbb{R}^{nd^2 \times nd^2}$ and the vector 
$c_{\text{GN}} \in \mathbb{R}^{nd^2}$ are defined as
follows (for simplicity we have omitted the explicit dependence
of $\{G_i\}_{i \in \mathcal{V}}$ and $\{G_{ij}\}_{(i,j) \in \mathcal{E}}$):
$$H_{\text{GN}} = [\bar{H}_{ij}],$$
where $\bar{H}_{ij} \in \mathbb{R}^{d^2 \times d^2}$ for all $i,j$.
When $i \neq j$, $\bar{H}_{ij}$ is defined by
\begin{align*}
& \bar{H}_{ij} = \\ 
& \begin{cases}
0
\end{cases} & \text{ if } 
\begin{cases}
(j \not\in \mathcal{N}_i), \\
(i \not\in \mathcal{N}_j),\\
\end{cases} \\ \vspace{2mm}
& \begin{cases}
-((G_i^{-1}G_j)\otimes(G_i^{-T}G_i^{-1}))
\end{cases} & \text{ if } 
\begin{cases}
(j \in \mathcal{N}_i), \\
(i \not\in \mathcal{N}_j),\\
\end{cases} \\ \vspace{2mm}
& \begin{cases}
-((G_i^TG_j^{-T})\otimes(G_j^{-T}G_j^{-1}))
\end{cases} & \text{ if } 
\begin{cases}
(j \not\in \mathcal{N}_i), \\
(i \in \mathcal{N}_j),\\
\end{cases} \\ \vspace{2mm}
& \begin{cases}
-((G_i^{-1}G_j)\otimes(G_i^{-T}G_i^{-1}))   \\
-~((G_i^TG_j^{-T})\otimes(G_j^{-T}G_j^{-1}))
\end{cases} & \text{ if } 
\begin{cases}
(j \in \mathcal{N}_i), \\
(i \in \mathcal{N}_j). \\
\end{cases} \\ \vspace{2mm}
\end{align*}
When $i = j$, $\bar{H}_{ii}$ is defined by 
\begin{align*}
\bar{H}_{ii}  = &  \sum_{j \in \mathcal{N}_i}((G_i^{-1}G_jG_j^TG_i^{-T}) \otimes (G_i^{-T}G_i^{-1})) + \\
& \sum_{\{j: i \in \mathcal{N}_j\}}I_d\otimes(G_j^{-T}G_j^{-1}).
\end{align*}
Now, $c_{\text{GN}} = [c_1^T, c_2^T, \ldots, c_n^T]^T$, where 
$c_i \in \mathbb{R}^{d^2}$ for all $i$. The $c_i$ are defined
by
\begin{align*}
c_i = & \sum_{j \in \mathcal{N}_i}((G_i^{-1}G_j)\otimes G_i^{-T})\text{vec}(G_{ij} - G_i^{-1}G_j) - \\
& \sum_{\{j: i \in \mathcal{N}_j\}}I_d\otimes G_j^{-T}\text{vec}(G_{ji} - G_j^{-1}G_i).
\end{align*}

\subsection*{Algorithm 3}
\begin{enumerate}
\item Run Algorithm 2 and let $\{G^*_i\}_{i \in \mathcal{V}}$
bet the collection of matrices obtained in step (2) of that 
algorithm. \\
\item Let $\mathbb{R}^{d \times d} \ni E_i^* = 0$ for $i = 1,2,\ldots, n$. \\
 
\item \textbf{repeat:} \\
\begin{enumerate}
\item $G_i \rightarrow G_i + E_i^* \text{ for all } i,$ \\
\item Update the $E_i^*$ by \eqref{eq:gn1} and \eqref{eq:gn2},
i.e., $$\text{vec}(U_2(\{{E}_i^*\}_{i \in \mathcal{V}})) = x,$$
where $x$ is the solution to \eqref{eq:gn2}.
\end{enumerate}
\end{enumerate}
The stoping criteria in step (3) of Algorithm 5 could be that 
the improvement of the cost function is smaller than a certain 
threshold for two consecutive iterations, or it could be
that a certain number of iterations have been executed etc.
It should be noted that $H_{\text{GN}}$ is both positive definite and 
sparse. In order to solve \eqref{eq:gn2} one can use for example the Conjugate Gradient method~\cite{luenberger1973introduction,forsgren2015connection}.

\subsection{Affine and Euclidean transformations}
In this subsection we consider affine and Euclidean 
transformations. These transformations are linear 
when homogenous coordinates are used. To be more precise,
an element in $\text{Aff}(d,\mathbb{R})$ is a matrix 
$$G = \begin{bmatrix}
Q & {t} \\
0 & 1
\end{bmatrix},$$
where $Q \in GL(d,\mathbb{R})$, ${t} \in \mathbb{R}^d$ and $1$ is a 
scalar. Its inverse is given by
$$G^{-1} = \begin{bmatrix}
Q^{-1} & -Q^{-1}{t} \\
0 & 1
\end{bmatrix}.$$ 
Euclidean transformations, $E(d)$, is a special case of affine
transformations where the matrix $Q \in \mathcal{O}(d)$.

For any connected graph $\mathcal{G} = (\mathcal{V}, \mathcal{E})$  (due to Lemma~\ref{lem:1}), if and only if the 
collection $\{{G}^*_{ij}\}_{(i,j) \in \mathcal{E}}$
is transitively consistent and contains only 
affine transformations, there is a unique (up
to transformation from the left by affine transformations) collection
$\{G_i\}_{i \in \mathcal{V}}$ of affine 
transformations such that 
$$G^*_{ij} = G_i^{-1}G_j.$$
Each $G_i$ is given by 
$$G = \begin{bmatrix}
Q_i & {t}_i \\
0 & 1
\end{bmatrix},$$
and each $G^*_{ij}$ is given by 
$$G^*_{ij} = \begin{bmatrix}
Q_i^{-1}Q_j & Q_i^{-1}({t}_j - {t}_i) \\
0 & 1
\end{bmatrix}.$$

Now, let $\{{G}^*_{ij}\}_{(i,j) \in \mathcal{E}}$ 
be a collection of matrices in $\text{Aff}(d,\mathbb{R})$ that are
not necessarily transitively consistent. It holds that (by a slight abuse of notation)
\begin{align}
\nonumber
g(\{G_i\}_{i \in \mathcal{V}})  = ~& \sum\limits_{(i,j) \in \mathcal{E}}\frac{1}{2}\|G_{ij} - G^{-1}_iG_j \|_F^2 \\
\nonumber
 = ~& g(\{Q_i\}_{i \in \mathcal{V}})  \\ 
\label{eq:31}
 & ~+ \sum\limits_{(i,j) \in \mathcal{E}}\frac{1}{2}\|{t}_{ij} - Q_i^{-1}({t}_j - {t}_i) \|_F^2,
\end{align}
where ${t}_{ij}$ is the translational
part of the transformation $G_{ij}$. 
We see that there is a special structure of \eqref{eq:31}, where the cost function consists of two parts. The first part is only a function of the $Q_i$, whereas the second part is a function of both rotations and translations. 

Define the following optimization problem 
\begin{align*}
& \text{(P8)} \quad 
\begin{cases}
~\min\limits_{\{{t}_i\}_{i \in \mathcal{V}}}\sum\limits_{(i,j) \in \mathcal{E}}\frac{1}{2}\|{t}_{ij} - Q_i^{-1}({t}_j - {t}_i) \|_F^2.
 \end{cases}
\end{align*}
Let 
\begin{align*}
c_{\text{Aff}} & = [c_1^T, c_2^T, \ldots, c_n^T]^T,
\end{align*}
where 
\begin{align*}
c_i & = \sum_{j \in \mathcal{N}_i}Q_i^{-T}{t}_{ij} - \sum_{\{j: i \in \mathcal{N}_j\}}Q_j^{-T}{t}_{ji}.
\end{align*}

Let 
$$H_{\text{Aff}} = [\tilde{H}_{ij}],$$
where $\tilde{H}_{ij} \in \mathbb{R}^{(d-1) \times (d-1)}$ for all $i,j$.
When $i \neq j$, $\tilde{H}_{ij}$ is defined by
\begin{align*}
& \tilde{H}_{ij} = \\ 
& \begin{cases}
0
\end{cases} & \text{ if } 
\begin{cases}
(j \not\in \mathcal{N}_i), \\
(i \not\in \mathcal{N}_j),\\
\end{cases} \\ \vspace{2mm}
& \begin{cases}
-Q_i^{-T}Q_i^{-1}
\end{cases} & \text{ if } 
\begin{cases}
(j \in \mathcal{N}_i), \\
(i \not\in \mathcal{N}_j),\\
\end{cases} \\ \vspace{2mm}
& \begin{cases}
-Q_j^{-T}Q_j	^{-1}
\end{cases} & \text{ if } 
\begin{cases}
(j \not\in \mathcal{N}_i), \\
(i \in \mathcal{N}_j),\\
\end{cases} \\ \vspace{2mm}
& \begin{cases}
-Q_i^{-T}Q_i	^{-1}   \\
-~Q_j^{-T}Q_j	^{-1}
\end{cases} & \text{ if } 
\begin{cases}
(j \in \mathcal{N}_i), \\
(i \in \mathcal{N}_j). \\
\end{cases} \\ \vspace{2mm}
\end{align*}
When $i = j$, $\tilde{H}_{ii}$ is defined by 
\begin{align*}
\tilde{H}_{ii}  = &  \sum_{j \in \mathcal{N}_i}Q_i^{-T}Q_i^{-1} +
\sum_{\{j: i \in \mathcal{N}_j\}}Q_j^{-T}Q_j	^{-1}.
\end{align*}
The matrix $H_{\text{Aff}}$ and the vector $c_{\text{Aff}}$ depend on 
 $\{G_i\}_{i \in \mathcal{V}}$ and $\{G_{ij}\}_{(i,j) \in \mathcal{E}}$.

The solutions to the problem (P8) is given 
by the elements in the set
\begin{align*}
\{\{{t}_i\}_{i \in \mathcal{V}}: H_{\text{Aff}}[{t}_1^T, {t}_2^T, \ldots, {t}_n^T ]^T
& = -c_{\text{Aff}}\}.
\end{align*}

The Gauss-Newton method developed in Section~\ref{sec:gn}, i.e., Algorithm 3,
can be adapted to the case of affine transformations.
Now we require that 
$$E_i \odot B = E_i \text{ for all } i, \text{ where } B = \begin{bmatrix}
1_d1_d^T & 1_d \\
0 & 0
\end{bmatrix}.$$
We remind the reader that $\odot$ denotes element-wise multiplication.
In each iteration (in the modified step (3) of Algorithm 3)
the collection $\{E_i^{*}\}_{i \in \mathcal{V}}$
is obtained by
\begin{equation}\label{eq:gn3}
\text{vec}(U_2(\{{E}^{*}_i\}_{i \in \mathcal{V}})) = x,
\end{equation}
where $x = Xv$, $v$ is obtained by the solution to
\begin{align}
\label{eq:gn4}
& X^TH_{\text{GN}}Xv = -X^Tc_{\text{GN}},
\end{align}
and $X \in \mathbb{R}^{nd^2 \times n(d-1)d}$ is defined below.
$$X = I_n \otimes (I_d \otimes \bar{B}),$$
where $$\bar{B} = \begin{bmatrix}
I_{d-1} \\
0 
\end{bmatrix} \in \mathbb{R}^{d \times (d-1)}.$$

Now we present the following algorithm for 
affine transformations.
\subsection*{Algorithm 4}
\begin{enumerate}
\item 
Run Algorithm 2 for the collection 
$\{Q_{ij}\}_{(i,j) \in \mathcal{E}}$ and let $\{Q_i\}_{i \in \mathcal{V}}$ be the 
matrices obtained in step (2) of that 
algorithm.\\
 \item Solve problem (P8) for the $t_i$ using the $Q_i$
 from (1) and let $$G_i = \begin{bmatrix}
 Q_i & {t}_i \\
 0 & 1
 \end{bmatrix} \text{ for all } i.$$ \\
 \item Let $\mathbb{R}^{d \times d} \ni E_i^* = 0$ for $i = 1,2,\ldots, n$. \\
 \item \textbf{while} a stoping criteria has not been met: \\
\begin{enumerate}
\item $G_i \rightarrow G_i + E_i^* \text{ for all } i,$ \\
\item Update the $E_i^*$ by \eqref{eq:gn3} and \eqref{eq:gn4},
i.e., $$\text{vec}(U_2(\{{E}_i^*\}_{i \in \mathcal{V}})) = x,$$
where $x$ is the solution to \eqref{eq:gn3}.
\end{enumerate}
\end{enumerate}

\begin{rem}
There are many variations of Algorithm 4 that
can be employed. The most simple one is to omit steps (3) and (4). Another one is to run the Gauss-Newton 
method (Algorithm 3) for the $Q_i$ matrices after 
step (2). The expression in \eqref{eq:31} can also
be changed to include weights. For example, if the orthogonal matrices are closer to be transitively consistent than the translations, the first part of the expression, i.e., $g(\{Q_i\}_{i \in \mathcal{V}})$, could be weighted with a positive weight larger than $1$.
\end{rem}

\begin{rem}
After a slight modification, Algorithm 4 can be used
for Euclidean transformations instead of affine ones. 
In this case Algorithm 5 (see Section~\ref{sec:44}) is used in (1)  to generate the $Q_i$ transformations instead of Algorithm 2. Numerical simulations (see Section~\ref{sec:results}) show
that this is a good method in comparison to Algorithm 1 or Algorithm 2 (where the matrices are finally projected onto the set of Euclidean transformations $E(d)$).
\end{rem}

\section{Orthogonal matrices}\label{sec:44}
In this section problem (P2) is studied. For orthogonal matrices the 
objective functions $f$ and $g$ are equivalent. The Gauss-Newton method (Algorithm 3) is hence not necessary.
Furthermore, the orthogonal matrices is an important class of matrices,
not the least in dimension $d = 3$.

We begin by formulating the following result.

\begin{prop}\label{prop:3}
For the connected graph $\mathcal{G} = (\mathcal{V}, \mathcal{E})$, let $\{{G}_{ij}\}_{(i,j) \in \mathcal{E}}$ be a collection of 
matrices in $GL(d, \mathbb{R})$. 
Let $\{G_i\}_{i \in \mathcal{V}}$ be a collection of matrices obtained from Algorithm 2. Let $\{G^*_i\}_{i \in \mathcal{V}}$ be a 
collection of matrices solving the optimization problem (P2).
It holds that 
$$f(\{\sqrt{n}G^{-1}_i\}_{i \in \mathcal{V}}) \leq g(\{G^*_i\}_{i \in \mathcal{V}}).$$
\end{prop}

\noindent\emph{Proof: } \\
It is easy to verify that for orthogonal matrices,
(P3) is a relaxation of (P2) when $Q = nI$. Now, (Proposition~\ref{lem:9}) the solution to (P3) with $Q = nI$ is provided 
by the matrices obtained by Algorithm 2 after scaling 
by $\frac{1}{\sqrt{n}}$. 
\hfill $\blacksquare$\newline

Let us now extend Algorithm 1 (Algorithm 2) in the following way.

\subsection*{Algorithm 5}
\begin{enumerate}
\item 
Same as in Algorithm 1 (same as in Algorithm 2).\\
 \item Same as in Algorithm 1 (same as in Algorithm 2).\\
\item Let $G_i^{**}$ be the projection of $G^*_j$ onto $O(d)$, 
i.e., 
$$G_i^{**} = QV^T,$$
where $QDV^T$ is the SVD of $G_i^{*}$. Let
$$G^*_{ij} =  G_i^{T**}G_j^{**}.$$
The collection $\{G^*_{ij}\}_{(i,j) \in \mathcal{E}}$ is the final 
transitively consistent collection.
\end{enumerate}

Proposition~\ref{prop:3} can now be used to provide performance 
guarantees. An upper bound on the closeness to optimality
 is given by
 \begin{equation}\label{eq:5}
 g(\{G^{**}_i\}_{i \in \mathcal{V}})) - f(\sqrt{n}G^{*-1}_i\}_{i \in \mathcal{V}}),
 \end{equation}
where the $G^{*}_i$ are obtained from Algorithm 2 and the $G^{**}_i$
are obtained from Algorithm 5 -- assuming the first two steps are the same as in 
Algorithm 2. 

If the ${G}_{ij}$ are also elements in $O(d)$, the difference in \eqref{eq:5}
is almost tight. For example, in the case when $d = 3$, $n = 100$, and the $G_{ij}$ are
generated from $G_i$-matrices and $R_{ij}$-matrices matrices by
$G_{ij} = G_i^{-1}G_jR_{ij}$  ($R_{ij}$ is an orthogonal matrix with geodesic
distance to $I$ less or equal to $\pi/4$. It is generated by drawing a skew symmetric matrix from the uniform distribution over the closed ball with radius $\pi/4$ and then taking the matrix exponential of that matrix).
Let
\begin{align}
\nonumber
 	& h(\{G^{**}_i\}_{i \in \mathcal{V}}, \{G^{*-1}_i\}_{i \in \mathcal{V}}) \\
 	= & \frac{g(\{G^{**}_i\}_{i \in \mathcal{V}})) - f(\{\sqrt{n}G^{*-1}_i\}_{i \in \mathcal{V}})}{f(\{\sqrt{n}G^{*-1}_i\}_{i \in \mathcal{V}})}. \label{gapFcn}\\
 	\intertext{For $1000$ experiments we observe that}
 	& h(\{G^{**}_i\}_{i \in \mathcal{V}}, \{G^{*-1}_i\}_{i \in \mathcal{V}}) \leq 6*10^{-4}. \nonumber
 \end{align} 
This means that the solution obtained by Algorithm 5 is closer
than $0.06 \%$ to the global optimum of problem (P2). The graphs in these experiments were QSC and the adjacency matrices contained $100$ zero entries.

\subsection{Distributed algorithms}\label{sec:5}
In this subsection we show that Algorithm 5 
can be implemented in a distributed way. 
Besides the graph $\mathcal{G} = (\mathcal{V}, \mathcal{E})$, which describes
what transformations are available, another graph $\mathcal{G}^{\text{com}} = (\mathcal{V}, \mathcal{E}^{\text{com}})$
is used. It is always assumed $\mathcal{E} \subset \mathcal{E}^{\text{com}}$. The graph $\mathcal{G}^{\text{com}}$ is referred to 
as \emph{the communication graph}. The assumptions on the communication graph $\mathcal{G}^{\text{com}}$ differ between the two presented algorithms.

\subsubsection{Orthogonal matrices and QSC communication graph}\label{sec:dist:orthogonal}
Here it is
assumed that all transformations are orthogonal matrices, i.e., elements in $O(d)$. That 
is, the ${G}_{ij}$ matrices as well as the ${G}^*_{ij}$  matrices and the ${G}^*_{i}$ matrices are assumed 
to be elements in $O(d)$.

The algorithm will now be presented,
after which an explanation and
justification is provided. In this
algorithm it is assumed that $\mathcal{G} = \mathcal{G}^{\text{com}}$
is QSC. 
The notation $\mathcal{N}_i$ is used to denote $\mathcal{N}_i(\mathcal{G}) = \mathcal{N}_i(\mathcal{G}^{\text{com}})$.
\vspace{2mm}
 
\subsection*{Algorithm 6}
Let $$X(0) = [X_1^T(0), X_2^T(0), \ldots, X_n^T(0)]^T,$$
where, for all $i$, {the elements of the matrix $X_i(0) \in \mathbb{R}^{d \times d}$ are drawn from  
$\mathcal{U}(-0.5,0.5)$, i.e., the uniform distribution with the open interval $(-0.5,0.5)$ as support.}
Let $X_i(t)$ for $t \in \mathbb{N}$ be defined by the following 
distributed algorithm:
\begin{align*}
{X}_1(t+1) & = {X}_1(t) + \epsilon\sum_{j \in \mathcal{N}_1}({G}_{1j}X_j(t) - X_1(t)), \\
{X}_2(t+1) & = {X}_2(t) + \epsilon\sum_{j \in \mathcal{N}_2}({G}_{2j}X_j(t) - X_2(t)), \\
& \hspace{2mm} \vdots \\
{X}_n(t+1) & = {X}_n(t) + \epsilon\sum_{j \in \mathcal{N}_n}({G}_{nj}X_j(t) - X_n(t)),
\end{align*}
where $\epsilon > 0$.1
In compact notation this is written as 
\begin{equation}\label{algorithm:2}
{X}(t+1) = {X}(t) - \epsilon Z(\mathcal{G}, \{{G}_{ij}\}_{(i,j) \in \mathcal{E}}))X(t).
\end{equation}
For a sufficiently large $t$, let 
$G_{i}^{T*}$
be the projection of $X_i(t)$ onto $O(d)$,
and let $G_{ij}^* = G_{i}^{T*}G_{j}^{*}$ for all $i,j$.
It should be noted that if the spectral radius 
is not known, in practice it is enough to choose $\epsilon$ to something small.

\subsubsection*{Analysis of the algorithm}
In this section the theoretical analysis 
of the algorithm is provided. 
The first thing we need to guarantee is that the 
matrix 
$$I - \epsilon Z(\mathcal{G},\{{G}_{ij}\}_{(i,j) \in \mathcal{E}}),$$
appearing in the right-hand side of the discrete time linear
system \eqref{algorithm:2}, is critically 
stable in the linear dynamical systems sense. This means that all eigenvalues 
must be smaller than or equal to $1$ in absolute value and any Jordan-block 
corresponding to an eigenvalue whose
absolute value is $1$
must be one-dimensional \cite{lindquist1996introduction}. 

\begin{lem}
If $\mathcal{G} (\mathcal{V}, \mathcal{E})$ is QSC and $\epsilon > 0$ small enough
it holds that 
$$I - \epsilon Z(\mathcal{G},\{{G}_{ij}\}_{(i,j) \in \mathcal{E}})$$
is critically stable.
\end{lem}

\noindent\emph{Proof: }
According to Lemma~\ref{lem:10} it holds that $Z$ is critically stable 
and has no non-zero eigenvalues on the imaginary axis. This means that 
for $\epsilon$ small enough the eigenvalues of $\epsilon Z(\mathcal{G},\{{G}_{ij}\}_{(i,j) \in \mathcal{E}})$ are located in the closed unit disc centered at $-1$; 
the eigenvalues on the boundary are simple. 
\hfill $\blacksquare$\newline

\begin{rem}
Numerical simulations seem to indicate that in practice one can choose 
$$\epsilon \in \left (0, \frac{1}{\rho(Z(\mathcal{G}, \{{G}_{ij}\}_{(i,j) \in \mathcal{E}}))}\right ),$$
where $\rho(Z(\mathcal{G}, \{{G}_{ij}\}_{(i,j) \in \mathcal{E}}))$ is the spectral radius.
\end{rem}

Now we can deduce that if $\epsilon$ is chosen small enough, $X(t)$ 
converges (to something).  It easy to verify (Lemma~\ref{lem:2}) that
if the $G_{ij}$ were transitively consistent,
$X(t)$ {would converge with} exponential 
rate of convergence to 
$$\bar{X} = [\bar{x}_1, \bar{x}_2, \ldots, \bar{x}_d] = [\bar{X}_1^T, \bar{X}_2^T, \ldots, \bar{X}_n^T]^T,$$
where
$\bar{x}_i$ is the projection of the $i$th
column of $X(0)$ onto $\text{ker}(Z)$ and 
$\bar{X}_i \in \mathbb{R}^{d \times d}$ for
all $i$.
Since the $X_i(0)$ are drawn from the distribution 
{$(\mathcal{U}(-0.5,0.5))^{d \times d}$}, it is extremely unlikely that (probability zero)
$\bar{X}$ has not full rank.  If all the $\bar{X}_i$ 
are full rank matrices, 
$$G_{ij} = \bar{X}_i\bar{X}_j^T \quad \text{ for all } i,j.$$

Now, if the $G_{ij}$ are not transitively consistent, in general the $X_{i}(t)$ converge
to $0$, which is not favorable. However,
if the $G_{ij}$ are close to being transitively 
consistent, since the eigenvalues
of $Z$ are continuous in the $G_{ij}$, the 
$d$ smallest eigenvalues are significantly 
smaller in magnitude than the other eigenvalues;
also the $d$ smallest singular values are
significantly smaller than the other singular 
values. Up to rotation, the right-singular vectors corresponding 
to the $d$ smallest 
singular values are continuous in the 
$G_{ij}$, see Lemma~\ref{lem:11}.

Let the right-singular vectors corresponding 
to the $d$ smallest 
 singular values comprise the 
columns of the matrix $Y \in \mathbb{R}^{nd \times d}$.
The matrix $Y$ is equal 
to $V$ obtained in the first step of Algorithm 1
(up to transformation from the left). 
Now, as $t \rightarrow \infty$, under the assumption 
that the $G_{ij}$ are close to the $G^*_{ij}$, the columns of $X(t)$ 
converge to $\text{im}(Y)$ much faster than 
$X(t)$ converges to $0$. Thus, for $t$ large enough
$X(t)$ is approximately equal to $Y$ up 
to transformation from the left. This convergence
can be seen in Figure~\ref{distributedZOrthogonal}
for different choices of $n$, $d$, and magnitudes of 
noise.

The last step of the algorithm is 
 justified by Lemma~\ref{lem:2}.

\subsubsection{Orthogonal matrices and 
symmetric connected communication graph}\label{sec:5.2}
In this section a general distributed algorithm is 
presented, which works for $G_{ij}$ matrices in $GL(d, \mathbb{R})$, 
a directed connected  graph $\mathcal{G}$, and a symmetric communication graph $\mathcal{G}^{\text{com}}$. We will make the assumption 
that $\mathcal{G}^{\text{com}}$ is the union graph of 
$\mathcal{G}$ and $\bar{\mathcal{G}}$, i.e.,
$\mathcal{G}^{\text{com}} = (\mathcal{V}, \mathcal{E} \cup \bar{\mathcal{E}})$. The difference between Algorithm 6 and Algorithm 7 presented here, is that the $Z$-matrix
is used in the former, whereas the 
$H$-matrix is used in the latter. Here, different from Section~\ref{sec:dist:orthogonal}, it does not
hold that $\mathcal{N}_i(\mathcal{G}^{\text{com}}) = \mathcal{N}_i(\mathcal{G})$ for all $i$. When 
we write $\mathcal{N}_i$ this is shorthand for 
$\mathcal{N}_i(\mathcal{G}^{\text{com}})$.

\subsection*{Algorithm 7}
Let $$X(0) = [X_1^T(0), X_2^T(0), \ldots, X_n^T(0)]^T,$$
where
{the elements of the matrix $X_i(0) \in \mathbb{R}^{d \times d}$ are drawn from  
$\mathcal{U}(-0.5,0.5)$.}
Let $X_i(t)$ for $t \in \mathbb{N}$ be defined by the following 
distributed algorithm:
\begin{align*}
{X}_1(t+1) & = {X}_1(t) + \epsilon\sum_{j \in \mathcal{N}_1}Q_{1j}X_j(t) - V_{1j}X_1(t)), \\
{X}_2(t+1) & = {X}_2(t) + \epsilon\sum_{j \in \mathcal{N}_2}(Q_{2j}X_j(t) - V_{2j}X_2(t)), \\
& \hspace{2mm} \vdots \\
{X}_n(t+1) & = {X}_n(t) + \epsilon\sum_{j \in \mathcal{N}_n}(Q_{nj}X_j(t) - V_{nj}X_n(t)),
\end{align*}
where $$\epsilon \in \left (0, \frac{1}{\rho(H(\mathcal{G}, \{{{G}_{ij}}\}_{(i,j) \in \mathcal{E}}))} \right );$$ 
$\rho(H(\mathcal{G}, \{{{G}_{ij}}\}_{(i,j) \in \mathcal{E}}))$ is the spectral radius and
\begin{align*}
Q_{ij} & ~= G_{ij} + G_{ji}^T, \\
V_{ij} & ~=I_d + G_{ji}^TG_{ji}.
\end{align*}
In compact notation this is written as 
\begin{equation}\label{algorithm:222}
{X}(t+1) = {X}(t) - \epsilon H(\mathcal{G}, \{{{G}_{ij}}\}_{(i,j) \in \mathcal{E}}))X(t).
\end{equation}
For a sufficiently large $t$, let 
$G_{i}^{*{-1}}= X_i(t)$
and let $G_{ij}^* = G_{i}^{*{-1}}G_{j}^{*}$ for all $i,j$.
It should be noted that if the spectral radius 
is not known, in practice it is enough to choose $\epsilon$ to something small.

\begin{rem}
In the definitions of $Q_{ij}$ and $V_{ij}$,
in the case when $(j,i) \not\in \mathcal{G}$, 
the symbol $G_{ij}$ should be interpreted as 
the matrix in $\mathbb{R}^{d\times d}$ containing only 
zero-elements.
\end{rem}

\subsubsection*{Analysis of the algorithm}
The convergence analysis of this Algorithm is analogous 
simpler than that of Algorithm 6. The eigenvalues of
the $H$-matri are real since the matrix is symmetric. Instead of using Lemma~\ref{lem:10}, Lemma~\ref{lem:20} is used instead. 

\subsection{Gradient flow for orthogonal matrices}\label{sec:6}
Under the assumption that all the $G_{ij}$
are elements in $O(d)$, we here provide 
a method,
which will be used for comparison to our earlier 
methods. Results, along the lines of the ones presented in this
section, can be found in \cite{afsari2011riemannian,afsari2004some,sarlette2009autonomous}.

For all $i$, define the cost functions
$$g_i: (O(d))^{n} \rightarrow \mathbb{R}^+$$
by
\begin{align*}
g_i(G_1, G_2, \ldots, G_{n}) = 
\sum_{j \in \mathcal{N}_i}A_{ij}\|G_i{G}_{ij}G_j^T - I\|_F^2.
\end{align*}
The overall cost function  
$$g: (O(d))^{n} \rightarrow \mathbb{R}^+$$
is equal to $g$, i.e.,
\begin{align*}
g(\{G^*_i\}_{i \in \mathcal{V}}) = 
\sum_{i = 1}^{n}g_i(G_1, G_2, \ldots, G_{n}).
\end{align*}

The (negative) gradient flow on $(O(d))^n$ of $g$ is given
by 
\begin{align}
\label{eq:7}
\dot{G}_i(t) & ~= -
\sum_{j \in \mathcal{N}_i}A_{ij}((G_i(t){G}_{ij}G_j(t)^T)^T \\
\nonumber
&~ \: \: \: \: \: - (G_i(t){G}_{ij}G_j(t)^T))G_i(t) \\
\nonumber
&~ \: \: \: \: \:  -\sum_{\{j: i \in \mathcal{N}_j\}}A_{ij}((G_j(t){G}_{ij}^TG_i(t)^T)^T \\
\nonumber
&~ \: \: \: \: \: - (G_j(t){G}_{ij}^TG_j(t)^T))G_i(t), \\
\nonumber
&~\text{for all } i \in \mathcal{V}.
\end{align}

Now we present an algorithm, which improves on 
Algorithm 5. However, as will be seen in Section~\ref{sec:results},
this improvement is marginal. 

\subsection*{Algorithm 8}
\begin{enumerate}
\item 
Run Algorithm 5 and let $\{G_i^{**}\}$ be the 
orthogonal matrices obtained in step (3) of the 
algorithm.\\
 \item Solve \eqref{eq:7} numerically (for example by 
 using ode45 in Matlab) for a sufficiently large time  interval $[0,T]$ with the $G_i^{**}$ as 
 initial conditions.\\
 \item Let $\{G_i(t)^{-1}G_j(t)\}_{(i,j) \in \mathcal{V} \times \mathcal{V}}$ be the collection 
 of transitively consistent matrices.
\end{enumerate}

\begin{rem}
In step (3) of Algorithm 8, if the $G_i(T)$ are 
not elements of $O(d)$ (due to errors from numerical integration), they need to be projected onto $O(d)$.
\end{rem}

\section{Numerical verification}\label{sec:results}

In our experiments, we consider Algorithm 8 first. Subsequently, the centralised Z- and H-matrix methods are evaluated for different configurations. Eventually, the analogous distributed methods are used in our simulations. In order to compare the methods,
an assumption throughout this section is that the 
graph $\mathcal{G}$ -- describing what transformations are available -- is QSC.

\newcommand{\nDraws}{100}
\newcommand{\nRepetitions}{1}

\subsection{Generating graphs and transformations}\label{sec:72}
For each of the following experiments, the collection $\{G^*_{i}\}_{i \in \mathcal{V}}^n$ is generated by drawing random matrices in $O(d)$ \cite{bernard_cvpr}. From that, the (full) set of transitively consistent matrices $\{G^*_{ij} = G_{i}^{*-1} G^*_{j}\}_{i,j \in \mathcal{V}}$ is created. The noisy set of pairwise transformations $\{G_{ij}  \}_{i,j \in \mathcal{V}}$ is generated by adding element-wise Gaussian noise with zero mean and standard deviation $\sigma$ to each $G^*_{ij}$. After adding the element-wise Gaussian noise, the matrix is additionally projected onto $O(d)$.

Furthermore, a quasi-strongly connected (QSC) graph with graph density $\rho$ -- not mix up with the spectral radius of a matrix -- is generated in the 
following manner. For generating a minimum QSC graph $\mathcal{G} = (\mathcal{V}, \mathcal{E})$,
two lists are used. One list $\mathcal{L}^{\mathcal{G}}$ keeps track of the nodes that are already considered, and one list $\mathcal{L}^{\mathcal{\bar G}}$ keeps track of the nodes that have not been considered. By a \emph{minimum} QSC graph we mean a QSC graph that is a (spanning) tree, i.e., one with exactly $n-1$ edges. Initially, we set $\mathcal{V} = \{1,\ldots,n\}$, ${\mathcal{E} = \{(i,i) ~ : ~ i \in \mathcal{V}\}}$, $\mathcal{L}^{\mathcal{G}} = \{r\}$, where $r \in \mathcal{V}$ is a randomly selected node, and $\mathcal{L}^{\mathcal{\bar G}} = \{1,\ldots,n\} - \{r\}$. Then, the following procedure is repeated $n{-}1$ times: pick random nodes $i \in \mathcal{L}^{\mathcal{\bar G}}$ and $j \in \mathcal{L}^{\mathcal{G}}$, add the edge $(i,j)$ to $\mathcal{E}$, and update $\mathcal{L}^{\mathcal{\bar G}}$ and $\mathcal{L}^{\mathcal{G}}$ accordingly. After $n{-}1$ repetitions a (minimum) QSC graph has been generated. At this point we store the edge set $\mathcal{E}$ and call it $\mathcal{E}^{\mathit{QSC}}$. Next, random edges are added to $\mathcal{E}$ until the the density of the graph is larger than or equal to $\rho$, which is defined below.

We remind the reader that $A$ is the adjacency matrix of $\mathcal{G}$ with elements $A_{ij}$. The graph density $\rho(\mathcal{G})$ is defined by
\begin{align}\label{eq21}
  \rho(\mathcal{G}) = \left( \frac{1}{n^2-\vert \mathcal{E}^{\mathit{QSC}} \vert}\sum_{(i,j) \notin \mathcal{E}^{\mathit{QSC}}}^{n} A_{ij} \right) \enspace .
\end{align} 
The intuition behind the graph density is that it is the proportion of the number of present edges in $\mathcal{G}$ with respect to a fully connected graph (having $n^2$ edges) excluding the edges in $\mathcal{E}^{\mathit{QSC}}$. With that, $\rho=0$ denotes a minimum QSC graph, whereas $\rho=1$ denotes a fully connected graph. Generating random QSC graphs with different values of the parameter $\rho$ allows us to consider different degrees of missing transformations. 

Using the graph $\mathcal{G}$ with density $\rho$, the collection $\{G_{ij}\}_{(i,j) \in \mathcal{E}}$ is the one that is eventually used for the evaluation. In the simulations, for each individual sub-figure the simulations have been performed with $\nDraws$ random sets of orthogonal transformations (the transitively consistent ones and the synthetically generated noisy versions thereof) and QSC graphs. 
Shown in the sub-figures is the mean of all runs.

\subsection{Algorithm 8 -- orthogonal matrices}
\begin{figure}[!th] 
\centering
\includegraphics[scale=0.37]{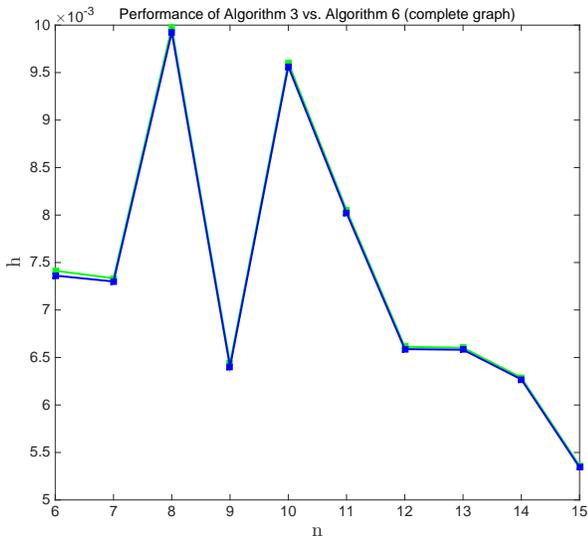}
\caption{Upper bounds on the optimality gap, i.e., $h$, for the solution of Algorithm 5, green line, and Algorithm 8, blue line. The graph is complete.}
 \label{fig_1}
\end{figure}

For $d = 3$, Figure~\ref{fig_1} shows upper bounds on the gap between the optimal value and the value of the objective function obtained by two methods -- Algorithm 5, green curve, and Algorithm 8, blue curve. In Algorithm 8, the initial states are given 
by the solution to Algorithm 5. The ODE in \eqref{eq:7} is solved numerically in Matlab by ode45. For each number of 
coordinate systems $n$,
100 simulations are conducted and averages are shown in
Figure~\ref{fig_1}. In each simulation a set of transitively consistent orthogonal matrices $\{G^*_{ij}\}_{(i,j) \in \mathcal{V} \times \mathcal{V}}$ are generated from a set
of orthogonal matrices $\{G^*_{i}\}_{i \in \mathcal{V}}$ 
according to the description in Section~\ref{sec:72} below. The graph $G$
used in each of the experiments is the complete graph.   

\begin{figure}[!th] 
\centering
\includegraphics[scale=0.37]{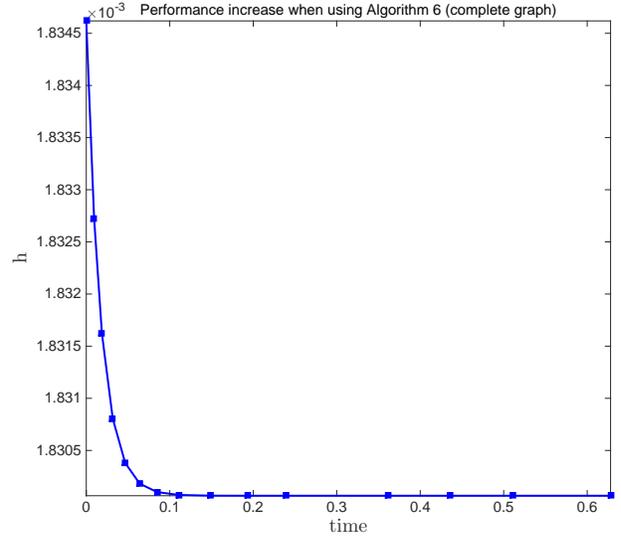}
\caption{Improvement of $h$ when Algorithm 8 is used. In this case when $n = 15$, $d =3$, and the graph is complete.}
 \label{fig_2}
\end{figure}

For a single numerical experiment, Figure~\ref{fig_2} shows the improvement of  $h$ when Algorithm 8 is used. One can see that Algorithm 5 generates matrices 
that are close to a local optimum of problem (P2).

It can be seen that only a marginal improvement can be made using the significantly more computationally expensive Algorithm 8. Due to the heavy computational burden, in the following simulations we omit Algorithm 8 and focus on the methods based on the Z- and H-matrix.

\subsection{Centralized methods for matrices in $O(d)$}
 In this set of experiments we compare the H-matrix method, the Z-matrix method and a (naive) reference-based method, where the latter serves as baseline for the comparison.

\subsubsection{The reference-based method}
For the reference-based method, a minimum QSC graph 
$\mathcal{G}^{\text{min-QSC}} = (\mathcal{V}, \mathcal{E}^{\text{min-QSC}} \subseteq \mathcal{E})$ 
with $n-1$ edges is randomly drawn as a subgraph of 
$\mathcal{G} = (\mathcal{V}, \mathcal{E})$. 
For that, all centers (see Def.~\ref{qscGraph}) of the graph $\mathcal{G}$ are initially determined by looking at the $n \times n$ distance matrix between all $n$ nodes. From the set of centers, a node $c$ is randomly selected. Since $\mathcal{G}$ is QSC, there is at least one such center. Let $\mathcal{L}^{\mathcal{G}^{\text{min-QSC}}} = \{c\}$ be the list of nodes that have already been considered and initialise $\mathcal{E}^{\text{min-QSC}} = \emptyset$. The following procedure is repeated until $\vert \mathcal{E}^{\text{min-QSC}} \vert = n-1$: randomly select a node $r \in \mathcal{L}^{\mathcal{G}^{\text{min-QSC}}}$, select a random node $r' \in  \{i : (i,r) \in \mathcal{E}\}$, if there is such an $r'$, add the edge $(r',r)$ to $\mathcal{E}^{\text{min-QSC}}$ and add $r'$ to $\mathcal{L}^{\mathcal{G}^{\text{min-QSC}}}$.

Per construction, the graph $\mathcal{G}^{\text{min-QSC}}$ is a spanning tree that contains a center. Thus, according to Lemma~\ref{lem:3}, the set $\{G_{ij}\}_{(i,j) \in \mathcal{E}^{\text{min-QSC}}}$ is transitively consistent for $\mathcal{G}^{\text{min-QSC}}$. W.l.o.g., by setting $G_c^* = I$ for the center $c$ of $\mathcal{G}^{\text{min-QSC}}$, all (other) $G_i^*$ are (uniquely) determined as 
\begin{align}
	G_i^* = G_j^* G_{ij}^{-1} \quad \text{for } i \neq c, ~(i,j) \in \mathcal{E}^{\text{min-QSC}} \,. \label{refGi}
\end{align}

To summarise, in the reference-based method a (random) rooted spanning tree graph is considered as subgraph of $\mathcal{G}$, i.e., all but $n-1$ relative transformations $G_{ij}$ (accounting for the transitive inconsistency) are discarded such that the remaining $n-1$ relative transformations are transitively consistent.

 \begin{figure*}[!t]
  \centerline{
    \subfigure{%
    % This file was created by matlab2tikz v0.4.7 running on MATLAB 8.4.
% Copyright (c) 2008--2014, Nico Schlömer <nico.schloemer@gmail.com>
% All rights reserved.
% Minimal pgfplots version: 1.3
% 
% The latest updates can be retrieved from
%   http://www.mathworks.com/matlabcentral/fileexchange/22022-matlab2tikz
% where you can also make suggestions and rate matlab2tikz.
% 
\begin{tikzpicture}

\begin{axis}[%
width=7cm,
height=3cm,
scale only axis,
separate axis lines,
every outer x axis line/.append style={white!15!black},
every x tick label/.append style={font=\color{white!15!black}},
xmin=0,
xmax=0.7,
xlabel={$\sigma$},
every outer y axis line/.append style={white!15!black},
every y tick label/.append style={font=\color{white!15!black}},
ymin=-0.1,
ymax=5.1,
ylabel={g'},
title style={font=\bfseries},
title={$n{=}30$, $d{=}3$, $\rho{=}0.5$}
]
\addplot [color=blue,solid,line width=1.0pt,forget plot]
  table[row sep=crcr]{%
0	5.80681667455282e-30\\
0.1	0.0284437038279075\\
0.2	0.112902714165976\\
0.3	0.279976280361374\\
0.4	0.627971167735205\\
0.5	1.12142279891524\\
0.6	1.62138611466591\\
0.7	2.05583829186495\\
};
\addplot [color=blue,mark size=1.4pt,only marks,mark=square*,mark options={solid,fill=blue},forget plot]
  table[row sep=crcr]{%
0	5.80681667455282e-30\\
0.1	0.0284437038279075\\
0.2	0.112902714165976\\
0.3	0.279976280361374\\
0.4	0.627971167735205\\
0.5	1.12142279891524\\
0.6	1.62138611466591\\
0.7	2.05583829186495\\
};
\addplot [color=green,solid,line width=1.0pt,forget plot]
  table[row sep=crcr]{%
0	4.55054309592263e-30\\
0.1	0.0266894654069442\\
0.2	0.106430027917321\\
0.3	0.267028663702478\\
0.4	0.610079900589041\\
0.5	1.1008050967236\\
0.6	1.59957855960705\\
0.7	2.03233929784081\\
};
\addplot [color=green,mark size=1.4pt,only marks,mark=square*,mark options={solid,fill=green},forget plot]
  table[row sep=crcr]{%
0	4.55054309592263e-30\\
0.1	0.0266894654069442\\
0.2	0.106430027917321\\
0.3	0.267028663702478\\
0.4	0.610079900589041\\
0.5	1.1008050967236\\
0.6	1.59957855960705\\
0.7	2.03233929784081\\
};
\addplot [color=black,solid,line width=1.0pt,forget plot]
  table[row sep=crcr]{%
0	2.00569072847376e-30\\
0.1	0.156046322729085\\
0.2	0.575419155278648\\
0.3	1.26684483550207\\
0.4	2.43154867984212\\
0.5	3.63732040545908\\
0.6	4.3419769128561\\
0.7	4.68333615219304\\
};
\addplot [color=black,mark size=1.4pt,only marks,mark=square*,mark options={solid,fill=black},forget plot]
  table[row sep=crcr]{%
0	2.00569072847376e-30\\
0.1	0.156046322729085\\
0.2	0.575419155278648\\
0.3	1.26684483550207\\
0.4	2.43154867984212\\
0.5	3.63732040545908\\
0.6	4.3419769128561\\
0.7	4.68333615219304\\
};
\end{axis}
\end{tikzpicture}%
%
 }  \hfil
   \subfigure{%
    % This file was created by matlab2tikz v0.4.7 running on MATLAB 8.4.
% Copyright (c) 2008--2014, Nico Schlömer <nico.schloemer@gmail.com>
% All rights reserved.
% Minimal pgfplots version: 1.3
% 
% The latest updates can be retrieved from
%   http://www.mathworks.com/matlabcentral/fileexchange/22022-matlab2tikz
% where you can also make suggestions and rate matlab2tikz.
% 
\begin{tikzpicture}

\begin{axis}[%
width=7cm,
height=3cm,
scale only axis,
separate axis lines,
every outer x axis line/.append style={white!15!black},
every x tick label/.append style={font=\color{white!15!black}},
xmin=0,
xmax=1,
xlabel={$\rho$},
every outer y axis line/.append style={white!15!black},
every y tick label/.append style={font=\color{white!15!black}},
ymin=-0.1,
ymax=1.5,
ylabel={g'},
title style={font=\bfseries},
title={$\sigma{=}0.3$, $n{=}30$, $d{=}3$}
]
\addplot [color=blue,solid,line width=1.0pt,forget plot]
  table[row sep=crcr]{%
0	6.32139042205196e-30\\
0.1	0.239260117751017\\
0.2	0.258739265161344\\
0.3	0.27105396872609\\
0.4	0.27444598297429\\
0.5	0.278095694101308\\
0.6	0.283180453611295\\
0.7	0.284192723563005\\
0.8	0.286037290764718\\
0.9	0.288291582869921\\
1	0.287266703898924\\
};
\addplot [color=blue,mark size=1.4pt,only marks,mark=square*,mark options={solid,fill=blue},forget plot]
  table[row sep=crcr]{%
0	6.32139042205196e-30\\
0.1	0.239260117751017\\
0.2	0.258739265161344\\
0.3	0.27105396872609\\
0.4	0.27444598297429\\
0.5	0.278095694101308\\
0.6	0.283180453611295\\
0.7	0.284192723563005\\
0.8	0.286037290764718\\
0.9	0.288291582869921\\
1	0.287266703898924\\
};
\addplot [color=green,solid,line width=1.0pt,forget plot]
  table[row sep=crcr]{%
0	2.24692433018977e-29\\
0.1	0.180805548922294\\
0.2	0.227679350991769\\
0.3	0.249450200474516\\
0.4	0.258207740321803\\
0.5	0.264629324868958\\
0.6	0.271640801985737\\
0.7	0.27478838719158\\
0.8	0.277412814261558\\
0.9	0.280759104030251\\
1	0.28025869430311\\
};
\addplot [color=green,mark size=1.4pt,only marks,mark=square*,mark options={solid,fill=green},forget plot]
  table[row sep=crcr]{%
0	2.24692433018977e-29\\
0.1	0.180805548922294\\
0.2	0.227679350991769\\
0.3	0.249450200474516\\
0.4	0.258207740321803\\
0.5	0.264629324868958\\
0.6	0.271640801985737\\
0.7	0.27478838719158\\
0.8	0.277412814261558\\
0.9	0.280759104030251\\
1	0.28025869430311\\
};
\addplot [color=black,solid,line width=1.0pt,forget plot]
  table[row sep=crcr]{%
0	6.54384687909502e-31\\
0.1	0.912110720560682\\
0.2	1.12585159477048\\
0.3	1.1140856005615\\
0.4	1.27462615513385\\
0.5	1.24912189391076\\
0.6	1.26266614354516\\
0.7	1.38060418283999\\
0.8	1.30135183332302\\
0.9	1.39988794603838\\
1	1.37244310641965\\
};
\addplot [color=black,mark size=1.4pt,only marks,mark=square*,mark options={solid,fill=black},forget plot]
  table[row sep=crcr]{%
0	6.54384687909502e-31\\
0.1	0.912110720560682\\
0.2	1.12585159477048\\
0.3	1.1140856005615\\
0.4	1.27462615513385\\
0.5	1.24912189391076\\
0.6	1.26266614354516\\
0.7	1.38060418283999\\
0.8	1.30135183332302\\
0.9	1.39988794603838\\
1	1.37244310641965\\
};
\end{axis}
\end{tikzpicture}%
%
 } 
    }
    \vspace{-0.4cm} 
  \centerline{
     \subfigure{%
    % This file was created by matlab2tikz v0.4.7 running on MATLAB 8.4.
% Copyright (c) 2008--2014, Nico Schlömer <nico.schloemer@gmail.com>
% All rights reserved.
% Minimal pgfplots version: 1.3
% 
% The latest updates can be retrieved from
%   http://www.mathworks.com/matlabcentral/fileexchange/22022-matlab2tikz
% where you can also make suggestions and rate matlab2tikz.
% 
\begin{tikzpicture}

\begin{axis}[%
width=7cm,
height=3cm,
scale only axis,
separate axis lines,
every outer x axis line/.append style={white!15!black},
every x tick label/.append style={font=\color{white!15!black}},
xmin=2,
xmax=20,
xlabel={$d$},
every outer y axis line/.append style={white!15!black},
every y tick label/.append style={font=\color{white!15!black}},
ymin=-0.1,
ymax=35.1,
ylabel={g'},
title style={font=\bfseries},
title={$\sigma{=}0.3$, $n{=}30$, $\rho{=}0.5$}
]
\addplot [color=blue,solid,line width=1.0pt,forget plot]
  table[row sep=crcr]{%
2	0.0934781094015174\\
4	0.560690028000594\\
6	1.41189402181289\\
8	2.64511293459192\\
10	4.1930366757779\\
12	6.02884242304004\\
14	8.04296942507072\\
16	10.1774555737549\\
18	12.4377846814881\\
20	14.7899325488877\\
};
\addplot [color=blue,mark size=1.4pt,only marks,mark=square*,mark options={solid,fill=blue},forget plot]
  table[row sep=crcr]{%
2	0.0934781094015174\\
4	0.560690028000594\\
6	1.41189402181289\\
8	2.64511293459192\\
10	4.1930366757779\\
12	6.02884242304004\\
14	8.04296942507072\\
16	10.1774555737549\\
18	12.4377846814881\\
20	14.7899325488877\\
};
\addplot [color=green,solid,line width=1.0pt,forget plot]
  table[row sep=crcr]{%
2	0.0885622172340653\\
4	0.53646720604058\\
6	1.36224666839895\\
8	2.56907356023743\\
10	4.09109453176731\\
12	5.90253683800283\\
14	7.88460764875725\\
16	9.99202480657338\\
18	12.2349188553887\\
20	14.5423485982679\\
};
\addplot [color=green,mark size=1.4pt,only marks,mark=square*,mark options={solid,fill=green},forget plot]
  table[row sep=crcr]{%
2	0.0885622172340653\\
4	0.53646720604058\\
6	1.36224666839895\\
8	2.56907356023743\\
10	4.09109453176731\\
12	5.90253683800283\\
14	7.88460764875725\\
16	9.99202480657338\\
18	12.2349188553887\\
20	14.5423485982679\\
};
\addplot [color=black,solid,line width=1.0pt,forget plot]
  table[row sep=crcr]{%
2	0.430366158400958\\
4	2.44583954876716\\
6	5.41203295323591\\
8	9.09147641254302\\
10	12.9982792159808\\
12	16.8466277417984\\
14	20.7496823945005\\
16	24.6214597624893\\
18	28.3926779320224\\
20	32.1816178795094\\
};
\addplot [color=black,mark size=1.4pt,only marks,mark=square*,mark options={solid,fill=black},forget plot]
  table[row sep=crcr]{%
2	0.430366158400958\\
4	2.44583954876716\\
6	5.41203295323591\\
8	9.09147641254302\\
10	12.9982792159808\\
12	16.8466277417984\\
14	20.7496823945005\\
16	24.6214597624893\\
18	28.3926779320224\\
20	32.1816178795094\\
};
\end{axis}
\end{tikzpicture}%
%
 } \hfil
     \subfigure{%
    % This file was created by matlab2tikz v0.4.7 running on MATLAB 8.4.
% Copyright (c) 2008--2014, Nico Schlömer <nico.schloemer@gmail.com>
% All rights reserved.
% Minimal pgfplots version: 1.3
% 
% The latest updates can be retrieved from
%   http://www.mathworks.com/matlabcentral/fileexchange/22022-matlab2tikz
% where you can also make suggestions and rate matlab2tikz.
% 
\begin{tikzpicture}

\begin{axis}[%
width=7cm,
height=3cm,
scale only axis,
separate axis lines,
every outer x axis line/.append style={white!15!black},
every x tick label/.append style={font=\color{white!15!black}},
xmin=10,
xmax=58,
xlabel={$n$},
every outer y axis line/.append style={white!15!black},
every y tick label/.append style={font=\color{white!15!black}},
ymin=0.1,
ymax=1.9,
ylabel={g'},
title style={font=\bfseries},
title={$\sigma{=}0.3$, $d{=}3$, $\rho{=}0.5$}
]
\addplot [color=blue,solid,line width=1.0pt,forget plot]
  table[row sep=crcr]{%
10	0.239310079743084\\
13	0.251009582593509\\
16	0.262057044154834\\
19	0.267552668937153\\
22	0.271507903166295\\
25	0.27328543178403\\
28	0.274830266288327\\
31	0.279459370238219\\
34	0.280533051688026\\
37	0.282105348571796\\
40	0.285157162951635\\
43	0.286187391992074\\
46	0.287127726152964\\
49	0.286987539878588\\
52	0.287400587609672\\
55	0.291175317244856\\
58	0.289242488544674\\
};
\addplot [color=blue,mark size=1.4pt,only marks,mark=square*,mark options={solid,fill=blue},forget plot]
  table[row sep=crcr]{%
10	0.239310079743084\\
13	0.251009582593509\\
16	0.262057044154834\\
19	0.267552668937153\\
22	0.271507903166295\\
25	0.27328543178403\\
28	0.274830266288327\\
31	0.279459370238219\\
34	0.280533051688026\\
37	0.282105348571796\\
40	0.285157162951635\\
43	0.286187391992074\\
46	0.287127726152964\\
49	0.286987539878588\\
52	0.287400587609672\\
55	0.291175317244856\\
58	0.289242488544674\\
};
\addplot [color=green,solid,line width=1.0pt,forget plot]
  table[row sep=crcr]{%
10	0.210373019354256\\
13	0.225781773014693\\
16	0.239725124263962\\
19	0.248824845326889\\
22	0.253983342008155\\
25	0.257869157315674\\
28	0.261316657579025\\
31	0.266525593169429\\
34	0.268774867275232\\
37	0.271053833744457\\
40	0.274979306769303\\
43	0.276647854105173\\
46	0.277965215714809\\
49	0.278469745459795\\
52	0.279328008042616\\
55	0.283502847502291\\
58	0.28196192051413\\
};
\addplot [color=green,mark size=1.4pt,only marks,mark=square*,mark options={solid,fill=green},forget plot]
  table[row sep=crcr]{%
10	0.210373019354256\\
13	0.225781773014693\\
16	0.239725124263962\\
19	0.248824845326889\\
22	0.253983342008155\\
25	0.257869157315674\\
28	0.261316657579025\\
31	0.266525593169429\\
34	0.268774867275232\\
37	0.271053833744457\\
40	0.274979306769303\\
43	0.276647854105173\\
46	0.277965215714809\\
49	0.278469745459795\\
52	0.279328008042616\\
55	0.283502847502291\\
58	0.28196192051413\\
};
\addplot [color=black,solid,line width=1.0pt,forget plot]
  table[row sep=crcr]{%
10	0.732670980459889\\
13	0.939165421825102\\
16	0.978355565905724\\
19	1.08444108774288\\
22	1.15123911548022\\
25	1.13705872523512\\
28	1.23138706674476\\
31	1.35080254891932\\
34	1.2754403931655\\
37	1.36506467060651\\
40	1.43947277415658\\
43	1.45859557467519\\
46	1.44410042583584\\
49	1.52998620538724\\
52	1.56870763325772\\
55	1.56044908406325\\
58	1.66104652433173\\
};
\addplot [color=black,mark size=1.4pt,only marks,mark=square*,mark options={solid,fill=black},forget plot]
  table[row sep=crcr]{%
10	0.732670980459889\\
13	0.939165421825102\\
16	0.978355565905724\\
19	1.08444108774288\\
22	1.15123911548022\\
25	1.13705872523512\\
28	1.23138706674476\\
31	1.35080254891932\\
34	1.2754403931655\\
37	1.36506467060651\\
40	1.43947277415658\\
43	1.45859557467519\\
46	1.44410042583584\\
49	1.52998620538724\\
52	1.56870763325772\\
55	1.56044908406325\\
58	1.66104652433173\\
};
\end{axis}
\end{tikzpicture}%
%
 } 
    }  
    \vspace{-0.4cm}
\caption{Normalised error in \eqref{gPrimeFcn} for the reference-based method (black), the Z-matrix method (blue) and the H-matrix method (green) when considering transformations in $O(d)$. 
In each sub-figure, a different parameter varies along the horizontal axis.}
    \label{centralizedZHrefOrthogonal} 
\end{figure*}
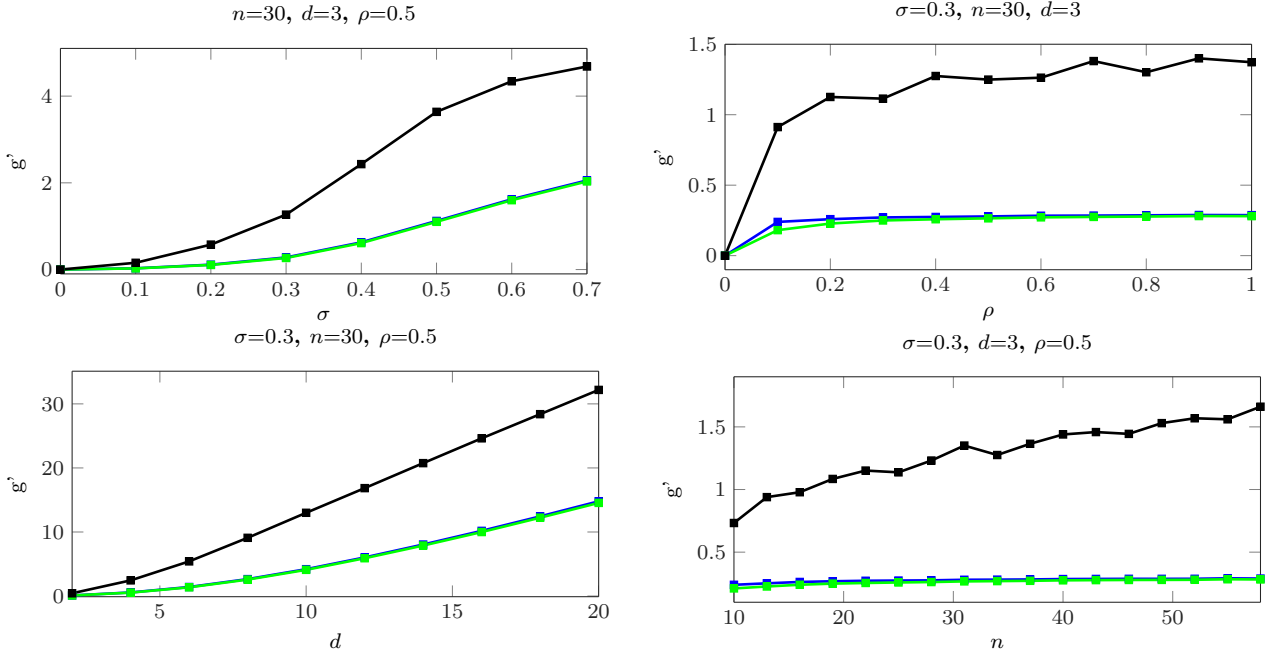 

 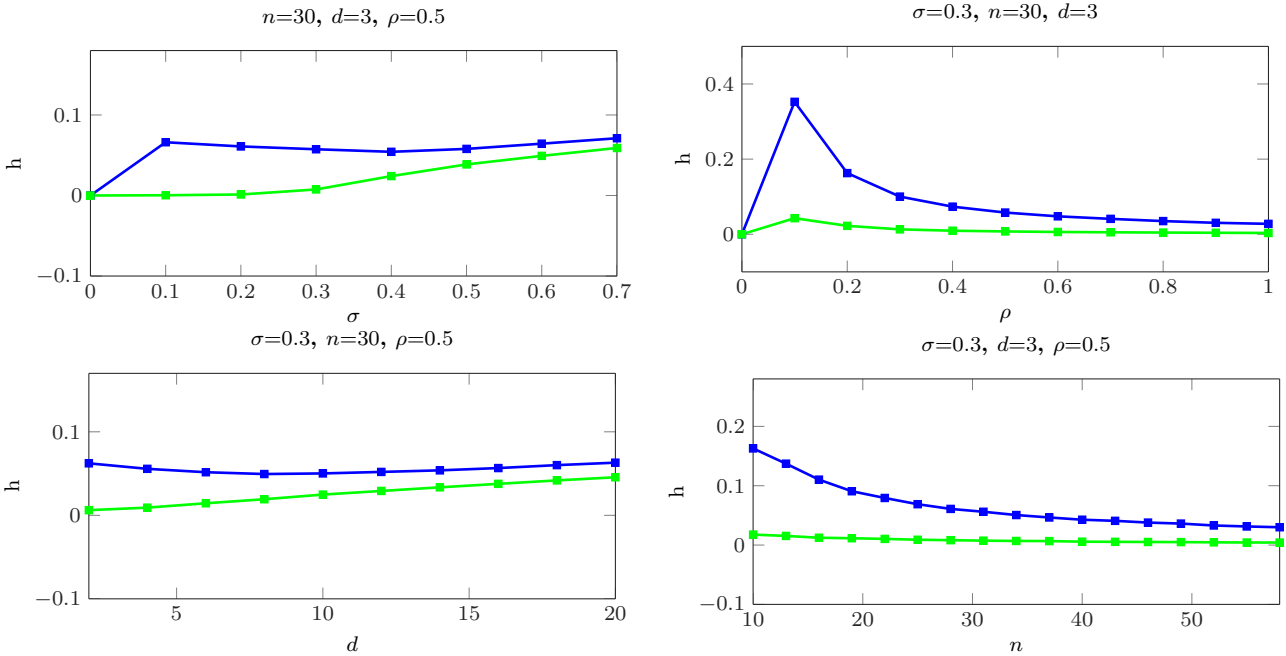
\begin{figure*}[!t]
  \centerline{
    \subfigure{%
    % This file was created by matlab2tikz v0.4.7 running on MATLAB 8.4.
% Copyright (c) 2008--2014, Nico Schlömer <nico.schloemer@gmail.com>
% All rights reserved.
% Minimal pgfplots version: 1.3
% 
% The latest updates can be retrieved from
%   http://www.mathworks.com/matlabcentral/fileexchange/22022-matlab2tikz
% where you can also make suggestions and rate matlab2tikz.
% 
\begin{tikzpicture}

\begin{axis}[%
width=7cm,
height=3cm,
scale only axis,
separate axis lines,
every outer x axis line/.append style={white!15!black},
every x tick label/.append style={font=\color{white!15!black}},
xmin=0,
xmax=0.7,
xlabel={$\sigma$},
every outer y axis line/.append style={white!15!black},
every y tick label/.append style={font=\color{white!15!black}},
ymin=-0.1,
ymax=0.18,
ylabel={h},
title style={font=\bfseries},
title={$n{=}30$, $d{=}3$, $\rho{=}0.5$}
]
\addplot [color=blue,solid,line width=1.0pt,forget plot]
  table[row sep=crcr]{%
0	0\\
0.1	0.066124303560276\\
0.2	0.0609219945332867\\
0.3	0.0574142994778704\\
0.4	0.0542464683382993\\
0.5	0.0578663315146045\\
0.6	0.0643689967233379\\
0.7	0.0710965581231556\\
};
\addplot [color=blue,mark size=1.4pt,only marks,mark=square*,mark options={solid,fill=blue},forget plot]
  table[row sep=crcr]{%
0	0\\
0.1	0.066124303560276\\
0.2	0.0609219945332867\\
0.3	0.0574142994778704\\
0.4	0.0542464683382993\\
0.5	0.0578663315146045\\
0.6	0.0643689967233379\\
0.7	0.0710965581231556\\
};
\addplot [color=green,solid,line width=1.0pt,forget plot]
  table[row sep=crcr]{%
0	0\\
0.1	0.000301430283936537\\
0.2	0.00125423066502772\\
0.3	0.00752131691677922\\
0.4	0.024085301318409\\
0.5	0.0385884513280378\\
0.6	0.0492285851761768\\
0.7	0.0589873156007689\\
};
\addplot [color=green,mark size=1.4pt,only marks,mark=square*,mark options={solid,fill=green},forget plot]
  table[row sep=crcr]{%
0	0\\
0.1	0.000301430283936537\\
0.2	0.00125423066502772\\
0.3	0.00752131691677922\\
0.4	0.024085301318409\\
0.5	0.0385884513280378\\
0.6	0.0492285851761768\\
0.7	0.0589873156007689\\
};
\end{axis}
\end{tikzpicture}%
%
 }  \hfil
   \subfigure{%
    % This file was created by matlab2tikz v0.4.7 running on MATLAB 8.4.
% Copyright (c) 2008--2014, Nico Schlömer <nico.schloemer@gmail.com>
% All rights reserved.
% Minimal pgfplots version: 1.3
% 
% The latest updates can be retrieved from
%   http://www.mathworks.com/matlabcentral/fileexchange/22022-matlab2tikz
% where you can also make suggestions and rate matlab2tikz.
% 
\begin{tikzpicture}

\begin{axis}[%
width=7cm,
height=3cm,
scale only axis,
separate axis lines,
every outer x axis line/.append style={white!15!black},
every x tick label/.append style={font=\color{white!15!black}},
xmin=0,
xmax=1,
xlabel={$\rho$},
every outer y axis line/.append style={white!15!black},
every y tick label/.append style={font=\color{white!15!black}},
ymin=-0.1,
ymax=0.5,
ylabel={h},
title style={font=\bfseries},
title={$\sigma{=}0.3$, $n{=}30$, $d{=}3$}
]
\addplot [color=blue,solid,line width=1.0pt,forget plot]
  table[row sep=crcr]{%
0	0\\
0.1	0.3524532512469\\
0.2	0.162803897033157\\
0.3	0.100282037864151\\
0.4	0.073544197283047\\
0.5	0.0578284871501858\\
0.6	0.0477601704587374\\
0.7	0.0410471303402835\\
0.8	0.0351114858690291\\
0.9	0.0305029872897975\\
1	0.0277600024933631\\
};
\addplot [color=blue,mark size=1.4pt,only marks,mark=square*,mark options={solid,fill=blue},forget plot]
  table[row sep=crcr]{%
0	0\\
0.1	0.3524532512469\\
0.2	0.162803897033157\\
0.3	0.100282037864151\\
0.4	0.073544197283047\\
0.5	0.0578284871501858\\
0.6	0.0477601704587374\\
0.7	0.0410471303402835\\
0.8	0.0351114858690291\\
0.9	0.0305029872897975\\
1	0.0277600024933631\\
};
\addplot [color=green,solid,line width=1.0pt,forget plot]
  table[row sep=crcr]{%
0	0\\
0.1	0.0428565367599553\\
0.2	0.0224588613203508\\
0.3	0.0130840106290948\\
0.4	0.00951607249615897\\
0.5	0.00756827702948972\\
0.6	0.0061321963867867\\
0.7	0.00530921389060938\\
0.8	0.00457412248743973\\
0.9	0.00401917601133287\\
1	0.0035099040162992\\
};
\addplot [color=green,mark size=1.4pt,only marks,mark=square*,mark options={solid,fill=green},forget plot]
  table[row sep=crcr]{%
0	0\\
0.1	0.0428565367599553\\
0.2	0.0224588613203508\\
0.3	0.0130840106290948\\
0.4	0.00951607249615897\\
0.5	0.00756827702948972\\
0.6	0.0061321963867867\\
0.7	0.00530921389060938\\
0.8	0.00457412248743973\\
0.9	0.00401917601133287\\
1	0.0035099040162992\\
};
\end{axis}
\end{tikzpicture}%
%
 } 
    }
    \vspace{-0.4cm} 
  \centerline{
     \subfigure{%
    % This file was created by matlab2tikz v0.4.7 running on MATLAB 8.4.
% Copyright (c) 2008--2014, Nico Schlömer <nico.schloemer@gmail.com>
% All rights reserved.
% Minimal pgfplots version: 1.3
% 
% The latest updates can be retrieved from
%   http://www.mathworks.com/matlabcentral/fileexchange/22022-matlab2tikz
% where you can also make suggestions and rate matlab2tikz.
% 
\begin{tikzpicture}

\begin{axis}[%
width=7cm,
height=3cm,
scale only axis,
separate axis lines,
every outer x axis line/.append style={white!15!black},
every x tick label/.append style={font=\color{white!15!black}},
xmin=2,
xmax=20,
xlabel={$d$},
every outer y axis line/.append style={white!15!black},
every y tick label/.append style={font=\color{white!15!black}},
ymin=-0.1,
ymax=0.17,
ylabel={h},
title style={font=\bfseries},
title={$\sigma{=}0.3$, $n{=}30$, $\rho{=}0.5$}
]
\addplot [color=blue,solid,line width=1.0pt,forget plot]
  table[row sep=crcr]{%
2	0.0622572764144456\\
4	0.0556617175387847\\
6	0.0515980709152342\\
8	0.049436470599741\\
10	0.0501651524638652\\
12	0.0519650541197466\\
14	0.0538470602224498\\
16	0.0565902829961123\\
18	0.060105663760456\\
20	0.0630440788590463\\
};
\addplot [color=blue,mark size=1.4pt,only marks,mark=square*,mark options={solid,fill=blue},forget plot]
  table[row sep=crcr]{%
2	0.0622572764144456\\
4	0.0556617175387847\\
6	0.0515980709152342\\
8	0.049436470599741\\
10	0.0501651524638652\\
12	0.0519650541197466\\
14	0.0538470602224498\\
16	0.0565902829961123\\
18	0.060105663760456\\
20	0.0630440788590463\\
};
\addplot [color=green,solid,line width=1.0pt,forget plot]
  table[row sep=crcr]{%
2	0.00617833741047696\\
4	0.00918472531005869\\
6	0.0144322184323952\\
8	0.0192624179919178\\
10	0.0248302795248757\\
12	0.0291848336801928\\
14	0.033522450054112\\
16	0.037698121558589\\
18	0.0418081205705602\\
20	0.0455452559918439\\
};
\addplot [color=green,mark size=1.4pt,only marks,mark=square*,mark options={solid,fill=green},forget plot]
  table[row sep=crcr]{%
2	0.00617833741047696\\
4	0.00918472531005869\\
6	0.0144322184323952\\
8	0.0192624179919178\\
10	0.0248302795248757\\
12	0.0291848336801928\\
14	0.033522450054112\\
16	0.037698121558589\\
18	0.0418081205705602\\
20	0.0455452559918439\\
};
\end{axis}
\end{tikzpicture}%
%
 } \hfil
     \subfigure{%
    % This file was created by matlab2tikz v0.4.7 running on MATLAB 8.4.
% Copyright (c) 2008--2014, Nico Schlömer <nico.schloemer@gmail.com>
% All rights reserved.
% Minimal pgfplots version: 1.3
% 
% The latest updates can be retrieved from
%   http://www.mathworks.com/matlabcentral/fileexchange/22022-matlab2tikz
% where you can also make suggestions and rate matlab2tikz.
% 
\begin{tikzpicture}

\begin{axis}[%
width=7cm,
height=3cm,
scale only axis,
separate axis lines,
every outer x axis line/.append style={white!15!black},
every x tick label/.append style={font=\color{white!15!black}},
xmin=10,
xmax=58,
xlabel={$n$},
every outer y axis line/.append style={white!15!black},
every y tick label/.append style={font=\color{white!15!black}},
ymin=-0.1,
ymax=0.28,
ylabel={h},
title style={font=\bfseries},
title={$\sigma{=}0.3$, $d{=}3$, $\rho{=}0.5$}
]
\addplot [color=blue,solid,line width=1.0pt,forget plot]
  table[row sep=crcr]{%
10	0.162848205440892\\
13	0.13715647763006\\
16	0.110251567176308\\
19	0.0906752418425698\\
22	0.0793457873722539\\
25	0.0688369893519325\\
28	0.0609080556179771\\
31	0.0560991135546446\\
34	0.0505711945458967\\
37	0.0464673559360037\\
40	0.0426963595671245\\
43	0.0408100531878773\\
46	0.0378667387538143\\
49	0.0361950310625396\\
52	0.032962339917849\\
55	0.0313215412160703\\
58	0.030006286722232\\
};
\addplot [color=blue,mark size=1.4pt,only marks,mark=square*,mark options={solid,fill=blue},forget plot]
  table[row sep=crcr]{%
10	0.162848205440892\\
13	0.13715647763006\\
16	0.110251567176308\\
19	0.0906752418425698\\
22	0.0793457873722539\\
25	0.0688369893519325\\
28	0.0609080556179771\\
31	0.0560991135546446\\
34	0.0505711945458967\\
37	0.0464673559360037\\
40	0.0426963595671245\\
43	0.0408100531878773\\
46	0.0378667387538143\\
49	0.0361950310625396\\
52	0.032962339917849\\
55	0.0313215412160703\\
58	0.030006286722232\\
};
\addplot [color=green,solid,line width=1.0pt,forget plot]
  table[row sep=crcr]{%
10	0.0176367427695644\\
13	0.0153893526056427\\
16	0.0123320537769597\\
19	0.0114890677030245\\
22	0.0103634984813469\\
25	0.00893814642013352\\
28	0.00820110689354867\\
31	0.00742977817812294\\
34	0.00688377122144852\\
37	0.00669292857720621\\
40	0.00568781805217147\\
43	0.00544189782825823\\
46	0.00516614996117205\\
49	0.00499000688119422\\
52	0.00457387466748022\\
55	0.00428931833635636\\
58	0.0042427562611328\\
};
\addplot [color=green,mark size=1.4pt,only marks,mark=square*,mark options={solid,fill=green},forget plot]
  table[row sep=crcr]{%
10	0.0176367427695644\\
13	0.0153893526056427\\
16	0.0123320537769597\\
19	0.0114890677030245\\
22	0.0103634984813469\\
25	0.00893814642013352\\
28	0.00820110689354867\\
31	0.00742977817812294\\
34	0.00688377122144852\\
37	0.00669292857720621\\
40	0.00568781805217147\\
43	0.00544189782825823\\
46	0.00516614996117205\\
49	0.00499000688119422\\
52	0.00457387466748022\\
55	0.00428931833635636\\
58	0.0042427562611328\\
};
\end{axis}
\end{tikzpicture}%
%
 } 
    }  
    \vspace{-0.4cm}
\caption{Gap function in \eqref{gapFcn} for the Z-matrix method (blue) and the H-matrix method (green) when considering transformations in $O(d)$. 
In each sub-figure, a different parameter varies along the horizontal axis.}
    \label{gapPlotZHOrthogonal} 
\end{figure*} 

In Fig.~\ref{centralizedZHrefOrthogonal} the results of the experiments are shown. On the vertical axis, a normalised version of the function in \eqref{gFcn}, defined by
\begin{align}
	g'(\{G_{i}\}_{i \in \mathcal{V}}) = \frac{1}{\vert \mathcal{E} \vert}\sum_{(i.j) \in \mathcal{E}} \| G_{ij} - G_i^{-1} G_j\|_F^2 \label{gPrimeFcn}
\end{align}
is used. Each sub-figure shows a different varying parameter on the horizontal axis. The title of each sub-figure indicates the fixed parameters.

It can be seen that in all cases the Z-matrix approach is nearly as good as the H-matrix approach when looking at orthogonal transformations. However, as anticipated, the reference-based method performs worse than both proposed methods. For the case of different degrees of noise (Fig.~\ref{centralizedZHrefOrthogonal}, top left) it can be seen that the total error increases with increasing noise. Similarly, in the case of different dimensions (Fig.~\ref{centralizedZHrefOrthogonal}, bottom left), the error increases with increasing dimensionality. This can be explained by the fact that the Frobenius norm in \eqref{gPrimeFcn} sums over $d^2$ values. For various values of the graph density (Fig.~\ref{centralizedZHrefOrthogonal}, top right), the error for the H- and Z-matrix method is approximately constant (apart from the case of a rooted spanning tree at $\rho = 0$, according to Lemma~\ref{lem:3}.). 

\subsection{Centralized methods for matrices in $GL(d, \mathbb{R})$}
In this set of experiments we compare the H-matrix method and the Z-matrix method. 

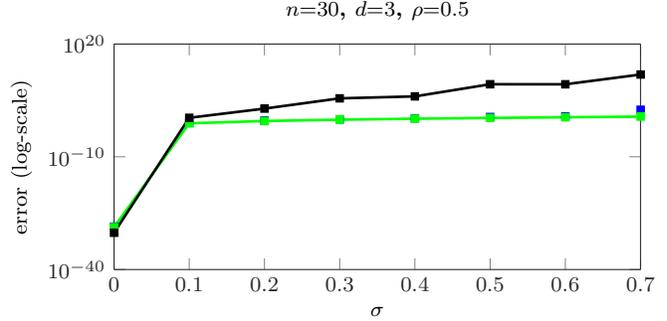
\begin{figure}
  \centerline{
    \subfigure{%
    % This file was created by matlab2tikz v0.4.7 running on MATLAB 8.4.
% Copyright (c) 2008--2014, Nico Schlömer <nico.schloemer@gmail.com>
% All rights reserved.
% Minimal pgfplots version: 1.3
% 
% The latest updates can be retrieved from
%   http://www.mathworks.com/matlabcentral/fileexchange/22022-matlab2tikz
% where you can also make suggestions and rate matlab2tikz.
% 
\begin{tikzpicture}

\begin{axis}[%
width=7cm,
height=3cm,
scale only axis,
separate axis lines,
every outer x axis line/.append style={white!15!black},
every x tick label/.append style={font=\color{white!15!black}},
xmin=0,
xmax=0.7,
xlabel={$\sigma$},
every outer y axis line/.append style={white!15!black},
every y tick label/.append style={font=\color{white!15!black}},
ymode=log,
ymin=1e-40,
ymax=1e+20,
yminorticks=true,
ylabel={error (log-scale)},
title style={font=\bfseries},
title={$n{=}30$, $d{=}3$, $\rho{=}0.5$}
]
\addplot [color=blue,mark size=1.4pt,only marks,mark=square*,mark options={solid,fill=blue},forget plot]
  table[row sep=crcr]{%
0	2.5311341973753e-29\\
0.1	0.0986057472102722\\
0.2	0.406841706646399\\
0.3	0.931430852131955\\
0.4	1.70792657583622\\
0.5	3.8660004479839\\
0.6	5.5190856102739\\
0.7	324.619262843574\\
};
\addplot [color=green,solid,line width=1.0pt,forget plot]
  table[row sep=crcr]{%
0	1.88326313057218e-29\\
0.1	0.0810429261793862\\
0.2	0.327259012366548\\
0.3	0.74859734829273\\
0.4	1.3511033561427\\
0.5	2.17230042463724\\
0.6	3.29277711869461\\
0.7	4.82956398495021\\
};
\addplot [color=green,mark size=1.4pt,only marks,mark=square*,mark options={solid,fill=green},forget plot]
  table[row sep=crcr]{%
0	1.88326313057218e-29\\
0.1	0.0810429261793862\\
0.2	0.327259012366548\\
0.3	0.74859734829273\\
0.4	1.3511033561427\\
0.5	2.17230042463724\\
0.6	3.29277711869461\\
0.7	4.82956398495021\\
};
\addplot [color=black,solid,line width=1.0pt,forget plot]
  table[row sep=crcr]{%
0	6.08822556710495e-31\\
0.1	2.35866112872715\\
0.2	664.536034663529\\
0.3	354025.745840133\\
0.4	1177292.84759338\\
0.5	2008894681.04251\\
0.6	1948962417.18984\\
0.7	775115589233.978\\
};
\addplot [color=black,mark size=1.4pt,only marks,mark=square*,mark options={solid,fill=black},forget plot]
  table[row sep=crcr]{%
0	6.08822556710495e-31\\
0.1	2.35866112872715\\
0.2	664.536034663529\\
0.3	354025.745840133\\
0.4	1177292.84759338\\
0.5	2008894681.04251\\
0.6	1948962417.18984\\
0.7	775115589233.978\\
};
\end{axis}
\end{tikzpicture}%
%
 }  
  }  
   \vspace{-0.4cm}
\caption{Normalised error in \eqref{gPrimeFcn} on the horizontal axis shown as log-scale for the Z-matrix method (blue), the H-matrix method (green) and the reference-based method (black) when considering transformations in $GL(d, \mathbb{R})$. Note that the blue and green line coincide. }
    \label{linearRefProblem} 
\end{figure}

Using the reference-based method for the case of linear transformations is problematic because this method inverts the matrices $G_{ij}$ for $(i,j) \in \mathcal{E}^{\text{min-QSC}}$ (see \eqref{gPrimeFcn}). Therefore, for reasonably large noise, it is likely that there is some $(i,j) \in \mathcal{E}^{\text{min-QSC}}$  where $G_{ij}$ is ill-conditioned, resulting in the corresponding term in $g'$ blowing up. In Fig.~\ref{linearRefProblem} this problem is illustrated, where the horizontal axis is shown in log-scale. The lines of the Z- and H-matrix methods almost coincide, so only the green line of the H-matrix method is visible. The reference-based method's (black) line results in extremely large errors. Due to this reason, and since we have already shown that for the  case of orthogonal transformations the reference-based method is inferior, in the following the reference based method is 
not used in the comparisons.

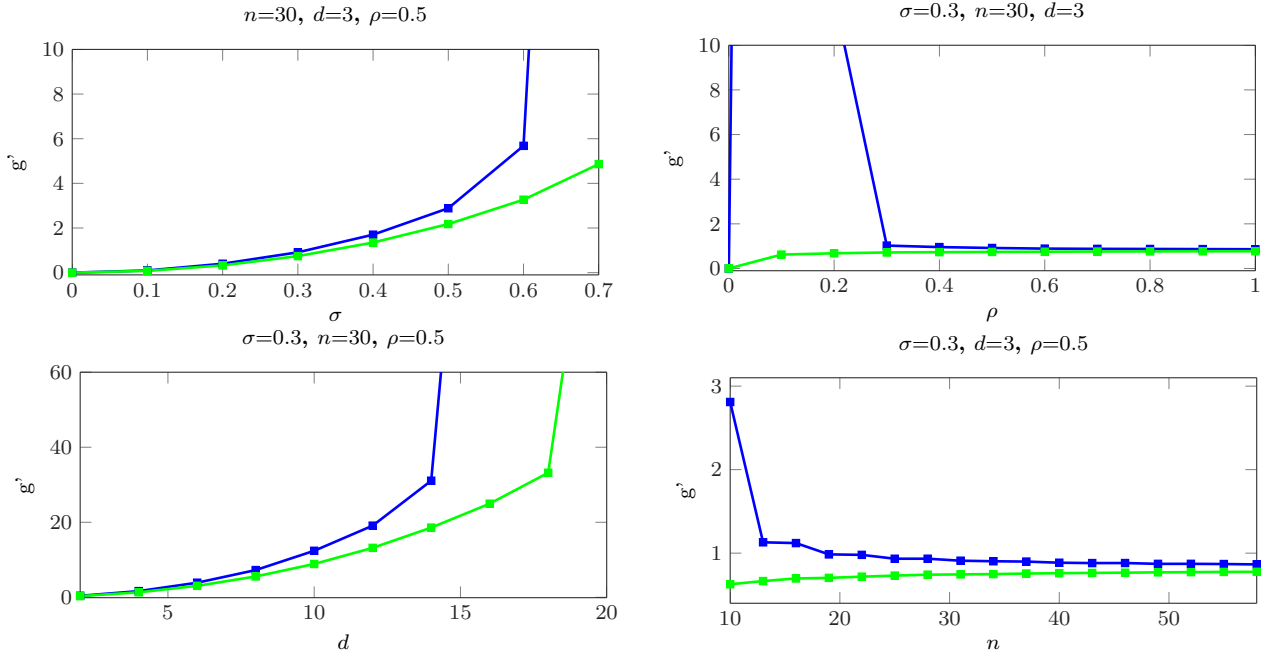
\begin{figure*}[!t] 
  \centerline{
    \subfigure{%
    % This file was created by matlab2tikz v0.4.7 running on MATLAB 8.4.
% Copyright (c) 2008--2014, Nico Schlömer <nico.schloemer@gmail.com>
% All rights reserved.
% Minimal pgfplots version: 1.3
% 
% The latest updates can be retrieved from
%   http://www.mathworks.com/matlabcentral/fileexchange/22022-matlab2tikz
% where you can also make suggestions and rate matlab2tikz.
% 
\begin{tikzpicture}

\begin{axis}[%
width=7cm,
height=3cm,
scale only axis,
separate axis lines,
every outer x axis line/.append style={white!15!black},
every x tick label/.append style={font=\color{white!15!black}},
xmin=0,
xmax=0.7,
xlabel={$\sigma$},
every outer y axis line/.append style={white!15!black},
every y tick label/.append style={font=\color{white!15!black}},
ymin=-0.1,
ymax=10,
ylabel={g'},
title style={font=\bfseries},
title={$n{=}30$, $d{=}3$, $\rho{=}0.5$}
]
\addplot [color=blue,solid,line width=1.0pt,forget plot]
  table[row sep=crcr]{%
0	3.12699193505831e-29\\
0.1	0.0983229156593044\\
0.2	0.39964883941617\\
0.3	0.910871710954999\\
0.4	1.70115788628408\\
0.5	2.88561855147042\\
0.6	5.6806006598977\\
0.7	65.4817975480593\\
};
\addplot [color=blue,mark size=1.4pt,only marks,mark=square*,mark options={solid,fill=blue},forget plot]
  table[row sep=crcr]{%
0	3.12699193505831e-29\\
0.1	0.0983229156593044\\
0.2	0.39964883941617\\
0.3	0.910871710954999\\
0.4	1.70115788628408\\
0.5	2.88561855147042\\
0.6	5.6806006598977\\
0.7	65.4817975480593\\
};
\addplot [color=green,solid,line width=1.0pt,forget plot]
  table[row sep=crcr]{%
0	1.80233942057178e-29\\
0.1	0.0811554647791873\\
0.2	0.325423403073005\\
0.3	0.74216287608608\\
0.4	1.34891024658168\\
0.5	2.17712096347347\\
0.6	3.26788078584297\\
0.7	4.8623254542242\\
};
\addplot [color=green,mark size=1.4pt,only marks,mark=square*,mark options={solid,fill=green},forget plot]
  table[row sep=crcr]{%
0	1.80233942057178e-29\\
0.1	0.0811554647791873\\
0.2	0.325423403073005\\
0.3	0.74216287608608\\
0.4	1.34891024658168\\
0.5	2.17712096347347\\
0.6	3.26788078584297\\
0.7	4.8623254542242\\
};
\end{axis}
\end{tikzpicture}%
%
 } \hfil 
   \subfigure{%
    % This file was created by matlab2tikz v0.4.7 running on MATLAB 8.4.
% Copyright (c) 2008--2014, Nico Schlömer <nico.schloemer@gmail.com>
% All rights reserved.
% Minimal pgfplots version: 1.3
% 
% The latest updates can be retrieved from
%   http://www.mathworks.com/matlabcentral/fileexchange/22022-matlab2tikz
% where you can also make suggestions and rate matlab2tikz.
% 
\begin{tikzpicture}

\begin{axis}[%
width=7cm,
height=3cm,
scale only axis,
separate axis lines,
every outer x axis line/.append style={white!15!black},
every x tick label/.append style={font=\color{white!15!black}},
xmin=0,
xmax=1,
xlabel={$\rho$},
every outer y axis line/.append style={white!15!black},
every y tick label/.append style={font=\color{white!15!black}},
ymin=-0.1,
ymax=10,
ylabel={g'},
title style={font=\bfseries},
title={$\sigma{=}0.3$, $n{=}30$, $d{=}3$}
]
\addplot [color=blue,solid,line width=1.0pt,forget plot]
  table[row sep=crcr]{%
0	1.04404921157565e-27\\
0.1	192.515407274243\\
0.2	12.1684718184638\\
0.3	1.0257406903641\\
0.4	0.958515701030384\\
0.5	0.918448549557185\\
0.6	0.890260461718291\\
0.7	0.879690895583096\\
0.8	0.874822980495786\\
0.9	0.8681181592315\\
1	0.862611801635616\\
};
\addplot [color=blue,mark size=1.4pt,only marks,mark=square*,mark options={solid,fill=blue},forget plot]
  table[row sep=crcr]{%
0	1.04404921157565e-27\\
0.1	192.515407274243\\
0.2	12.1684718184638\\
0.3	1.0257406903641\\
0.4	0.958515701030384\\
0.5	0.918448549557185\\
0.6	0.890260461718291\\
0.7	0.879690895583096\\
0.8	0.874822980495786\\
0.9	0.8681181592315\\
1	0.862611801635616\\
};
\addplot [color=green,solid,line width=1.0pt,forget plot]
  table[row sep=crcr]{%
0	6.8286298999097e-27\\
0.1	0.622581305905365\\
0.2	0.68074905906868\\
0.3	0.717795013379693\\
0.4	0.731652528744499\\
0.5	0.742240117416676\\
0.6	0.751677872132101\\
0.7	0.759872366818966\\
0.8	0.76605954633394\\
0.9	0.770894495823145\\
1	0.772677158425669\\
};
\addplot [color=green,mark size=1.4pt,only marks,mark=square*,mark options={solid,fill=green},forget plot]
  table[row sep=crcr]{%
0	6.8286298999097e-27\\
0.1	0.622581305905365\\
0.2	0.68074905906868\\
0.3	0.717795013379693\\
0.4	0.731652528744499\\
0.5	0.742240117416676\\
0.6	0.751677872132101\\
0.7	0.759872366818966\\
0.8	0.76605954633394\\
0.9	0.770894495823145\\
1	0.772677158425669\\
};
\end{axis}
\end{tikzpicture}%
%
 } 
    } \vspace{-0.4cm}
  \centerline{
     \subfigure{%
    % This file was created by matlab2tikz v0.4.7 running on MATLAB 8.4.
% Copyright (c) 2008--2014, Nico Schlömer <nico.schloemer@gmail.com>
% All rights reserved.
% Minimal pgfplots version: 1.3
% 
% The latest updates can be retrieved from
%   http://www.mathworks.com/matlabcentral/fileexchange/22022-matlab2tikz
% where you can also make suggestions and rate matlab2tikz.
% 
\begin{tikzpicture}

\begin{axis}[%
width=7cm,
height=3cm,
scale only axis,
separate axis lines,
every outer x axis line/.append style={white!15!black},
every x tick label/.append style={font=\color{white!15!black}},
xmin=2,
xmax=20,
xlabel={$d$},
every outer y axis line/.append style={white!15!black},
every y tick label/.append style={font=\color{white!15!black}},
ymin=-0.1,
ymax=60,
ylabel={g'},
title style={font=\bfseries},
title={$\sigma{=}0.3$, $n{=}30$, $\rho{=}0.5$}
]
\addplot [color=blue,solid,line width=1.0pt,forget plot]
  table[row sep=crcr]{%
2	0.404276104322468\\
4	1.65108988917339\\
6	3.8940353330284\\
8	7.26114312114853\\
10	12.3999244471246\\
12	19.1029402601858\\
14	31.0803216858139\\
16	203.559572357374\\
18	6162.71751988608\\
20	1660.69557509698\\
};
\addplot [color=blue,mark size=1.4pt,only marks,mark=square*,mark options={solid,fill=blue},forget plot]
  table[row sep=crcr]{%
2	0.404276104322468\\
4	1.65108988917339\\
6	3.8940353330284\\
8	7.26114312114853\\
10	12.3999244471246\\
12	19.1029402601858\\
14	31.0803216858139\\
16	203.559572357374\\
18	6162.71751988608\\
20	1660.69557509698\\
};
\addplot [color=green,solid,line width=1.0pt,forget plot]
  table[row sep=crcr]{%
2	0.327959686523842\\
4	1.33362017080348\\
6	3.06402299826808\\
8	5.55917933813512\\
10	8.88078945866167\\
12	13.182714495086\\
14	18.5754808229707\\
16	24.9590018297355\\
18	33.2011538432843\\
20	139.710417652175\\
};
\addplot [color=green,mark size=1.4pt,only marks,mark=square*,mark options={solid,fill=green},forget plot]
  table[row sep=crcr]{%
2	0.327959686523842\\
4	1.33362017080348\\
6	3.06402299826808\\
8	5.55917933813512\\
10	8.88078945866167\\
12	13.182714495086\\
14	18.5754808229707\\
16	24.9590018297355\\
18	33.2011538432843\\
20	139.710417652175\\
};
\end{axis}
\end{tikzpicture}%
%
 } \hfil
     \subfigure{%
    % This file was created by matlab2tikz v0.4.7 running on MATLAB 8.4.
% Copyright (c) 2008--2014, Nico Schlömer <nico.schloemer@gmail.com>
% All rights reserved.
% Minimal pgfplots version: 1.3
% 
% The latest updates can be retrieved from
%   http://www.mathworks.com/matlabcentral/fileexchange/22022-matlab2tikz
% where you can also make suggestions and rate matlab2tikz.
% 
\begin{tikzpicture}

\begin{axis}[%
width=7cm,
height=3cm,
scale only axis,
separate axis lines,
every outer x axis line/.append style={white!15!black},
every x tick label/.append style={font=\color{white!15!black}},
xmin=10,
xmax=58,
xlabel={$n$},
every outer y axis line/.append style={white!15!black},
every y tick label/.append style={font=\color{white!15!black}},
ymin=0.4,
ymax=3.1,
ylabel={g'},
title style={font=\bfseries},
title={$\sigma{=}0.3$, $d{=}3$, $\rho{=}0.5$}
]
\addplot [color=blue,solid,line width=1.0pt,forget plot]
  table[row sep=crcr]{%
10	2.80983430702168\\
13	1.12971562398437\\
16	1.12057603320301\\
19	0.985179403547327\\
22	0.979460188470179\\
25	0.933162442631108\\
28	0.933513732539373\\
31	0.909483265026514\\
34	0.903001607171553\\
37	0.898962835108604\\
40	0.885454207048188\\
43	0.880986007972856\\
46	0.881606688078487\\
49	0.87168412971748\\
52	0.872015515585187\\
55	0.869655948455972\\
58	0.864958326292365\\
};
\addplot [color=blue,mark size=1.4pt,only marks,mark=square*,mark options={solid,fill=blue},forget plot]
  table[row sep=crcr]{%
10	2.80983430702168\\
13	1.12971562398437\\
16	1.12057603320301\\
19	0.985179403547327\\
22	0.979460188470179\\
25	0.933162442631108\\
28	0.933513732539373\\
31	0.909483265026514\\
34	0.903001607171553\\
37	0.898962835108604\\
40	0.885454207048188\\
43	0.880986007972856\\
46	0.881606688078487\\
49	0.87168412971748\\
52	0.872015515585187\\
55	0.869655948455972\\
58	0.864958326292365\\
};
\addplot [color=green,solid,line width=1.0pt,forget plot]
  table[row sep=crcr]{%
10	0.627365930668897\\
13	0.665627741620765\\
16	0.696770519173221\\
19	0.704129506373818\\
22	0.718315513099444\\
25	0.7310596875744\\
28	0.742172713394742\\
31	0.745406264178671\\
34	0.749184552806228\\
37	0.755921004301858\\
40	0.760094274407286\\
43	0.763317245705451\\
46	0.765931938580455\\
49	0.770128133643531\\
52	0.772396157418625\\
55	0.775370117921994\\
58	0.776547076498317\\
};
\addplot [color=green,mark size=1.4pt,only marks,mark=square*,mark options={solid,fill=green},forget plot]
  table[row sep=crcr]{%
10	0.627365930668897\\
13	0.665627741620765\\
16	0.696770519173221\\
19	0.704129506373818\\
22	0.718315513099444\\
25	0.7310596875744\\
28	0.742172713394742\\
31	0.745406264178671\\
34	0.749184552806228\\
37	0.755921004301858\\
40	0.760094274407286\\
43	0.763317245705451\\
46	0.765931938580455\\
49	0.770128133643531\\
52	0.772396157418625\\
55	0.775370117921994\\
58	0.776547076498317\\
};
\end{axis}
\end{tikzpicture}%
%
 } 
    }  \vspace{-0.4cm}
\caption{Normalised error in \eqref{gPrimeFcn} for the Z-matrix method (blue) and the H-matrix method (green) when considering transformations in $GL(d, \mathbb{R})$. 
In each sub-figure, a different parameter varies along the horizontal axis.}
    \label{centralizedZHrefLinear} 
\end{figure*} 

\begin{figure}[!th] 
\centering
\includegraphics[scale=0.46]{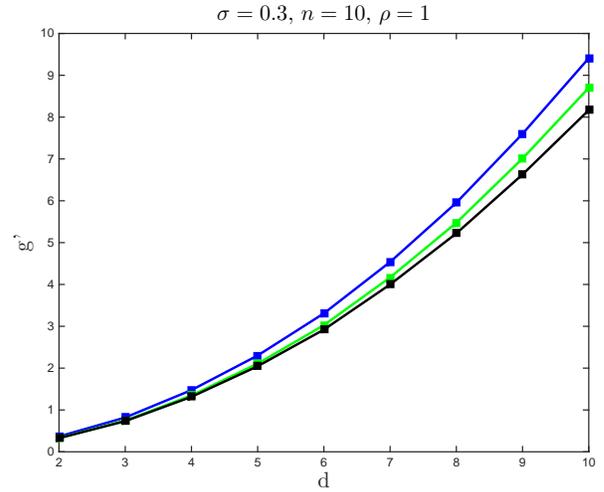}
\caption{Performance of the $Z$-matrix method (blue), the $H$-matrix method (green), and the Gauss-Newton method with the solution of the $H$-matrix method as initialization (black).}
 \label{fig_333}
\end{figure}
For the complete graph case, in Figure~\ref{fig_333} the improved performance the Gauss-Newton method, i.e., Algorithm 3, (with the solution of the $H$-matrix method as initialization)
is shown. The Gauss-Newton method run 5 iterations, but from inspection it
could be deduced that the main convergence occurs already after two iterations.

In Fig.~\ref{centralizedZHrefLinear}, the comparisons of the $H$-matrix method and the $Z$-matrix method are shown. It can be seen that for small noise (Fig.~\ref{centralizedZHrefLinear}, top left) both methods are comparable, whereas for a larger amount of noise the H-matrix method is able to obtain a smaller error. Similarly, for transformations with small dimensionality (Fig.~\ref{centralizedZHrefLinear}, bottom left), both methods are comparable whereas for larger dimensions the gap between both approaches increases. On the contrary, (Fig.~\ref{centralizedZHrefLinear}, top right) illustrates that with increasing graph density the line of the Z-matrix method approaches that of the H-matrix method (apart from the spanning tree case when $\rho = 0$, analogous to  the orthogonal transformation experiments). This indicates that the H-matrix method performs better than the Z-matrix method if there is only little information available. A similar observation can be made for various $n$ (Fig.~\ref{centralizedZHrefLinear}, bottom right).
For each subfigure, 100 simulations for a certain
configuration of $\sigma$, $n$, $d$, and $\rho$ 
are shown.

\subsection{Methods for affine and Euclidean transformations}
\begin{figure*}[t!]
    \centering
    \begin{subfigure}
        \centering
        \includegraphics[height=6.5cm]{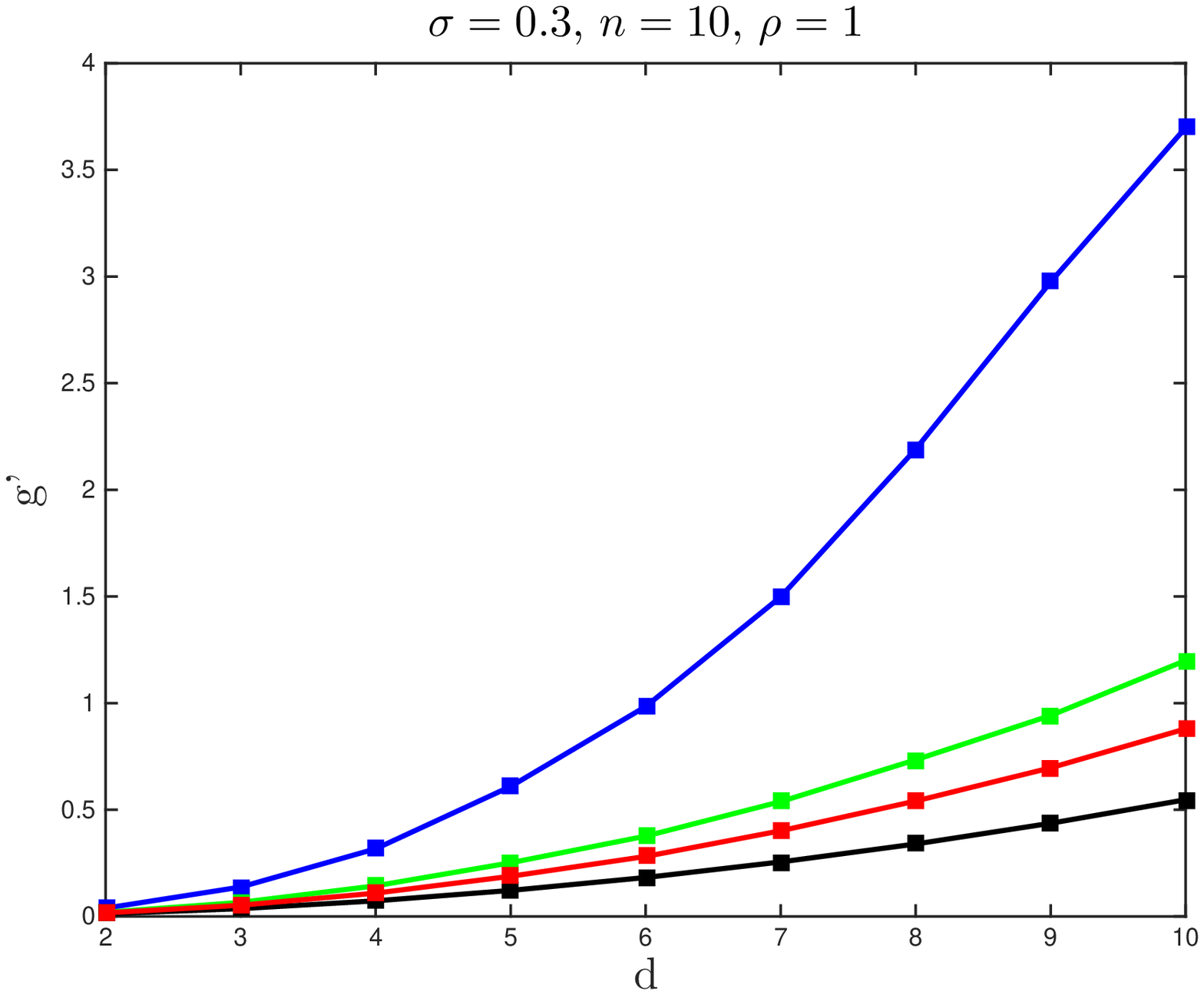}
    \end{subfigure}%
    \begin{subfigure}
        \centering
        \includegraphics[height=6.5cm]{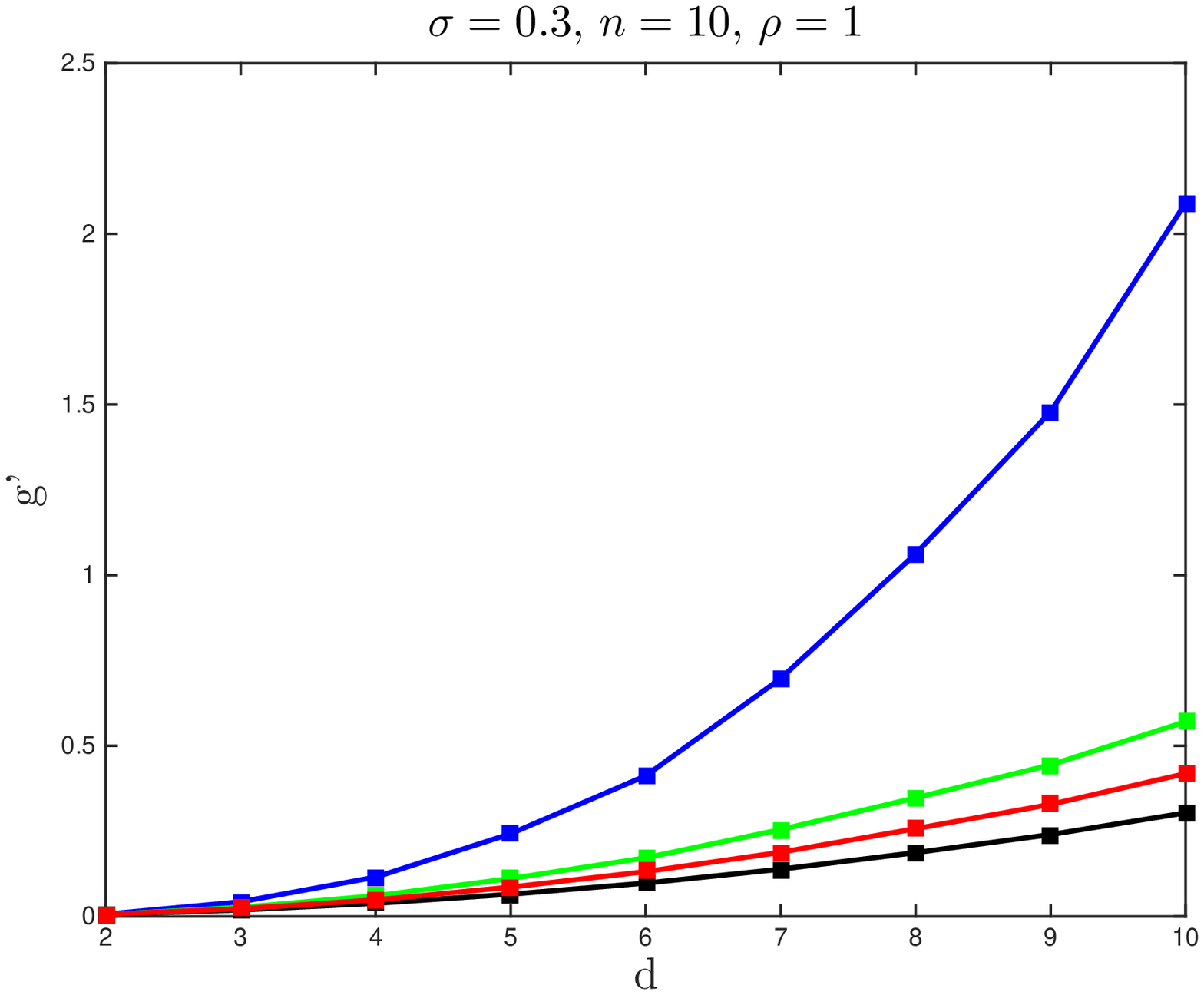}
    \end{subfigure}%
    \caption{Left figure: Performance of the methods for affine transformations. The $Z$-matrix method (blue), the $H$-matrix method (green), the first two steps of Algorithm 4 (red), and Algorithm 4 (black). Right figure: Performance of the methods for Euclidean transformations. The $Z$-matrix method (blue), the $H$-matrix method (green), the first two steps of Algorithm 4  where Algorithm 5 has been used instead of Algorithm 2 (red), and Algorithm 4 where Algorithm 5 has been used instead of Algorithm 2 (black). 
    }
    \label{fig:affine_euclidean}
\end{figure*}

In Figure~\ref{fig:affine_euclidean} -- for affine and Euclidean transformations -- a comparison
between four different methods can be found. The $G_{ij}$ transformations are affine respective Euclidean, but only the Algorithm 4 methods (red and black) preserve this property. In the bottom right figure the $G_i$ transformations obtained in the $Z$-matrix method respective the $H$-matrix method have been projected onto $E(d)$, i.e., the set of Euclidean transformations. The orthogonal matrix part of the $G_{ij}$ transformations were generated according to the description above. The elements in the transnational vectors were drawn from the uniform distribution over $(-5,5)$ and additional element-wise noise was added.

\subsection{Distributed methods}
\begin{figure*}
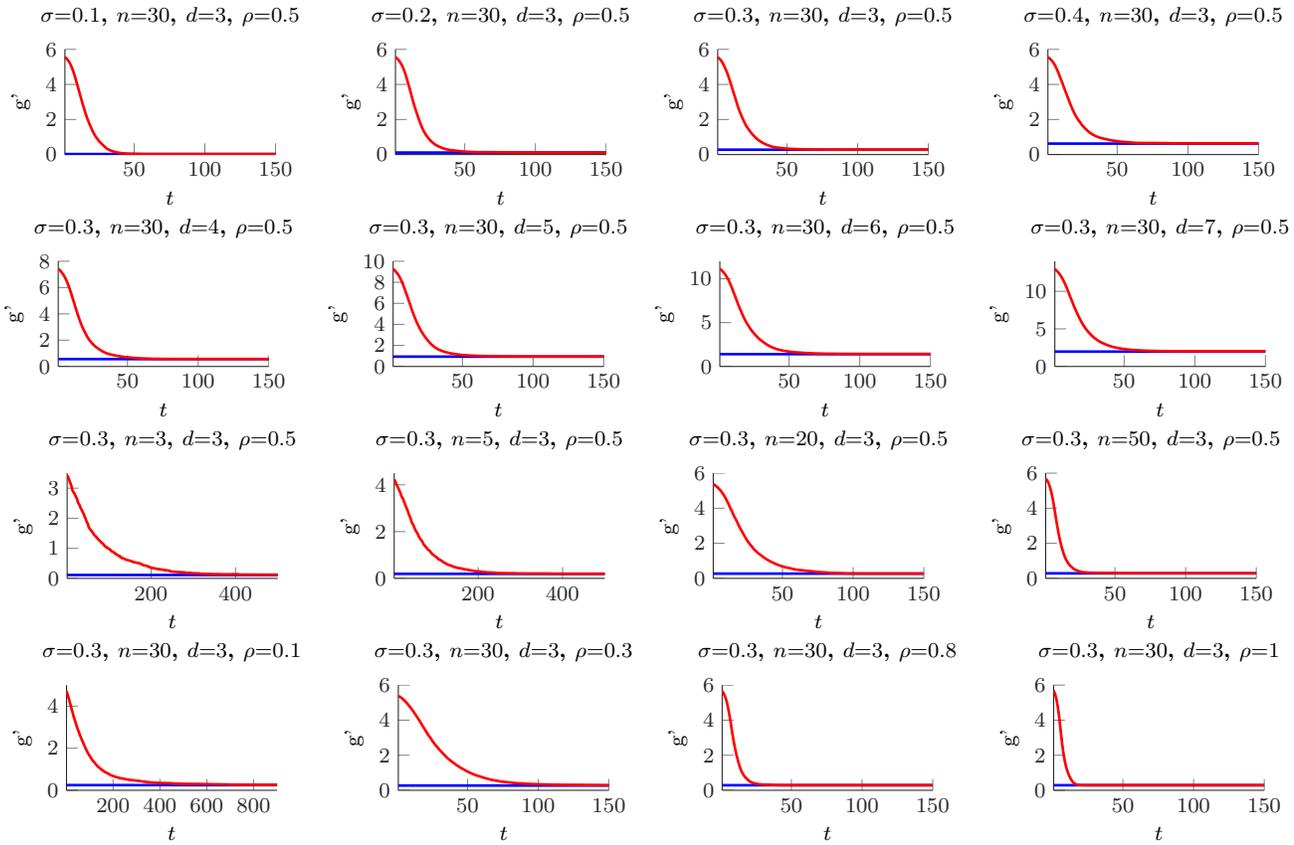

  \centerline{
    \subfigure{%
    % This file was created by matlab2tikz v0.4.7 running on MATLAB 8.4.
% Copyright (c) 2008--2014, Nico Schlömer <nico.schloemer@gmail.com>
% All rights reserved.
% Minimal pgfplots version: 1.3
% 
% The latest updates can be retrieved from
%   http://www.mathworks.com/matlabcentral/fileexchange/22022-matlab2tikz
% where you can also make suggestions and rate matlab2tikz.
% 
\begin{tikzpicture}

\begin{axis}[%
width=2.8cm,
height=1.4cm,
scale only axis,
every outer x axis line/.append style={white!15!black},
every x tick label/.append style={font=\color{white!15!black}},
xmin=1,
xmax=150,
xlabel={$t$},
every outer y axis line/.append style={white!15!black},
every y tick label/.append style={font=\color{white!15!black}},
ymin=0,
ymax=6,
ylabel={g'},
title style={font=\bfseries},
title={$\sigma{=}0.1$, $n{=}30$, $d{=}3$, $\rho{=}0.5$},
axis x line*=bottom,
axis y line*=left
]
\addplot [color=blue,solid,line width=1.0pt,forget plot]
  table[row sep=crcr]{%
1	0.028010054879727\\
2	0.028010054879727\\
3	0.028010054879727\\
4	0.028010054879727\\
5	0.028010054879727\\
6	0.028010054879727\\
7	0.028010054879727\\
8	0.028010054879727\\
9	0.028010054879727\\
10	0.028010054879727\\
11	0.028010054879727\\
12	0.028010054879727\\
13	0.028010054879727\\
14	0.028010054879727\\
15	0.028010054879727\\
16	0.028010054879727\\
17	0.028010054879727\\
18	0.028010054879727\\
19	0.028010054879727\\
20	0.028010054879727\\
21	0.028010054879727\\
22	0.028010054879727\\
23	0.028010054879727\\
24	0.028010054879727\\
25	0.028010054879727\\
26	0.028010054879727\\
27	0.028010054879727\\
28	0.028010054879727\\
29	0.028010054879727\\
30	0.028010054879727\\
31	0.028010054879727\\
32	0.028010054879727\\
33	0.028010054879727\\
34	0.028010054879727\\
35	0.028010054879727\\
36	0.028010054879727\\
37	0.028010054879727\\
38	0.028010054879727\\
39	0.028010054879727\\
40	0.028010054879727\\
41	0.028010054879727\\
42	0.028010054879727\\
43	0.028010054879727\\
44	0.028010054879727\\
45	0.028010054879727\\
46	0.028010054879727\\
47	0.028010054879727\\
48	0.028010054879727\\
49	0.028010054879727\\
50	0.028010054879727\\
51	0.028010054879727\\
52	0.028010054879727\\
53	0.028010054879727\\
54	0.028010054879727\\
55	0.028010054879727\\
56	0.028010054879727\\
57	0.028010054879727\\
58	0.028010054879727\\
59	0.028010054879727\\
60	0.028010054879727\\
61	0.028010054879727\\
62	0.028010054879727\\
63	0.028010054879727\\
64	0.028010054879727\\
65	0.028010054879727\\
66	0.028010054879727\\
67	0.028010054879727\\
68	0.028010054879727\\
69	0.028010054879727\\
70	0.028010054879727\\
71	0.028010054879727\\
72	0.028010054879727\\
73	0.028010054879727\\
74	0.028010054879727\\
75	0.028010054879727\\
76	0.028010054879727\\
77	0.028010054879727\\
78	0.028010054879727\\
79	0.028010054879727\\
80	0.028010054879727\\
81	0.028010054879727\\
82	0.028010054879727\\
83	0.028010054879727\\
84	0.028010054879727\\
85	0.028010054879727\\
86	0.028010054879727\\
87	0.028010054879727\\
88	0.028010054879727\\
89	0.028010054879727\\
90	0.028010054879727\\
91	0.028010054879727\\
92	0.028010054879727\\
93	0.028010054879727\\
94	0.028010054879727\\
95	0.028010054879727\\
96	0.028010054879727\\
97	0.028010054879727\\
98	0.028010054879727\\
99	0.028010054879727\\
100	0.028010054879727\\
101	0.028010054879727\\
102	0.028010054879727\\
103	0.028010054879727\\
104	0.028010054879727\\
105	0.028010054879727\\
106	0.028010054879727\\
107	0.028010054879727\\
108	0.028010054879727\\
109	0.028010054879727\\
110	0.028010054879727\\
111	0.028010054879727\\
112	0.028010054879727\\
113	0.028010054879727\\
114	0.028010054879727\\
115	0.028010054879727\\
116	0.028010054879727\\
117	0.028010054879727\\
118	0.028010054879727\\
119	0.028010054879727\\
120	0.028010054879727\\
121	0.028010054879727\\
122	0.028010054879727\\
123	0.028010054879727\\
124	0.028010054879727\\
125	0.028010054879727\\
126	0.028010054879727\\
127	0.028010054879727\\
128	0.028010054879727\\
129	0.028010054879727\\
130	0.028010054879727\\
131	0.028010054879727\\
132	0.028010054879727\\
133	0.028010054879727\\
134	0.028010054879727\\
135	0.028010054879727\\
136	0.028010054879727\\
137	0.028010054879727\\
138	0.028010054879727\\
139	0.028010054879727\\
140	0.028010054879727\\
141	0.028010054879727\\
142	0.028010054879727\\
143	0.028010054879727\\
144	0.028010054879727\\
145	0.028010054879727\\
146	0.028010054879727\\
147	0.028010054879727\\
148	0.028010054879727\\
149	0.028010054879727\\
150	0.028010054879727\\
};
\addplot [color=red,solid,line width=1.0pt,forget plot]
  table[row sep=crcr]{%
1	5.56677562490417\\
2	5.48867529256189\\
3	5.40054984937319\\
4	5.28310445289731\\
5	5.12890410030279\\
6	4.94785065603521\\
7	4.72327694789607\\
8	4.4738807500907\\
9	4.18882128950525\\
10	3.91353649837592\\
11	3.61454386677342\\
12	3.33130873992616\\
13	3.03533016295579\\
14	2.76253513923536\\
15	2.49871648507022\\
16	2.25050348126235\\
17	2.02124926836137\\
18	1.818880334801\\
19	1.631702352866\\
20	1.44857436261254\\
21	1.26962389312252\\
22	1.12520382771787\\
23	0.989743464781294\\
24	0.888422995613523\\
25	0.775359134718454\\
26	0.688766086235572\\
27	0.60986689061672\\
28	0.53468078670718\\
29	0.455408354066627\\
30	0.385354161399799\\
31	0.322557651612076\\
32	0.280923838445605\\
33	0.242157930637512\\
34	0.210560243433036\\
35	0.182610612690289\\
36	0.160533588946246\\
37	0.136990831475536\\
38	0.121643344953454\\
39	0.103420408378897\\
40	0.0937911318590755\\
41	0.0816685439833932\\
42	0.0722106435893778\\
43	0.0656396431817505\\
44	0.0609260307504243\\
45	0.0569775211870439\\
46	0.0510780309155024\\
47	0.0481984815621486\\
48	0.0457643937303849\\
49	0.0418919206630245\\
50	0.0401170556560008\\
51	0.03861688839701\\
52	0.0373486269103898\\
53	0.0362762364543531\\
54	0.0335032255568534\\
55	0.0327279303855444\\
56	0.0320712711540479\\
57	0.0315143105968237\\
58	0.0310410573216526\\
59	0.0306380451922827\\
60	0.0302939733548873\\
61	0.0299993919747648\\
62	0.0297464257235477\\
63	0.0295285302603418\\
64	0.0293402785513376\\
65	0.0291771746636567\\
66	0.0290354930055023\\
67	0.0289121410626715\\
68	0.0288045436265852\\
69	0.0287105464063075\\
70	0.0286283368245908\\
71	0.0285563797534613\\
72	0.0284933659687611\\
73	0.028438171200929\\
74	0.0283898238240647\\
75	0.0283474794408424\\
76	0.0283104008661677\\
77	0.0282779422661802\\
78	0.0282495364528041\\
79	0.0282246845540555\\
80	0.0282029474688231\\
81	0.0281839386690384\\
82	0.028167318033108\\
83	0.0281527864857921\\
84	0.0281400812862279\\
85	0.0281289718526415\\
86	0.028119256044269\\
87	0.0281107568421937\\
88	0.0281033193844959\\
89	0.028096808319797\\
90	0.0280911054487636\\
91	0.0280861076266802\\
92	0.0280817249026289\\
93	0.0280778788726453\\
94	0.0280745012257771\\
95	0.0280715324634063\\
96	0.0280689207736025\\
97	0.0280666210436705\\
98	0.0280645939954436\\
99	0.0280628054292417\\
100	0.0280612255637398\\
101	0.0280598284602605\\
102	0.0280585915212024\\
103	0.028057495053434\\
104	0.0280565218885109\\
105	0.0280556570525161\\
106	0.0280548874791767\\
107	0.0280542017606763\\
108	0.0280535899312717\\
109	0.0280530432794306\\
110	0.0280525541847516\\
111	0.0280521159764037\\
112	0.0280517228102453\\
113	0.0280513695621497\\
114	0.0280510517353896\\
115	0.0280507653802132\\
116	0.0280505070239939\\
117	0.0280502736105462\\
118	0.0280500624473911\\
119	0.0280498711599134\\
120	0.0280496976514968\\
121	0.0280495400688414\\
122	0.028049396771779\\
123	0.028049266306988\\
124	0.0280491473850952\\
125	0.0280490388607145\\
126	0.0280489397150378\\
127	0.0280488490406418\\
128	0.0280487660282195\\
129	0.0280486899549856\\
130	0.0280486201745367\\
131	0.028048556107977\\
132	0.0280484972361457\\
133	0.0280484430928038\\
134	0.0280483932586555\\
135	0.028048347356099\\
136	0.0280483050446115\\
137	0.0280482660166893\\
138	0.0280482299942715\\
139	0.0280481967255869\\
140	0.0280481659823706\\
141	0.0280481375574042\\
142	0.0280481112623402\\
143	0.0280480869257745\\
144	0.0280480643915361\\
145	0.0280480435171702\\
146	0.0280480241725875\\
147	0.0280480062388627\\
148	0.0280479896071636\\
149	0.028047974177794\\
150	0.0280479598593392\\
};
\end{axis}
\end{tikzpicture}%
%
 } \hfil 
   \subfigure{%
    % This file was created by matlab2tikz v0.4.7 running on MATLAB 8.4.
% Copyright (c) 2008--2014, Nico Schlömer <nico.schloemer@gmail.com>
% All rights reserved.
% Minimal pgfplots version: 1.3
% 
% The latest updates can be retrieved from
%   http://www.mathworks.com/matlabcentral/fileexchange/22022-matlab2tikz
% where you can also make suggestions and rate matlab2tikz.
% 
\begin{tikzpicture}

\begin{axis}[%
width=2.8cm,
height=1.4cm,
scale only axis,
every outer x axis line/.append style={white!15!black},
every x tick label/.append style={font=\color{white!15!black}},
xmin=1,
xmax=150,
xlabel={$t$},
every outer y axis line/.append style={white!15!black},
every y tick label/.append style={font=\color{white!15!black}},
ymin=0,
ymax=6,
ylabel={g'},
title style={font=\bfseries},
title={$\sigma{=}0.2$, $n{=}30$, $d{=}3$, $\rho{=}0.5$},
axis x line*=bottom,
axis y line*=left
]
\addplot [color=blue,solid,line width=1.0pt,forget plot]
  table[row sep=crcr]{%
1	0.113155916019526\\
2	0.113155916019526\\
3	0.113155916019526\\
4	0.113155916019526\\
5	0.113155916019526\\
6	0.113155916019526\\
7	0.113155916019526\\
8	0.113155916019526\\
9	0.113155916019526\\
10	0.113155916019526\\
11	0.113155916019526\\
12	0.113155916019526\\
13	0.113155916019526\\
14	0.113155916019526\\
15	0.113155916019526\\
16	0.113155916019526\\
17	0.113155916019526\\
18	0.113155916019526\\
19	0.113155916019526\\
20	0.113155916019526\\
21	0.113155916019526\\
22	0.113155916019526\\
23	0.113155916019526\\
24	0.113155916019526\\
25	0.113155916019526\\
26	0.113155916019526\\
27	0.113155916019526\\
28	0.113155916019526\\
29	0.113155916019526\\
30	0.113155916019526\\
31	0.113155916019526\\
32	0.113155916019526\\
33	0.113155916019526\\
34	0.113155916019526\\
35	0.113155916019526\\
36	0.113155916019526\\
37	0.113155916019526\\
38	0.113155916019526\\
39	0.113155916019526\\
40	0.113155916019526\\
41	0.113155916019526\\
42	0.113155916019526\\
43	0.113155916019526\\
44	0.113155916019526\\
45	0.113155916019526\\
46	0.113155916019526\\
47	0.113155916019526\\
48	0.113155916019526\\
49	0.113155916019526\\
50	0.113155916019526\\
51	0.113155916019526\\
52	0.113155916019526\\
53	0.113155916019526\\
54	0.113155916019526\\
55	0.113155916019526\\
56	0.113155916019526\\
57	0.113155916019526\\
58	0.113155916019526\\
59	0.113155916019526\\
60	0.113155916019526\\
61	0.113155916019526\\
62	0.113155916019526\\
63	0.113155916019526\\
64	0.113155916019526\\
65	0.113155916019526\\
66	0.113155916019526\\
67	0.113155916019526\\
68	0.113155916019526\\
69	0.113155916019526\\
70	0.113155916019526\\
71	0.113155916019526\\
72	0.113155916019526\\
73	0.113155916019526\\
74	0.113155916019526\\
75	0.113155916019526\\
76	0.113155916019526\\
77	0.113155916019526\\
78	0.113155916019526\\
79	0.113155916019526\\
80	0.113155916019526\\
81	0.113155916019526\\
82	0.113155916019526\\
83	0.113155916019526\\
84	0.113155916019526\\
85	0.113155916019526\\
86	0.113155916019526\\
87	0.113155916019526\\
88	0.113155916019526\\
89	0.113155916019526\\
90	0.113155916019526\\
91	0.113155916019526\\
92	0.113155916019526\\
93	0.113155916019526\\
94	0.113155916019526\\
95	0.113155916019526\\
96	0.113155916019526\\
97	0.113155916019526\\
98	0.113155916019526\\
99	0.113155916019526\\
100	0.113155916019526\\
101	0.113155916019526\\
102	0.113155916019526\\
103	0.113155916019526\\
104	0.113155916019526\\
105	0.113155916019526\\
106	0.113155916019526\\
107	0.113155916019526\\
108	0.113155916019526\\
109	0.113155916019526\\
110	0.113155916019526\\
111	0.113155916019526\\
112	0.113155916019526\\
113	0.113155916019526\\
114	0.113155916019526\\
115	0.113155916019526\\
116	0.113155916019526\\
117	0.113155916019526\\
118	0.113155916019526\\
119	0.113155916019526\\
120	0.113155916019526\\
121	0.113155916019526\\
122	0.113155916019526\\
123	0.113155916019526\\
124	0.113155916019526\\
125	0.113155916019526\\
126	0.113155916019526\\
127	0.113155916019526\\
128	0.113155916019526\\
129	0.113155916019526\\
130	0.113155916019526\\
131	0.113155916019526\\
132	0.113155916019526\\
133	0.113155916019526\\
134	0.113155916019526\\
135	0.113155916019526\\
136	0.113155916019526\\
137	0.113155916019526\\
138	0.113155916019526\\
139	0.113155916019526\\
140	0.113155916019526\\
141	0.113155916019526\\
142	0.113155916019526\\
143	0.113155916019526\\
144	0.113155916019526\\
145	0.113155916019526\\
146	0.113155916019526\\
147	0.113155916019526\\
148	0.113155916019526\\
149	0.113155916019526\\
150	0.113155916019526\\
};
\addplot [color=red,solid,line width=1.0pt,forget plot]
  table[row sep=crcr]{%
1	5.56486466797821\\
2	5.49281836476935\\
3	5.41074485544619\\
4	5.30673871437156\\
5	5.17432078744488\\
6	5.02111196498546\\
7	4.82622564644837\\
8	4.60372553069309\\
9	4.33628927172446\\
10	4.05150130826427\\
11	3.77448557652619\\
12	3.47355897421694\\
13	3.17409884906187\\
14	2.88832331382611\\
15	2.62029997132865\\
16	2.388670442313\\
17	2.14299825325026\\
18	1.91326225183252\\
19	1.69848356259456\\
20	1.51311726196888\\
21	1.33680118426444\\
22	1.19690838475932\\
23	1.07154369361423\\
24	0.967230114946257\\
25	0.865719417754956\\
26	0.777545661474472\\
27	0.703390495418179\\
28	0.638122491906222\\
29	0.575815476867019\\
30	0.527246234146809\\
31	0.48657226581841\\
32	0.444271100979372\\
33	0.415409071220444\\
34	0.378772544104088\\
35	0.348585708586037\\
36	0.324879653957361\\
37	0.301675993650789\\
38	0.284569797692798\\
39	0.265657558080315\\
40	0.253130820411364\\
41	0.243370161680951\\
42	0.236359535235171\\
43	0.227027904175533\\
44	0.217273486771339\\
45	0.21134747568449\\
46	0.201522632758971\\
47	0.192319407243411\\
48	0.183965439416692\\
49	0.173963720559607\\
50	0.171173263547926\\
51	0.166727644883613\\
52	0.158632873773217\\
53	0.155084747648953\\
54	0.153158281822173\\
55	0.148611099322851\\
56	0.146152026567753\\
57	0.144736400812435\\
58	0.141496621746153\\
59	0.138891208408814\\
60	0.138467764374558\\
61	0.137689576719501\\
62	0.137202169369006\\
63	0.136590558540364\\
64	0.136170468237957\\
65	0.135832821223145\\
66	0.135562927570976\\
67	0.135348413665653\\
68	0.13517896528185\\
69	0.135046051825209\\
70	0.134942659623679\\
71	0.134863048689282\\
72	0.133673020842766\\
73	0.133624671528987\\
74	0.132801039436121\\
75	0.131730497330421\\
76	0.131037028156377\\
77	0.131008177847067\\
78	0.128397817997351\\
79	0.127072199378657\\
80	0.125568310474804\\
81	0.123967958977751\\
82	0.123991136369989\\
83	0.121915047539655\\
84	0.12184628328674\\
85	0.118572481899656\\
86	0.118421437875272\\
87	0.118264619304285\\
88	0.118100678084575\\
89	0.117928930787933\\
90	0.117749403162062\\
91	0.117562803250796\\
92	0.117370436664141\\
93	0.117174079483187\\
94	0.116975825240582\\
95	0.116777922326107\\
96	0.116582616870152\\
97	0.116392013714332\\
98	0.116207964756246\\
99	0.116031990143327\\
100	0.115865233923374\\
101	0.115708452263837\\
102	0.115562029594467\\
103	0.115426016253155\\
104	0.11530018049745\\
105	0.115184067980551\\
106	0.115077062737849\\
107	0.114978445080751\\
108	0.114887443253643\\
109	0.11480327705108\\
110	0.114725192681536\\
111	0.114652488954016\\
112	0.114584535372108\\
113	0.11452078300152\\
114	0.114460769100087\\
115	0.114404116527561\\
116	0.114350528935899\\
117	0.11429978271046\\
118	0.114251716601401\\
119	0.114206219951367\\
120	0.114163220379791\\
121	0.114122671713014\\
122	0.114084542843624\\
123	0.114048808060112\\
124	0.114015439216615\\
125	0.113984399927327\\
126	0.113955641790366\\
127	0.113929102490441\\
128	0.113904705512842\\
129	0.113882361130314\\
130	0.113861968298967\\
131	0.11384341711304\\
132	0.113826591510526\\
133	0.113811371980832\\
134	0.11379763809093\\
135	0.11378527070939\\
136	0.113774153862632\\
137	0.113764176201771\\
138	0.113755232090773\\
139	0.113747222348192\\
140	0.11374005468712\\
141	0.113733643903214\\
142	0.113727911860787\\
143	0.113722787323645\\
144	0.113718205672187\\
145	0.113714108542151\\
146	0.113710443414163\\
147	0.113707163177359\\
148	0.113704225685051\\
149	0.113701593315846\\
150	0.113699232549849\\
};
\end{axis}
\end{tikzpicture}%
%
 }\hfil  
   \subfigure{%
    % This file was created by matlab2tikz v0.4.7 running on MATLAB 8.4.
% Copyright (c) 2008--2014, Nico Schlömer <nico.schloemer@gmail.com>
% All rights reserved.
% Minimal pgfplots version: 1.3
% 
% The latest updates can be retrieved from
%   http://www.mathworks.com/matlabcentral/fileexchange/22022-matlab2tikz
% where you can also make suggestions and rate matlab2tikz.
% 
\begin{tikzpicture}

\begin{axis}[%
width=2.8cm,
height=1.4cm,
scale only axis,
every outer x axis line/.append style={white!15!black},
every x tick label/.append style={font=\color{white!15!black}},
xmin=1,
xmax=150,
xlabel={$t$},
every outer y axis line/.append style={white!15!black},
every y tick label/.append style={font=\color{white!15!black}},
ymin=0,
ymax=6,
ylabel={g'},
title style={font=\bfseries},
title={$\sigma{=}0.3$, $n{=}30$, $d{=}3$, $\rho{=}0.5$},
axis x line*=bottom,
axis y line*=left
]
\addplot [color=blue,solid,line width=1.0pt,forget plot]
  table[row sep=crcr]{%
1	0.277498021378859\\
2	0.277498021378859\\
3	0.277498021378859\\
4	0.277498021378859\\
5	0.277498021378859\\
6	0.277498021378859\\
7	0.277498021378859\\
8	0.277498021378859\\
9	0.277498021378859\\
10	0.277498021378859\\
11	0.277498021378859\\
12	0.277498021378859\\
13	0.277498021378859\\
14	0.277498021378859\\
15	0.277498021378859\\
16	0.277498021378859\\
17	0.277498021378859\\
18	0.277498021378859\\
19	0.277498021378859\\
20	0.277498021378859\\
21	0.277498021378859\\
22	0.277498021378859\\
23	0.277498021378859\\
24	0.277498021378859\\
25	0.277498021378859\\
26	0.277498021378859\\
27	0.277498021378859\\
28	0.277498021378859\\
29	0.277498021378859\\
30	0.277498021378859\\
31	0.277498021378859\\
32	0.277498021378859\\
33	0.277498021378859\\
34	0.277498021378859\\
35	0.277498021378859\\
36	0.277498021378859\\
37	0.277498021378859\\
38	0.277498021378859\\
39	0.277498021378859\\
40	0.277498021378859\\
41	0.277498021378859\\
42	0.277498021378859\\
43	0.277498021378859\\
44	0.277498021378859\\
45	0.277498021378859\\
46	0.277498021378859\\
47	0.277498021378859\\
48	0.277498021378859\\
49	0.277498021378859\\
50	0.277498021378859\\
51	0.277498021378859\\
52	0.277498021378859\\
53	0.277498021378859\\
54	0.277498021378859\\
55	0.277498021378859\\
56	0.277498021378859\\
57	0.277498021378859\\
58	0.277498021378859\\
59	0.277498021378859\\
60	0.277498021378859\\
61	0.277498021378859\\
62	0.277498021378859\\
63	0.277498021378859\\
64	0.277498021378859\\
65	0.277498021378859\\
66	0.277498021378859\\
67	0.277498021378859\\
68	0.277498021378859\\
69	0.277498021378859\\
70	0.277498021378859\\
71	0.277498021378859\\
72	0.277498021378859\\
73	0.277498021378859\\
74	0.277498021378859\\
75	0.277498021378859\\
76	0.277498021378859\\
77	0.277498021378859\\
78	0.277498021378859\\
79	0.277498021378859\\
80	0.277498021378859\\
81	0.277498021378859\\
82	0.277498021378859\\
83	0.277498021378859\\
84	0.277498021378859\\
85	0.277498021378859\\
86	0.277498021378859\\
87	0.277498021378859\\
88	0.277498021378859\\
89	0.277498021378859\\
90	0.277498021378859\\
91	0.277498021378859\\
92	0.277498021378859\\
93	0.277498021378859\\
94	0.277498021378859\\
95	0.277498021378859\\
96	0.277498021378859\\
97	0.277498021378859\\
98	0.277498021378859\\
99	0.277498021378859\\
100	0.277498021378859\\
101	0.277498021378859\\
102	0.277498021378859\\
103	0.277498021378859\\
104	0.277498021378859\\
105	0.277498021378859\\
106	0.277498021378859\\
107	0.277498021378859\\
108	0.277498021378859\\
109	0.277498021378859\\
110	0.277498021378859\\
111	0.277498021378859\\
112	0.277498021378859\\
113	0.277498021378859\\
114	0.277498021378859\\
115	0.277498021378859\\
116	0.277498021378859\\
117	0.277498021378859\\
118	0.277498021378859\\
119	0.277498021378859\\
120	0.277498021378859\\
121	0.277498021378859\\
122	0.277498021378859\\
123	0.277498021378859\\
124	0.277498021378859\\
125	0.277498021378859\\
126	0.277498021378859\\
127	0.277498021378859\\
128	0.277498021378859\\
129	0.277498021378859\\
130	0.277498021378859\\
131	0.277498021378859\\
132	0.277498021378859\\
133	0.277498021378859\\
134	0.277498021378859\\
135	0.277498021378859\\
136	0.277498021378859\\
137	0.277498021378859\\
138	0.277498021378859\\
139	0.277498021378859\\
140	0.277498021378859\\
141	0.277498021378859\\
142	0.277498021378859\\
143	0.277498021378859\\
144	0.277498021378859\\
145	0.277498021378859\\
146	0.277498021378859\\
147	0.277498021378859\\
148	0.277498021378859\\
149	0.277498021378859\\
150	0.277498021378859\\
};
\addplot [color=red,solid,line width=1.0pt,forget plot]
  table[row sep=crcr]{%
1	5.55951592523087\\
2	5.48236765510905\\
3	5.4001722039316\\
4	5.29485244471415\\
5	5.14561839502678\\
6	4.97734638236091\\
7	4.79661832633217\\
8	4.58904513686634\\
9	4.36111060549875\\
10	4.10425971278653\\
11	3.847996792244\\
12	3.58614673682248\\
13	3.33275635188499\\
14	3.07325608279713\\
15	2.8079460308243\\
16	2.57425455615667\\
17	2.36829050943515\\
18	2.17253448545018\\
19	1.99036470600313\\
20	1.82177065139398\\
21	1.65747474730626\\
22	1.53906383954346\\
23	1.41461173494647\\
24	1.3104554658337\\
25	1.2079633370712\\
26	1.11443295539053\\
27	1.03270939748468\\
28	0.954601132373004\\
29	0.877243804210686\\
30	0.815462713957401\\
31	0.757468316767073\\
32	0.710238317644227\\
33	0.661397361973255\\
34	0.625122082720754\\
35	0.585492830362346\\
36	0.549142614383896\\
37	0.518517830067986\\
38	0.487486774250689\\
39	0.465355205671885\\
40	0.437628732999038\\
41	0.421647781519966\\
42	0.407631431851025\\
43	0.398983196440544\\
44	0.391963647609563\\
45	0.378383193204275\\
46	0.371920313541951\\
47	0.363692410708838\\
48	0.354163467580326\\
49	0.341824282321938\\
50	0.338875892530403\\
51	0.332333142594255\\
52	0.329720398344151\\
53	0.324720797935599\\
54	0.319050005309472\\
55	0.315171077927064\\
56	0.311716849952792\\
57	0.310294862047383\\
58	0.307193161548259\\
59	0.30608959388737\\
60	0.303488582187488\\
61	0.301444906297434\\
62	0.300750493170075\\
63	0.298379589045158\\
64	0.295392565434036\\
65	0.29145843286968\\
66	0.291035159482672\\
67	0.28904576693599\\
68	0.287138473335887\\
69	0.286870971874147\\
70	0.283498726855732\\
71	0.283294879698165\\
72	0.283118530867559\\
73	0.28296564817621\\
74	0.282832870862311\\
75	0.282717393561\\
76	0.281007237838929\\
77	0.280916814210627\\
78	0.28083796957134\\
79	0.28076926741654\\
80	0.280709468232715\\
81	0.280657493201566\\
82	0.280612395949922\\
83	0.280573340880927\\
84	0.280539586705485\\
85	0.280510473932174\\
86	0.280485415253684\\
87	0.280463887967552\\
88	0.280445427768412\\
89	0.280429623431653\\
90	0.280416112063758\\
91	0.28040457471775\\
92	0.280394732263123\\
93	0.28038634146172\\
94	0.28037919123933\\
95	0.280373099163116\\
96	0.280367908142656\\
97	0.280363483372016\\
98	0.280359709525435\\
99	0.280356488212444\\
100	0.28035373569124\\
101	0.2803513808329\\
102	0.280349363324036\\
103	0.280347632091911\\
104	0.280346143933775\\
105	0.280344862331081\\
106	0.280343756429066\\
107	0.280342800162725\\
108	0.280341971511188\\
109	0.28034125186389\\
110	0.280340625483377\\
111	0.280340079051188\\
112	0.280339601284782\\
113	0.280339182614966\\
114	0.28033881491464\\
115	0.280338491270963\\
116	0.280338205794112\\
117	0.280337953456859\\
118	0.280337729960021\\
119	0.280337531619618\\
120	0.280337355272206\\
121	0.280337198195441\\
122	0.280337058041357\\
123	0.280336932780306\\
124	0.280336820653783\\
125	0.280336720134705\\
126	0.280336629893908\\
127	0.280336548771853\\
128	0.280336475754696\\
129	0.280336409954013\\
130	0.280336350589587\\
131	0.280336296974774\\
132	0.280336248504033\\
133	0.280336204642279\\
134	0.280336164915768\\
135	0.28033612890429\\
136	0.280336096234452\\
137	0.280336066573896\\
138	0.280336039626308\\
139	0.280336015127096\\
140	0.280335992839645\\
141	0.280335972552063\\
142	0.280335954074343\\
143	0.28033593723589\\
144	0.28033592188335\\
145	0.280335907878715\\
146	0.280335895097653\\
147	0.28033588342804\\
148	0.280335872768663\\
149	0.280335863028078\\
150	0.280335854123592\\
};
\end{axis}
\end{tikzpicture}%
%
 }  \hfil
   \subfigure{%
    % This file was created by matlab2tikz v0.4.7 running on MATLAB 8.4.
% Copyright (c) 2008--2014, Nico Schlömer <nico.schloemer@gmail.com>
% All rights reserved.
% Minimal pgfplots version: 1.3
% 
% The latest updates can be retrieved from
%   http://www.mathworks.com/matlabcentral/fileexchange/22022-matlab2tikz
% where you can also make suggestions and rate matlab2tikz.
% 
\begin{tikzpicture}

\begin{axis}[%
width=2.8cm,
height=1.4cm,
scale only axis,
every outer x axis line/.append style={white!15!black},
every x tick label/.append style={font=\color{white!15!black}},
xmin=1,
xmax=150,
xlabel={$t$},
every outer y axis line/.append style={white!15!black},
every y tick label/.append style={font=\color{white!15!black}},
ymin=0,
ymax=6,
ylabel={g'},
title style={font=\bfseries},
title={$\sigma{=}0.4$, $n{=}30$, $d{=}3$, $\rho{=}0.5$},
axis x line*=bottom,
axis y line*=left
]
\addplot [color=blue,solid,line width=1.0pt,forget plot]
  table[row sep=crcr]{%
1	0.625376566061339\\
2	0.625376566061339\\
3	0.625376566061339\\
4	0.625376566061339\\
5	0.625376566061339\\
6	0.625376566061339\\
7	0.625376566061339\\
8	0.625376566061339\\
9	0.625376566061339\\
10	0.625376566061339\\
11	0.625376566061339\\
12	0.625376566061339\\
13	0.625376566061339\\
14	0.625376566061339\\
15	0.625376566061339\\
16	0.625376566061339\\
17	0.625376566061339\\
18	0.625376566061339\\
19	0.625376566061339\\
20	0.625376566061339\\
21	0.625376566061339\\
22	0.625376566061339\\
23	0.625376566061339\\
24	0.625376566061339\\
25	0.625376566061339\\
26	0.625376566061339\\
27	0.625376566061339\\
28	0.625376566061339\\
29	0.625376566061339\\
30	0.625376566061339\\
31	0.625376566061339\\
32	0.625376566061339\\
33	0.625376566061339\\
34	0.625376566061339\\
35	0.625376566061339\\
36	0.625376566061339\\
37	0.625376566061339\\
38	0.625376566061339\\
39	0.625376566061339\\
40	0.625376566061339\\
41	0.625376566061339\\
42	0.625376566061339\\
43	0.625376566061339\\
44	0.625376566061339\\
45	0.625376566061339\\
46	0.625376566061339\\
47	0.625376566061339\\
48	0.625376566061339\\
49	0.625376566061339\\
50	0.625376566061339\\
51	0.625376566061339\\
52	0.625376566061339\\
53	0.625376566061339\\
54	0.625376566061339\\
55	0.625376566061339\\
56	0.625376566061339\\
57	0.625376566061339\\
58	0.625376566061339\\
59	0.625376566061339\\
60	0.625376566061339\\
61	0.625376566061339\\
62	0.625376566061339\\
63	0.625376566061339\\
64	0.625376566061339\\
65	0.625376566061339\\
66	0.625376566061339\\
67	0.625376566061339\\
68	0.625376566061339\\
69	0.625376566061339\\
70	0.625376566061339\\
71	0.625376566061339\\
72	0.625376566061339\\
73	0.625376566061339\\
74	0.625376566061339\\
75	0.625376566061339\\
76	0.625376566061339\\
77	0.625376566061339\\
78	0.625376566061339\\
79	0.625376566061339\\
80	0.625376566061339\\
81	0.625376566061339\\
82	0.625376566061339\\
83	0.625376566061339\\
84	0.625376566061339\\
85	0.625376566061339\\
86	0.625376566061339\\
87	0.625376566061339\\
88	0.625376566061339\\
89	0.625376566061339\\
90	0.625376566061339\\
91	0.625376566061339\\
92	0.625376566061339\\
93	0.625376566061339\\
94	0.625376566061339\\
95	0.625376566061339\\
96	0.625376566061339\\
97	0.625376566061339\\
98	0.625376566061339\\
99	0.625376566061339\\
100	0.625376566061339\\
101	0.625376566061339\\
102	0.625376566061339\\
103	0.625376566061339\\
104	0.625376566061339\\
105	0.625376566061339\\
106	0.625376566061339\\
107	0.625376566061339\\
108	0.625376566061339\\
109	0.625376566061339\\
110	0.625376566061339\\
111	0.625376566061339\\
112	0.625376566061339\\
113	0.625376566061339\\
114	0.625376566061339\\
115	0.625376566061339\\
116	0.625376566061339\\
117	0.625376566061339\\
118	0.625376566061339\\
119	0.625376566061339\\
120	0.625376566061339\\
121	0.625376566061339\\
122	0.625376566061339\\
123	0.625376566061339\\
124	0.625376566061339\\
125	0.625376566061339\\
126	0.625376566061339\\
127	0.625376566061339\\
128	0.625376566061339\\
129	0.625376566061339\\
130	0.625376566061339\\
131	0.625376566061339\\
132	0.625376566061339\\
133	0.625376566061339\\
134	0.625376566061339\\
135	0.625376566061339\\
136	0.625376566061339\\
137	0.625376566061339\\
138	0.625376566061339\\
139	0.625376566061339\\
140	0.625376566061339\\
141	0.625376566061339\\
142	0.625376566061339\\
143	0.625376566061339\\
144	0.625376566061339\\
145	0.625376566061339\\
146	0.625376566061339\\
147	0.625376566061339\\
148	0.625376566061339\\
149	0.625376566061339\\
150	0.625376566061339\\
};
\addplot [color=red,solid,line width=1.0pt,forget plot]
  table[row sep=crcr]{%
1	5.56190338168981\\
2	5.49597496499845\\
3	5.42294146115717\\
4	5.33025141803094\\
5	5.22299006985183\\
6	5.08973657316527\\
7	4.93029549878713\\
8	4.74234764861593\\
9	4.54184881854023\\
10	4.33544672770566\\
11	4.12535190897209\\
12	3.90855597518046\\
13	3.68849600666151\\
14	3.48132357715984\\
15	3.27011787942873\\
16	3.0531420690133\\
17	2.84559324005195\\
18	2.65204131203338\\
19	2.47273695132625\\
20	2.31505242403416\\
21	2.16140082245805\\
22	2.03351742637096\\
23	1.91853004461907\\
24	1.81145464159781\\
25	1.68994461066785\\
26	1.6073848388439\\
27	1.51496748959603\\
28	1.44059925963929\\
29	1.36006739462833\\
30	1.2877553462009\\
31	1.23246029849058\\
32	1.17793335257927\\
33	1.12983890564205\\
34	1.08333623080797\\
35	1.0427774356416\\
36	1.01317752197835\\
37	0.981199710539067\\
38	0.956408771233148\\
39	0.937069272269914\\
40	0.913518511693523\\
41	0.896378212595832\\
42	0.881229667717595\\
43	0.860735047164886\\
44	0.84145366571021\\
45	0.827308180321682\\
46	0.80912575870814\\
47	0.798496849832254\\
48	0.784841962267468\\
49	0.777308525586277\\
50	0.763791637121529\\
51	0.75359794268501\\
52	0.742696162568231\\
53	0.735995322884769\\
54	0.727138954034996\\
55	0.718157299843715\\
56	0.712797684042866\\
57	0.708426909935912\\
58	0.698946393204632\\
59	0.693182900350673\\
60	0.686250102239742\\
61	0.682568187935069\\
62	0.673185390254945\\
63	0.665350369788913\\
64	0.664279656039098\\
65	0.663581923652666\\
66	0.66208374662203\\
67	0.659380146369092\\
68	0.658322014622705\\
69	0.657470697844879\\
70	0.657358157398199\\
71	0.657697801513889\\
72	0.657705873238929\\
73	0.657811158142602\\
74	0.658753979830554\\
75	0.658323277529844\\
76	0.657647703490717\\
77	0.656365383862393\\
78	0.655451354743461\\
79	0.655343413885519\\
80	0.654424109480942\\
81	0.65440566992433\\
82	0.653995464357142\\
83	0.652922164474588\\
84	0.652586548082772\\
85	0.652693264610817\\
86	0.648199166432283\\
87	0.644589884924551\\
88	0.643100087083723\\
89	0.640623514112528\\
90	0.639758805088018\\
91	0.637939202122297\\
92	0.63743348694965\\
93	0.637046209584092\\
94	0.636743392462586\\
95	0.636502031237968\\
96	0.636306392208659\\
97	0.63614550536384\\
98	0.636011543658223\\
99	0.635898791508778\\
100	0.635802985621486\\
101	0.63572088713221\\
102	0.635649997677529\\
103	0.635588366104725\\
104	0.635534453320273\\
105	0.635487035283516\\
106	0.635445131675178\\
107	0.6354079523252\\
108	0.635374856270879\\
109	0.635345320043062\\
110	0.635318912864176\\
111	0.635295277136472\\
112	0.635274113051519\\
113	0.635255166453342\\
114	0.635238219293233\\
115	0.635223082158568\\
116	0.635209588462631\\
117	0.635197589960925\\
118	0.635186953320583\\
119	0.635177557518433\\
120	0.635169291883741\\
121	0.635162054635513\\
122	0.635155751792946\\
123	0.63515029636206\\
124	0.635145607722126\\
125	0.635141611152822\\
126	0.635138237457391\\
127	0.635135422648712\\
128	0.635133107674595\\
129	0.635131238165907\\
130	0.635129764196796\\
131	0.635128640050476\\
132	0.635127823987094\\
133	0.635127278012291\\
134	0.635126967646535\\
135	0.635126861696094\\
136	0.635126932027031\\
137	0.63512715334376\\
138	0.635127502973658\\
139	0.635127960659154\\
140	0.635128508358413\\
141	0.635129130055585\\
142	0.635129811581263\\
143	0.635130540443602\\
144	0.635131305670309\\
145	0.635132097661542\\
146	0.635132908053581\\
147	0.635133729593008\\
148	0.635134556021045\\
149	0.635135381967601\\
150	0.635136202854527\\
};
\end{axis}
\end{tikzpicture}%
%
 }
    } \vspace{-0.4cm}
 	\centerline{
    \subfigure{%
    % This file was created by matlab2tikz v0.4.7 running on MATLAB 8.4.
% Copyright (c) 2008--2014, Nico Schlömer <nico.schloemer@gmail.com>
% All rights reserved.
% Minimal pgfplots version: 1.3
% 
% The latest updates can be retrieved from
%   http://www.mathworks.com/matlabcentral/fileexchange/22022-matlab2tikz
% where you can also make suggestions and rate matlab2tikz.
% 
\begin{tikzpicture}

\begin{axis}[%
width=2.8cm,
height=1.4cm,
scale only axis,
every outer x axis line/.append style={white!15!black},
every x tick label/.append style={font=\color{white!15!black}},
xmin=1,
xmax=150,
xlabel={$t$},
every outer y axis line/.append style={white!15!black},
every y tick label/.append style={font=\color{white!15!black}},
ymin=0,
ymax=8,
ylabel={g'},
title style={font=\bfseries},
title={$\sigma{=}0.3$, $n{=}30$, $d{=}4$, $\rho{=}0.5$},
axis x line*=bottom,
axis y line*=left
]
\addplot [color=blue,solid,line width=1.0pt,forget plot]
  table[row sep=crcr]{%
1	0.561417905112387\\
2	0.561417905112387\\
3	0.561417905112387\\
4	0.561417905112387\\
5	0.561417905112387\\
6	0.561417905112387\\
7	0.561417905112387\\
8	0.561417905112387\\
9	0.561417905112387\\
10	0.561417905112387\\
11	0.561417905112387\\
12	0.561417905112387\\
13	0.561417905112387\\
14	0.561417905112387\\
15	0.561417905112387\\
16	0.561417905112387\\
17	0.561417905112387\\
18	0.561417905112387\\
19	0.561417905112387\\
20	0.561417905112387\\
21	0.561417905112387\\
22	0.561417905112387\\
23	0.561417905112387\\
24	0.561417905112387\\
25	0.561417905112387\\
26	0.561417905112387\\
27	0.561417905112387\\
28	0.561417905112387\\
29	0.561417905112387\\
30	0.561417905112387\\
31	0.561417905112387\\
32	0.561417905112387\\
33	0.561417905112387\\
34	0.561417905112387\\
35	0.561417905112387\\
36	0.561417905112387\\
37	0.561417905112387\\
38	0.561417905112387\\
39	0.561417905112387\\
40	0.561417905112387\\
41	0.561417905112387\\
42	0.561417905112387\\
43	0.561417905112387\\
44	0.561417905112387\\
45	0.561417905112387\\
46	0.561417905112387\\
47	0.561417905112387\\
48	0.561417905112387\\
49	0.561417905112387\\
50	0.561417905112387\\
51	0.561417905112387\\
52	0.561417905112387\\
53	0.561417905112387\\
54	0.561417905112387\\
55	0.561417905112387\\
56	0.561417905112387\\
57	0.561417905112387\\
58	0.561417905112387\\
59	0.561417905112387\\
60	0.561417905112387\\
61	0.561417905112387\\
62	0.561417905112387\\
63	0.561417905112387\\
64	0.561417905112387\\
65	0.561417905112387\\
66	0.561417905112387\\
67	0.561417905112387\\
68	0.561417905112387\\
69	0.561417905112387\\
70	0.561417905112387\\
71	0.561417905112387\\
72	0.561417905112387\\
73	0.561417905112387\\
74	0.561417905112387\\
75	0.561417905112387\\
76	0.561417905112387\\
77	0.561417905112387\\
78	0.561417905112387\\
79	0.561417905112387\\
80	0.561417905112387\\
81	0.561417905112387\\
82	0.561417905112387\\
83	0.561417905112387\\
84	0.561417905112387\\
85	0.561417905112387\\
86	0.561417905112387\\
87	0.561417905112387\\
88	0.561417905112387\\
89	0.561417905112387\\
90	0.561417905112387\\
91	0.561417905112387\\
92	0.561417905112387\\
93	0.561417905112387\\
94	0.561417905112387\\
95	0.561417905112387\\
96	0.561417905112387\\
97	0.561417905112387\\
98	0.561417905112387\\
99	0.561417905112387\\
100	0.561417905112387\\
101	0.561417905112387\\
102	0.561417905112387\\
103	0.561417905112387\\
104	0.561417905112387\\
105	0.561417905112387\\
106	0.561417905112387\\
107	0.561417905112387\\
108	0.561417905112387\\
109	0.561417905112387\\
110	0.561417905112387\\
111	0.561417905112387\\
112	0.561417905112387\\
113	0.561417905112387\\
114	0.561417905112387\\
115	0.561417905112387\\
116	0.561417905112387\\
117	0.561417905112387\\
118	0.561417905112387\\
119	0.561417905112387\\
120	0.561417905112387\\
121	0.561417905112387\\
122	0.561417905112387\\
123	0.561417905112387\\
124	0.561417905112387\\
125	0.561417905112387\\
126	0.561417905112387\\
127	0.561417905112387\\
128	0.561417905112387\\
129	0.561417905112387\\
130	0.561417905112387\\
131	0.561417905112387\\
132	0.561417905112387\\
133	0.561417905112387\\
134	0.561417905112387\\
135	0.561417905112387\\
136	0.561417905112387\\
137	0.561417905112387\\
138	0.561417905112387\\
139	0.561417905112387\\
140	0.561417905112387\\
141	0.561417905112387\\
142	0.561417905112387\\
143	0.561417905112387\\
144	0.561417905112387\\
145	0.561417905112387\\
146	0.561417905112387\\
147	0.561417905112387\\
148	0.561417905112387\\
149	0.561417905112387\\
150	0.561417905112387\\
};
\addplot [color=red,solid,line width=1.0pt,forget plot]
  table[row sep=crcr]{%
1	7.42585125046701\\
2	7.33277546904133\\
3	7.21688915317926\\
4	7.07733547353242\\
5	6.91929563505208\\
6	6.71905375236111\\
7	6.49390995978419\\
8	6.22278286017869\\
9	5.92723952412915\\
10	5.60388369962632\\
11	5.24893093480852\\
12	4.91854052790194\\
13	4.56252055405656\\
14	4.2229847095203\\
15	3.88444251943176\\
16	3.57044956576086\\
17	3.28963760382661\\
18	3.01700459725206\\
19	2.78562290410996\\
20	2.57273919835678\\
21	2.36315620202318\\
22	2.17474745959768\\
23	2.01950124912455\\
24	1.87658336112471\\
25	1.75058949547828\\
26	1.64300110204078\\
27	1.54623875344129\\
28	1.45474841582301\\
29	1.37983278491704\\
30	1.3078736664979\\
31	1.24107831416062\\
32	1.17260743473768\\
33	1.1147680772361\\
34	1.05695748725528\\
35	1.01424189988988\\
36	0.978999974294062\\
37	0.943764123734543\\
38	0.916339034039106\\
39	0.889561861502456\\
40	0.860575790033879\\
41	0.843614423568856\\
42	0.821041055475897\\
43	0.8053065134999\\
44	0.785293475340662\\
45	0.764784104381348\\
46	0.750935090535409\\
47	0.735988896763008\\
48	0.724886446177701\\
49	0.709853976883005\\
50	0.699922356628886\\
51	0.690282826430883\\
52	0.679273700676945\\
53	0.666445899718126\\
54	0.655677743047344\\
55	0.648527389073558\\
56	0.639662612351815\\
57	0.629684634695341\\
58	0.625063555200703\\
59	0.619985815047144\\
60	0.615463000152407\\
61	0.612572368975299\\
62	0.608746462779932\\
63	0.600328706026716\\
64	0.597985015716771\\
65	0.594617895477867\\
66	0.591459672960972\\
67	0.587888161103746\\
68	0.586183638078739\\
69	0.583344465711293\\
70	0.581993060331445\\
71	0.57880968020838\\
72	0.577772963443241\\
73	0.576882258563136\\
74	0.576119219563643\\
75	0.575467591928648\\
76	0.574912783963488\\
77	0.574441687654342\\
78	0.571789649769991\\
79	0.571430605711146\\
80	0.571130733192995\\
81	0.570879828475526\\
82	0.570669685102099\\
83	0.570493592141124\\
84	0.570346006508561\\
85	0.570222321108795\\
86	0.570118690984309\\
87	0.570031897996421\\
88	0.569959243061103\\
89	0.569898459092422\\
90	0.569847639977913\\
91	0.569805182185402\\
92	0.569769736437664\\
93	0.569740167487375\\
94	0.569715520468179\\
95	0.569694992632674\\
96	0.569677909541646\\
97	0.569663704960499\\
98	0.569651903863794\\
99	0.569642108059073\\
100	0.569633984026217\\
101	0.569627252635126\\
102	0.569621680457742\\
103	0.569617072433778\\
104	0.56961326568541\\
105	0.569610124306414\\
106	0.569607534976801\\
107	0.569605403275972\\
108	0.569603650586118\\
109	0.56960221149377\\
110	0.56960103161124\\
111	0.569600065751541\\
112	0.569599276400587\\
113	0.56959863243913\\
114	0.569598108074293\\
115	0.569597681946888\\
116	0.569597336386014\\
117	0.569597056786986\\
118	0.569596831092475\\
119	0.569596649359947\\
120	0.569596503401255\\
121	0.569596386482492\\
122	0.569596293074175\\
123	0.569596218643423\\
124	0.569596159481191\\
125	0.569596112558707\\
126	0.569596075408277\\
127	0.569596046024377\\
128	0.569596022781642\\
129	0.569596004366927\\
130	0.569595989723062\\
131	0.569595978002351\\
132	0.569595968528141\\
133	0.569595960763126\\
134	0.569595954283218\\
135	0.569595948756039\\
136	0.56959594392326\\
137	0.5695959395861\\
138	0.569595935593463\\
139	0.56959593183224\\
140	0.569595928219408\\
141	0.569595924695605\\
142	0.569595921219923\\
143	0.569595917765704\\
144	0.569595914317158\\
145	0.569595910866653\\
146	0.569595907412559\\
147	0.569595903957538\\
148	0.569595900507201\\
149	0.569595897069065\\
150	0.569595893651746\\
};
\end{axis}
\end{tikzpicture}%
%
 } \hfil 
   \subfigure{%
    % This file was created by matlab2tikz v0.4.7 running on MATLAB 8.4.
% Copyright (c) 2008--2014, Nico Schlömer <nico.schloemer@gmail.com>
% All rights reserved.
% Minimal pgfplots version: 1.3
% 
% The latest updates can be retrieved from
%   http://www.mathworks.com/matlabcentral/fileexchange/22022-matlab2tikz
% where you can also make suggestions and rate matlab2tikz.
% 
\begin{tikzpicture}

\begin{axis}[%
width=2.8cm,
height=1.4cm,
scale only axis,
every outer x axis line/.append style={white!15!black},
every x tick label/.append style={font=\color{white!15!black}},
xmin=1,
xmax=150,
xlabel={$t$},
every outer y axis line/.append style={white!15!black},
every y tick label/.append style={font=\color{white!15!black}},
ymin=0,
ymax=10,
ylabel={g'},
title style={font=\bfseries},
title={$\sigma{=}0.3$, $n{=}30$, $d{=}5$, $\rho{=}0.5$},
axis x line*=bottom,
axis y line*=left
]
\addplot [color=blue,solid,line width=1.0pt,forget plot]
  table[row sep=crcr]{%
1	0.934833773404205\\
2	0.934833773404205\\
3	0.934833773404205\\
4	0.934833773404205\\
5	0.934833773404205\\
6	0.934833773404205\\
7	0.934833773404205\\
8	0.934833773404205\\
9	0.934833773404205\\
10	0.934833773404205\\
11	0.934833773404205\\
12	0.934833773404205\\
13	0.934833773404205\\
14	0.934833773404205\\
15	0.934833773404205\\
16	0.934833773404205\\
17	0.934833773404205\\
18	0.934833773404205\\
19	0.934833773404205\\
20	0.934833773404205\\
21	0.934833773404205\\
22	0.934833773404205\\
23	0.934833773404205\\
24	0.934833773404205\\
25	0.934833773404205\\
26	0.934833773404205\\
27	0.934833773404205\\
28	0.934833773404205\\
29	0.934833773404205\\
30	0.934833773404205\\
31	0.934833773404205\\
32	0.934833773404205\\
33	0.934833773404205\\
34	0.934833773404205\\
35	0.934833773404205\\
36	0.934833773404205\\
37	0.934833773404205\\
38	0.934833773404205\\
39	0.934833773404205\\
40	0.934833773404205\\
41	0.934833773404205\\
42	0.934833773404205\\
43	0.934833773404205\\
44	0.934833773404205\\
45	0.934833773404205\\
46	0.934833773404205\\
47	0.934833773404205\\
48	0.934833773404205\\
49	0.934833773404205\\
50	0.934833773404205\\
51	0.934833773404205\\
52	0.934833773404205\\
53	0.934833773404205\\
54	0.934833773404205\\
55	0.934833773404205\\
56	0.934833773404205\\
57	0.934833773404205\\
58	0.934833773404205\\
59	0.934833773404205\\
60	0.934833773404205\\
61	0.934833773404205\\
62	0.934833773404205\\
63	0.934833773404205\\
64	0.934833773404205\\
65	0.934833773404205\\
66	0.934833773404205\\
67	0.934833773404205\\
68	0.934833773404205\\
69	0.934833773404205\\
70	0.934833773404205\\
71	0.934833773404205\\
72	0.934833773404205\\
73	0.934833773404205\\
74	0.934833773404205\\
75	0.934833773404205\\
76	0.934833773404205\\
77	0.934833773404205\\
78	0.934833773404205\\
79	0.934833773404205\\
80	0.934833773404205\\
81	0.934833773404205\\
82	0.934833773404205\\
83	0.934833773404205\\
84	0.934833773404205\\
85	0.934833773404205\\
86	0.934833773404205\\
87	0.934833773404205\\
88	0.934833773404205\\
89	0.934833773404205\\
90	0.934833773404205\\
91	0.934833773404205\\
92	0.934833773404205\\
93	0.934833773404205\\
94	0.934833773404205\\
95	0.934833773404205\\
96	0.934833773404205\\
97	0.934833773404205\\
98	0.934833773404205\\
99	0.934833773404205\\
100	0.934833773404205\\
101	0.934833773404205\\
102	0.934833773404205\\
103	0.934833773404205\\
104	0.934833773404205\\
105	0.934833773404205\\
106	0.934833773404205\\
107	0.934833773404205\\
108	0.934833773404205\\
109	0.934833773404205\\
110	0.934833773404205\\
111	0.934833773404205\\
112	0.934833773404205\\
113	0.934833773404205\\
114	0.934833773404205\\
115	0.934833773404205\\
116	0.934833773404205\\
117	0.934833773404205\\
118	0.934833773404205\\
119	0.934833773404205\\
120	0.934833773404205\\
121	0.934833773404205\\
122	0.934833773404205\\
123	0.934833773404205\\
124	0.934833773404205\\
125	0.934833773404205\\
126	0.934833773404205\\
127	0.934833773404205\\
128	0.934833773404205\\
129	0.934833773404205\\
130	0.934833773404205\\
131	0.934833773404205\\
132	0.934833773404205\\
133	0.934833773404205\\
134	0.934833773404205\\
135	0.934833773404205\\
136	0.934833773404205\\
137	0.934833773404205\\
138	0.934833773404205\\
139	0.934833773404205\\
140	0.934833773404205\\
141	0.934833773404205\\
142	0.934833773404205\\
143	0.934833773404205\\
144	0.934833773404205\\
145	0.934833773404205\\
146	0.934833773404205\\
147	0.934833773404205\\
148	0.934833773404205\\
149	0.934833773404205\\
150	0.934833773404205\\
};
\addplot [color=red,solid,line width=1.0pt,forget plot]
  table[row sep=crcr]{%
1	9.26542693537142\\
2	9.14442528472567\\
3	9.00250422997947\\
4	8.83708814926363\\
5	8.63627237534677\\
6	8.39326487134915\\
7	8.11276162977925\\
8	7.79076699206541\\
9	7.43653827860177\\
10	7.05314073297955\\
11	6.67231102249975\\
12	6.27288548984872\\
13	5.89528710496354\\
14	5.48290216320546\\
15	5.10667663211481\\
16	4.76493252217426\\
17	4.42845630127087\\
18	4.12027225912118\\
19	3.84064055643507\\
20	3.57701982381606\\
21	3.34242822259015\\
22	3.12486741905225\\
23	2.91894055316011\\
24	2.731425614779\\
25	2.56172071057053\\
26	2.40136907332008\\
27	2.24502355902202\\
28	2.11027638686884\\
29	1.98958069797338\\
30	1.87891304118734\\
31	1.7818147530124\\
32	1.6965747084262\\
33	1.62342448220653\\
34	1.55625336060069\\
35	1.50486916738301\\
36	1.45670451783447\\
37	1.41098406837541\\
38	1.37080906047832\\
39	1.33200080605695\\
40	1.29972976912709\\
41	1.27123950767287\\
42	1.24268633925344\\
43	1.21993691384138\\
44	1.19797955300388\\
45	1.17981122495893\\
46	1.15816423907672\\
47	1.13903962556523\\
48	1.12050236057427\\
49	1.10611028374826\\
50	1.09039005606439\\
51	1.08182752436775\\
52	1.07112048419323\\
53	1.0610911467921\\
54	1.05104328102223\\
55	1.04244624560813\\
56	1.03621280370555\\
57	1.02985013909561\\
58	1.02351196240883\\
59	1.01875253244051\\
60	1.01462959940459\\
61	1.01006674040667\\
62	1.00560554809361\\
63	1.00202761271253\\
64	0.999763842273688\\
65	0.996665359999969\\
66	0.992027812173352\\
67	0.988993887996618\\
68	0.98680434081195\\
69	0.984724115798337\\
70	0.982097569292322\\
71	0.979874413013517\\
72	0.975861999389767\\
73	0.974899360565625\\
74	0.971602727703143\\
75	0.969425281864334\\
76	0.968621342229748\\
77	0.96642406826002\\
78	0.962759685988107\\
79	0.962040532602482\\
80	0.961384684314677\\
81	0.960786433515792\\
82	0.960240737332578\\
83	0.958622736874648\\
84	0.958049972936653\\
85	0.95618721251639\\
86	0.955678668980768\\
87	0.955217456487136\\
88	0.954801496066407\\
89	0.954428244036838\\
90	0.954094855249886\\
91	0.953798337592669\\
92	0.953535671263856\\
93	0.953303888616854\\
94	0.953100122445897\\
95	0.952921633830558\\
96	0.952765828834703\\
97	0.952630269850724\\
98	0.952512684206129\\
99	0.952410970590272\\
100	0.952323202924411\\
101	0.952247631158588\\
102	0.952182678749513\\
103	0.952126936956584\\
104	0.95207915641415\\
105	0.952038236628782\\
106	0.952003214109743\\
107	0.951973249801131\\
108	0.95194761638474\\
109	0.951925685897988\\
110	0.951906917985617\\
111	0.951890848991594\\
112	0.951877082005334\\
113	0.95186527790529\\
114	0.951855147391592\\
115	0.951846443964579\\
116	0.95183895778443\\
117	0.951832510335334\\
118	0.951826949813132\\
119	0.951822147155592\\
120	0.951817992638076\\
121	0.951814392962716\\
122	0.951811268775619\\
123	0.951808552553388\\
124	0.951806186806927\\
125	0.951804122556873\\
126	0.95180231804088\\
127	0.951800737618345\\
128	0.951799350842942\\
129	0.951798131677527\\
130	0.951797057829718\\
131	0.951796110189616\\
132	0.951795272353936\\
133	0.951794530223189\\
134	0.951793871660587\\
135	0.951793286203102\\
136	0.951792764816578\\
137	0.951792299688061\\
138	0.951791884049548\\
139	0.951791512028324\\
140	0.951791178519736\\
141	0.951790879078964\\
142	0.951790609828856\\
143	0.951790367381357\\
144	0.951790148770459\\
145	0.951789951394918\\
146	0.951789772969246\\
147	0.951789611481741\\
148	0.951789465158493\\
149	0.951789332432465\\
150	0.951789211916914\\
};
\end{axis}
\end{tikzpicture}%
%
 }\hfil  
   \subfigure{%
    % This file was created by matlab2tikz v0.4.7 running on MATLAB 8.4.
% Copyright (c) 2008--2014, Nico Schlömer <nico.schloemer@gmail.com>
% All rights reserved.
% Minimal pgfplots version: 1.3
% 
% The latest updates can be retrieved from
%   http://www.mathworks.com/matlabcentral/fileexchange/22022-matlab2tikz
% where you can also make suggestions and rate matlab2tikz.
% 
\begin{tikzpicture}

\begin{axis}[%
width=2.8cm,
height=1.4cm,
scale only axis,
every outer x axis line/.append style={white!15!black},
every x tick label/.append style={font=\color{white!15!black}},
xmin=1,
xmax=150,
xlabel={$t$},
every outer y axis line/.append style={white!15!black},
every y tick label/.append style={font=\color{white!15!black}},
ymin=0,
ymax=12,
ylabel={g'},
title style={font=\bfseries},
title={$\sigma{=}0.3$, $n{=}30$, $d{=}6$, $\rho{=}0.5$},
axis x line*=bottom,
axis y line*=left
]
\addplot [color=blue,solid,line width=1.0pt,forget plot]
  table[row sep=crcr]{%
1	1.40584910141926\\
2	1.40584910141926\\
3	1.40584910141926\\
4	1.40584910141926\\
5	1.40584910141926\\
6	1.40584910141926\\
7	1.40584910141926\\
8	1.40584910141926\\
9	1.40584910141926\\
10	1.40584910141926\\
11	1.40584910141926\\
12	1.40584910141926\\
13	1.40584910141926\\
14	1.40584910141926\\
15	1.40584910141926\\
16	1.40584910141926\\
17	1.40584910141926\\
18	1.40584910141926\\
19	1.40584910141926\\
20	1.40584910141926\\
21	1.40584910141926\\
22	1.40584910141926\\
23	1.40584910141926\\
24	1.40584910141926\\
25	1.40584910141926\\
26	1.40584910141926\\
27	1.40584910141926\\
28	1.40584910141926\\
29	1.40584910141926\\
30	1.40584910141926\\
31	1.40584910141926\\
32	1.40584910141926\\
33	1.40584910141926\\
34	1.40584910141926\\
35	1.40584910141926\\
36	1.40584910141926\\
37	1.40584910141926\\
38	1.40584910141926\\
39	1.40584910141926\\
40	1.40584910141926\\
41	1.40584910141926\\
42	1.40584910141926\\
43	1.40584910141926\\
44	1.40584910141926\\
45	1.40584910141926\\
46	1.40584910141926\\
47	1.40584910141926\\
48	1.40584910141926\\
49	1.40584910141926\\
50	1.40584910141926\\
51	1.40584910141926\\
52	1.40584910141926\\
53	1.40584910141926\\
54	1.40584910141926\\
55	1.40584910141926\\
56	1.40584910141926\\
57	1.40584910141926\\
58	1.40584910141926\\
59	1.40584910141926\\
60	1.40584910141926\\
61	1.40584910141926\\
62	1.40584910141926\\
63	1.40584910141926\\
64	1.40584910141926\\
65	1.40584910141926\\
66	1.40584910141926\\
67	1.40584910141926\\
68	1.40584910141926\\
69	1.40584910141926\\
70	1.40584910141926\\
71	1.40584910141926\\
72	1.40584910141926\\
73	1.40584910141926\\
74	1.40584910141926\\
75	1.40584910141926\\
76	1.40584910141926\\
77	1.40584910141926\\
78	1.40584910141926\\
79	1.40584910141926\\
80	1.40584910141926\\
81	1.40584910141926\\
82	1.40584910141926\\
83	1.40584910141926\\
84	1.40584910141926\\
85	1.40584910141926\\
86	1.40584910141926\\
87	1.40584910141926\\
88	1.40584910141926\\
89	1.40584910141926\\
90	1.40584910141926\\
91	1.40584910141926\\
92	1.40584910141926\\
93	1.40584910141926\\
94	1.40584910141926\\
95	1.40584910141926\\
96	1.40584910141926\\
97	1.40584910141926\\
98	1.40584910141926\\
99	1.40584910141926\\
100	1.40584910141926\\
101	1.40584910141926\\
102	1.40584910141926\\
103	1.40584910141926\\
104	1.40584910141926\\
105	1.40584910141926\\
106	1.40584910141926\\
107	1.40584910141926\\
108	1.40584910141926\\
109	1.40584910141926\\
110	1.40584910141926\\
111	1.40584910141926\\
112	1.40584910141926\\
113	1.40584910141926\\
114	1.40584910141926\\
115	1.40584910141926\\
116	1.40584910141926\\
117	1.40584910141926\\
118	1.40584910141926\\
119	1.40584910141926\\
120	1.40584910141926\\
121	1.40584910141926\\
122	1.40584910141926\\
123	1.40584910141926\\
124	1.40584910141926\\
125	1.40584910141926\\
126	1.40584910141926\\
127	1.40584910141926\\
128	1.40584910141926\\
129	1.40584910141926\\
130	1.40584910141926\\
131	1.40584910141926\\
132	1.40584910141926\\
133	1.40584910141926\\
134	1.40584910141926\\
135	1.40584910141926\\
136	1.40584910141926\\
137	1.40584910141926\\
138	1.40584910141926\\
139	1.40584910141926\\
140	1.40584910141926\\
141	1.40584910141926\\
142	1.40584910141926\\
143	1.40584910141926\\
144	1.40584910141926\\
145	1.40584910141926\\
146	1.40584910141926\\
147	1.40584910141926\\
148	1.40584910141926\\
149	1.40584910141926\\
150	1.40584910141926\\
};
\addplot [color=red,solid,line width=1.0pt,forget plot]
  table[row sep=crcr]{%
1	11.1435774265132\\
2	11.0125646420288\\
3	10.8455452861872\\
4	10.6613569084449\\
5	10.4432057484317\\
6	10.1827071676419\\
7	9.8861829835334\\
8	9.54649384526505\\
9	9.16629524618438\\
10	8.76523402195221\\
11	8.34262175736981\\
12	7.92552316335081\\
13	7.49358947718771\\
14	7.06669671409403\\
15	6.67184825850867\\
16	6.28631197374939\\
17	5.91603713052392\\
18	5.56929933693994\\
19	5.25094009624215\\
20	4.95368249591748\\
21	4.69250912746135\\
22	4.44161519493054\\
23	4.20355332425688\\
24	3.99520037460839\\
25	3.79837603239192\\
26	3.61581251679336\\
27	3.44596707875483\\
28	3.2765856052278\\
29	3.12894150252466\\
30	2.99343711849464\\
31	2.87258144471937\\
32	2.75660854648344\\
33	2.63966249930233\\
34	2.52763763312082\\
35	2.42967991363233\\
36	2.33662377973822\\
37	2.2546407019639\\
38	2.18334574056848\\
39	2.11131844168651\\
40	2.05121604405019\\
41	1.99173102267168\\
42	1.94377675074181\\
43	1.90338042765412\\
44	1.86797307722646\\
45	1.82944087196302\\
46	1.79998441300452\\
47	1.7770218602133\\
48	1.74947586345906\\
49	1.72558975964017\\
50	1.69889354740917\\
51	1.6772175733735\\
52	1.66259099671434\\
53	1.64146200434486\\
54	1.62785106584478\\
55	1.61158931916203\\
56	1.60156557506075\\
57	1.590468473665\\
58	1.57994400974471\\
59	1.56690280683623\\
60	1.55710592080477\\
61	1.54723154091294\\
62	1.53798304314387\\
63	1.52992335473939\\
64	1.523479433667\\
65	1.51520901634321\\
66	1.50738040225785\\
67	1.49884871559174\\
68	1.49119232222442\\
69	1.48593359163567\\
70	1.48036667927269\\
71	1.47674073487572\\
72	1.47351487279317\\
73	1.47066072790606\\
74	1.46745257307117\\
75	1.46341958975761\\
76	1.46136417184901\\
77	1.45817812285879\\
78	1.45576673989622\\
79	1.45218000670533\\
80	1.45059000801406\\
81	1.44920421293751\\
82	1.44799957985304\\
83	1.44695460925911\\
84	1.44434382632501\\
85	1.44354924144096\\
86	1.44286073755734\\
87	1.44226417549585\\
88	1.44174717921983\\
89	1.44129900202452\\
90	1.44091036584758\\
91	1.43904701997944\\
92	1.43872044393666\\
93	1.43843188979405\\
94	1.43817674370415\\
95	1.43795112017593\\
96	1.4377517278387\\
97	1.43757575075269\\
98	1.43742074682255\\
99	1.43728456435243\\
100	1.43716527691448\\
101	1.43706113558102\\
102	1.43697053643925\\
103	1.43689200044572\\
104	1.43682416227208\\
105	1.43676576488204\\
106	1.43671565706085\\
107	1.43667279182026\\
108	1.43663622434453\\
109	1.43660510879545\\
110	1.43657869378711\\
111	1.43655631666051\\
112	1.43653739685669\\
113	1.43652142874361\\
114	1.43650797423658\\
115	1.43649665549781\\
116	1.43648714793225\\
117	1.43647917362902\\
118	1.436472495338\\
119	1.43646691102415\\
120	1.43646224900611\\
121	1.43645836366129\\
122	1.43645513166336\\
123	1.4364524487091\\
124	1.43645022668714\\
125	1.43644839124064\\
126	1.43644687967713\\
127	1.43644563918197\\
128	1.43644462529537\\
129	1.43644380061701\\
130	1.43644313370626\\
131	1.43644259815001\\
132	1.43644217177338\\
133	1.43644183597237\\
134	1.43644157514998\\
135	1.43644137624013\\
136	1.43644122830598\\
137	1.4364411222013\\
138	1.43644105028495\\
139	1.43644100618039\\
140	1.43644098457312\\
141	1.43644098104013\\
142	1.43644099190638\\
143	1.43644101412401\\
144	1.43644104517076\\
145	1.43644108296446\\
146	1.43644112579127\\
147	1.43644117224525\\
148	1.43644122117764\\
149	1.43644127165425\\
150	1.43644132291964\\
};
\end{axis}
\end{tikzpicture}%
%
 }  \hfil
   \subfigure{%
    % This file was created by matlab2tikz v0.4.7 running on MATLAB 8.4.
% Copyright (c) 2008--2014, Nico Schlömer <nico.schloemer@gmail.com>
% All rights reserved.
% Minimal pgfplots version: 1.3
% 
% The latest updates can be retrieved from
%   http://www.mathworks.com/matlabcentral/fileexchange/22022-matlab2tikz
% where you can also make suggestions and rate matlab2tikz.
% 
\begin{tikzpicture}

\begin{axis}[%
width=2.8cm,
height=1.4cm,
scale only axis,
every outer x axis line/.append style={white!15!black},
every x tick label/.append style={font=\color{white!15!black}},
xmin=1,
xmax=150,
xlabel={$t$},
every outer y axis line/.append style={white!15!black},
every y tick label/.append style={font=\color{white!15!black}},
ymin=0,
ymax=14,
ylabel={g'},
title style={font=\bfseries},
title={$\sigma{=}0.3$, $n{=}30$, $d{=}7$, $\rho{=}0.5$},
axis x line*=bottom,
axis y line*=left
]
\addplot [color=blue,solid,line width=1.0pt,forget plot]
  table[row sep=crcr]{%
1	1.97732460807778\\
2	1.97732460807778\\
3	1.97732460807778\\
4	1.97732460807778\\
5	1.97732460807778\\
6	1.97732460807778\\
7	1.97732460807778\\
8	1.97732460807778\\
9	1.97732460807778\\
10	1.97732460807778\\
11	1.97732460807778\\
12	1.97732460807778\\
13	1.97732460807778\\
14	1.97732460807778\\
15	1.97732460807778\\
16	1.97732460807778\\
17	1.97732460807778\\
18	1.97732460807778\\
19	1.97732460807778\\
20	1.97732460807778\\
21	1.97732460807778\\
22	1.97732460807778\\
23	1.97732460807778\\
24	1.97732460807778\\
25	1.97732460807778\\
26	1.97732460807778\\
27	1.97732460807778\\
28	1.97732460807778\\
29	1.97732460807778\\
30	1.97732460807778\\
31	1.97732460807778\\
32	1.97732460807778\\
33	1.97732460807778\\
34	1.97732460807778\\
35	1.97732460807778\\
36	1.97732460807778\\
37	1.97732460807778\\
38	1.97732460807778\\
39	1.97732460807778\\
40	1.97732460807778\\
41	1.97732460807778\\
42	1.97732460807778\\
43	1.97732460807778\\
44	1.97732460807778\\
45	1.97732460807778\\
46	1.97732460807778\\
47	1.97732460807778\\
48	1.97732460807778\\
49	1.97732460807778\\
50	1.97732460807778\\
51	1.97732460807778\\
52	1.97732460807778\\
53	1.97732460807778\\
54	1.97732460807778\\
55	1.97732460807778\\
56	1.97732460807778\\
57	1.97732460807778\\
58	1.97732460807778\\
59	1.97732460807778\\
60	1.97732460807778\\
61	1.97732460807778\\
62	1.97732460807778\\
63	1.97732460807778\\
64	1.97732460807778\\
65	1.97732460807778\\
66	1.97732460807778\\
67	1.97732460807778\\
68	1.97732460807778\\
69	1.97732460807778\\
70	1.97732460807778\\
71	1.97732460807778\\
72	1.97732460807778\\
73	1.97732460807778\\
74	1.97732460807778\\
75	1.97732460807778\\
76	1.97732460807778\\
77	1.97732460807778\\
78	1.97732460807778\\
79	1.97732460807778\\
80	1.97732460807778\\
81	1.97732460807778\\
82	1.97732460807778\\
83	1.97732460807778\\
84	1.97732460807778\\
85	1.97732460807778\\
86	1.97732460807778\\
87	1.97732460807778\\
88	1.97732460807778\\
89	1.97732460807778\\
90	1.97732460807778\\
91	1.97732460807778\\
92	1.97732460807778\\
93	1.97732460807778\\
94	1.97732460807778\\
95	1.97732460807778\\
96	1.97732460807778\\
97	1.97732460807778\\
98	1.97732460807778\\
99	1.97732460807778\\
100	1.97732460807778\\
101	1.97732460807778\\
102	1.97732460807778\\
103	1.97732460807778\\
104	1.97732460807778\\
105	1.97732460807778\\
106	1.97732460807778\\
107	1.97732460807778\\
108	1.97732460807778\\
109	1.97732460807778\\
110	1.97732460807778\\
111	1.97732460807778\\
112	1.97732460807778\\
113	1.97732460807778\\
114	1.97732460807778\\
115	1.97732460807778\\
116	1.97732460807778\\
117	1.97732460807778\\
118	1.97732460807778\\
119	1.97732460807778\\
120	1.97732460807778\\
121	1.97732460807778\\
122	1.97732460807778\\
123	1.97732460807778\\
124	1.97732460807778\\
125	1.97732460807778\\
126	1.97732460807778\\
127	1.97732460807778\\
128	1.97732460807778\\
129	1.97732460807778\\
130	1.97732460807778\\
131	1.97732460807778\\
132	1.97732460807778\\
133	1.97732460807778\\
134	1.97732460807778\\
135	1.97732460807778\\
136	1.97732460807778\\
137	1.97732460807778\\
138	1.97732460807778\\
139	1.97732460807778\\
140	1.97732460807778\\
141	1.97732460807778\\
142	1.97732460807778\\
143	1.97732460807778\\
144	1.97732460807778\\
145	1.97732460807778\\
146	1.97732460807778\\
147	1.97732460807778\\
148	1.97732460807778\\
149	1.97732460807778\\
150	1.97732460807778\\
};
\addplot [color=red,solid,line width=1.0pt,forget plot]
  table[row sep=crcr]{%
1	12.9697534052573\\
2	12.801053291279\\
3	12.6097785930028\\
4	12.3889975358489\\
5	12.123183433624\\
6	11.8118107090538\\
7	11.4704242727672\\
8	11.084000118281\\
9	10.6529620700834\\
10	10.1879156521919\\
11	9.70973570187853\\
12	9.22695318160095\\
13	8.74998755656269\\
14	8.27695253064394\\
15	7.82040427062306\\
16	7.40242688191758\\
17	6.98232052129925\\
18	6.59107033621953\\
19	6.21954564977465\\
20	5.87266522808397\\
21	5.54794007948804\\
22	5.26442227164842\\
23	4.9973046359703\\
24	4.7510555447948\\
25	4.5328070509528\\
26	4.32808086849311\\
27	4.13541523089546\\
28	3.95409332915977\\
29	3.797560281467\\
30	3.65619459805211\\
31	3.51420219328422\\
32	3.38955743033921\\
33	3.27910416410582\\
34	3.18386852339577\\
35	3.09004011757521\\
36	3.01253549276997\\
37	2.93003813937569\\
38	2.85377549309503\\
39	2.78338639937977\\
40	2.71831750387104\\
41	2.65990504266183\\
42	2.60750717606678\\
43	2.56161558221781\\
44	2.51711151530734\\
45	2.47437822170942\\
46	2.44002061184306\\
47	2.40862352856257\\
48	2.37726765456163\\
49	2.34773747593928\\
50	2.32163222382521\\
51	2.2960250024923\\
52	2.27314608498415\\
53	2.25195293388608\\
54	2.23449572205391\\
55	2.2166073247005\\
56	2.2019552443059\\
57	2.19057581415263\\
58	2.17834488219623\\
59	2.16658819535628\\
60	2.15492402298252\\
61	2.14410398472108\\
62	2.1311958684439\\
63	2.12305267662383\\
64	2.11433371629691\\
65	2.10581468516949\\
66	2.10006648114209\\
67	2.09365846813366\\
68	2.08726630771365\\
69	2.08179580511205\\
70	2.07480424933386\\
71	2.07061183025348\\
72	2.06675279550332\\
73	2.06181174326652\\
74	2.05630354479576\\
75	2.05316072612338\\
76	2.05032068503951\\
77	2.04776772209791\\
78	2.04548390729616\\
79	2.04344971224495\\
80	2.04011421305119\\
81	2.03846228482449\\
82	2.03699918533824\\
83	2.03502404042491\\
84	2.0338439912841\\
85	2.03280154344191\\
86	2.03188197081628\\
87	2.03107181486937\\
88	2.03035885162012\\
89	2.02973203597383\\
90	2.02918143251408\\
91	2.0286981386729\\
92	2.02827420428985\\
93	2.02790255049583\\
94	2.02757689016175\\
95	2.02729165157692\\
96	2.0270419064897\\
97	2.02682330317004\\
98	2.02663200477427\\
99	2.0264646330219\\
100	2.02631821702304\\
101	2.02619014699928\\
102	2.02607813259558\\
103	2.02598016546282\\
104	2.02589448578819\\
105	2.02581955245447\\
106	2.02575401651681\\
107	2.02569669769579\\
108	2.02564656359775\\
109	2.02560271138819\\
110	2.02556435166078\\
111	2.02553079426265\\
112	2.02550143585596\\
113	2.0254757490154\\
114	2.02545327268035\\
115	2.02543360379921\\
116	2.02541639002091\\
117	2.0254013233049\\
118	2.02538813433569\\
119	2.02537658764172\\
120	2.02536647733027\\
121	2.02535762336094\\
122	2.02534986828988\\
123	2.02534307442532\\
124	2.02533712134236\\
125	2.02533190371165\\
126	2.02532732940219\\
127	2.02532331782356\\
128	2.0253197984772\\
129	2.02531670969042\\
130	2.02531399750982\\
131	2.02531161473412\\
132	2.02530952006865\\
133	2.02530767738615\\
134	2.02530605508048\\
135	2.02530462550139\\
136	2.02530336446032\\
137	2.02530225079801\\
138	2.02530126600638\\
139	2.02530039389764\\
140	2.02529962031484\\
141	2.02529893287853\\
142	2.02529832076509\\
143	2.02529777451269\\
144	2.0252972858514\\
145	2.02529684755457\\
146	2.02529645330857\\
147	2.02529609759885\\
148	2.02529577561013\\
149	2.02529548313893\\
150	2.0252952165171\\
};
\end{axis}
\end{tikzpicture}%
%
 }
    } \vspace{-0.4cm}
    \centerline{
    \subfigure{%
    \input{tikz/errorPlotDistributed-Zmatrix-orthogonal-k-3-d-3-sigma-0.3-graphDensity-0.5-nDraws-\nDraws.tikz}%
 } \hfil 
   \subfigure{%
    \input{tikz/errorPlotDistributed-Zmatrix-orthogonal-k-5-d-3-sigma-0.3-graphDensity-0.5-nDraws-\nDraws.tikz}%
 }\hfil  
   \subfigure{%
    % This file was created by matlab2tikz v0.4.7 running on MATLAB 8.4.
% Copyright (c) 2008--2014, Nico Schlömer <nico.schloemer@gmail.com>
% All rights reserved.
% Minimal pgfplots version: 1.3
% 
% The latest updates can be retrieved from
%   http://www.mathworks.com/matlabcentral/fileexchange/22022-matlab2tikz
% where you can also make suggestions and rate matlab2tikz.
% 
\begin{tikzpicture}

\begin{axis}[%
width=2.8cm,
height=1.4cm,
scale only axis,
every outer x axis line/.append style={white!15!black},
every x tick label/.append style={font=\color{white!15!black}},
xmin=1,
xmax=150,
xlabel={$t$},
every outer y axis line/.append style={white!15!black},
every y tick label/.append style={font=\color{white!15!black}},
ymin=0,
ymax=6,
ylabel={g'},
title style={font=\bfseries},
title={$\sigma{=}0.3$, $n{=}20$, $d{=}3$, $\rho{=}0.5$},
axis x line*=bottom,
axis y line*=left
]
\addplot [color=blue,solid,line width=1.0pt,forget plot]
  table[row sep=crcr]{%
1	0.271017121772222\\
2	0.271017121772222\\
3	0.271017121772222\\
4	0.271017121772222\\
5	0.271017121772222\\
6	0.271017121772222\\
7	0.271017121772222\\
8	0.271017121772222\\
9	0.271017121772222\\
10	0.271017121772222\\
11	0.271017121772222\\
12	0.271017121772222\\
13	0.271017121772222\\
14	0.271017121772222\\
15	0.271017121772222\\
16	0.271017121772222\\
17	0.271017121772222\\
18	0.271017121772222\\
19	0.271017121772222\\
20	0.271017121772222\\
21	0.271017121772222\\
22	0.271017121772222\\
23	0.271017121772222\\
24	0.271017121772222\\
25	0.271017121772222\\
26	0.271017121772222\\
27	0.271017121772222\\
28	0.271017121772222\\
29	0.271017121772222\\
30	0.271017121772222\\
31	0.271017121772222\\
32	0.271017121772222\\
33	0.271017121772222\\
34	0.271017121772222\\
35	0.271017121772222\\
36	0.271017121772222\\
37	0.271017121772222\\
38	0.271017121772222\\
39	0.271017121772222\\
40	0.271017121772222\\
41	0.271017121772222\\
42	0.271017121772222\\
43	0.271017121772222\\
44	0.271017121772222\\
45	0.271017121772222\\
46	0.271017121772222\\
47	0.271017121772222\\
48	0.271017121772222\\
49	0.271017121772222\\
50	0.271017121772222\\
51	0.271017121772222\\
52	0.271017121772222\\
53	0.271017121772222\\
54	0.271017121772222\\
55	0.271017121772222\\
56	0.271017121772222\\
57	0.271017121772222\\
58	0.271017121772222\\
59	0.271017121772222\\
60	0.271017121772222\\
61	0.271017121772222\\
62	0.271017121772222\\
63	0.271017121772222\\
64	0.271017121772222\\
65	0.271017121772222\\
66	0.271017121772222\\
67	0.271017121772222\\
68	0.271017121772222\\
69	0.271017121772222\\
70	0.271017121772222\\
71	0.271017121772222\\
72	0.271017121772222\\
73	0.271017121772222\\
74	0.271017121772222\\
75	0.271017121772222\\
76	0.271017121772222\\
77	0.271017121772222\\
78	0.271017121772222\\
79	0.271017121772222\\
80	0.271017121772222\\
81	0.271017121772222\\
82	0.271017121772222\\
83	0.271017121772222\\
84	0.271017121772222\\
85	0.271017121772222\\
86	0.271017121772222\\
87	0.271017121772222\\
88	0.271017121772222\\
89	0.271017121772222\\
90	0.271017121772222\\
91	0.271017121772222\\
92	0.271017121772222\\
93	0.271017121772222\\
94	0.271017121772222\\
95	0.271017121772222\\
96	0.271017121772222\\
97	0.271017121772222\\
98	0.271017121772222\\
99	0.271017121772222\\
100	0.271017121772222\\
101	0.271017121772222\\
102	0.271017121772222\\
103	0.271017121772222\\
104	0.271017121772222\\
105	0.271017121772222\\
106	0.271017121772222\\
107	0.271017121772222\\
108	0.271017121772222\\
109	0.271017121772222\\
110	0.271017121772222\\
111	0.271017121772222\\
112	0.271017121772222\\
113	0.271017121772222\\
114	0.271017121772222\\
115	0.271017121772222\\
116	0.271017121772222\\
117	0.271017121772222\\
118	0.271017121772222\\
119	0.271017121772222\\
120	0.271017121772222\\
121	0.271017121772222\\
122	0.271017121772222\\
123	0.271017121772222\\
124	0.271017121772222\\
125	0.271017121772222\\
126	0.271017121772222\\
127	0.271017121772222\\
128	0.271017121772222\\
129	0.271017121772222\\
130	0.271017121772222\\
131	0.271017121772222\\
132	0.271017121772222\\
133	0.271017121772222\\
134	0.271017121772222\\
135	0.271017121772222\\
136	0.271017121772222\\
137	0.271017121772222\\
138	0.271017121772222\\
139	0.271017121772222\\
140	0.271017121772222\\
141	0.271017121772222\\
142	0.271017121772222\\
143	0.271017121772222\\
144	0.271017121772222\\
145	0.271017121772222\\
146	0.271017121772222\\
147	0.271017121772222\\
148	0.271017121772222\\
149	0.271017121772222\\
150	0.271017121772222\\
};
\addplot [color=red,solid,line width=1.0pt,forget plot]
  table[row sep=crcr]{%
1	5.38033990319828\\
2	5.32383660951934\\
3	5.25960748470064\\
4	5.19744789495975\\
5	5.11293576138391\\
6	5.02672902063225\\
7	4.92021544260706\\
8	4.80760297925352\\
9	4.67760038270416\\
10	4.52518363101907\\
11	4.37208562631854\\
12	4.20360845599926\\
13	4.02016907041073\\
14	3.86555359419355\\
15	3.68199204090896\\
16	3.50946151677287\\
17	3.36722434571305\\
18	3.2000921657333\\
19	3.01973606011961\\
20	2.86455659222456\\
21	2.72658850896891\\
22	2.57699709820839\\
23	2.42583169880438\\
24	2.28876961164785\\
25	2.16516600721785\\
26	2.04659949373319\\
27	1.93964459900408\\
28	1.84483498645408\\
29	1.74610922043676\\
30	1.64381605269059\\
31	1.57210649688799\\
32	1.49030405390002\\
33	1.4183307353107\\
34	1.36810771455146\\
35	1.30819300933113\\
36	1.24105498278844\\
37	1.18839612849947\\
38	1.12505806846419\\
39	1.07601071376269\\
40	1.03241236933936\\
41	0.983721119271296\\
42	0.944180676511844\\
43	0.907435933715608\\
44	0.870013213334372\\
45	0.833056989447346\\
46	0.795605409150357\\
47	0.769347833701225\\
48	0.73133505698151\\
49	0.7073234689032\\
50	0.684000700424115\\
51	0.657859479324609\\
52	0.63697075500414\\
53	0.619486452551422\\
54	0.602987158358058\\
55	0.580837775504851\\
56	0.56568465979004\\
57	0.541153015945634\\
58	0.529121300662399\\
59	0.50927637115423\\
60	0.50002615210256\\
61	0.489674627315009\\
62	0.477963479563641\\
63	0.470453408554542\\
64	0.463633916956688\\
65	0.448614518695156\\
66	0.437268279705598\\
67	0.431233147827488\\
68	0.423997717815464\\
69	0.415608161859067\\
70	0.409140298433686\\
71	0.401684084194322\\
72	0.391304547652102\\
73	0.382627132963386\\
74	0.37522962225726\\
75	0.369301983452315\\
76	0.36223445113099\\
77	0.355804565724421\\
78	0.34875630358506\\
79	0.345144750093614\\
80	0.340460573368006\\
81	0.335360057913368\\
82	0.327579495174128\\
83	0.322104157880678\\
84	0.317827732796223\\
85	0.311836882092132\\
86	0.309669158595768\\
87	0.305337023851851\\
88	0.303466561529655\\
89	0.301830421956712\\
90	0.300395806362645\\
91	0.29913472678133\\
92	0.295060137936618\\
93	0.283950675513693\\
94	0.282817566719377\\
95	0.28183417520259\\
96	0.280978575316164\\
97	0.280232613794876\\
98	0.279581111717018\\
99	0.279011278366064\\
100	0.278512263945039\\
101	0.278074811492853\\
102	0.277690983031114\\
103	0.277353942288602\\
104	0.277057780766822\\
105	0.276797377058667\\
106	0.276568281766246\\
107	0.276366622270337\\
108	0.276189023071103\\
109	0.27603253852028\\
110	0.275894595571546\\
111	0.275772944756832\\
112	0.2756656180115\\
113	0.275570892268137\\
114	0.275487257953298\\
115	0.275413391679463\\
116	0.275348132543811\\
117	0.275290461538087\\
118	0.275239483647786\\
119	0.275194412279394\\
120	0.275154555704842\\
121	0.275119305255005\\
122	0.275088125030459\\
123	0.275060542929043\\
124	0.275036142816719\\
125	0.275014557691573\\
126	0.274995463710926\\
127	0.274978574968996\\
128	0.274963638927581\\
129	0.27495043241532\\
130	0.274938758122294\\
131	0.274928441526515\\
132	0.274919328197214\\
133	0.274911281427137\\
134	0.274904180152311\\
135	0.274897917123202\\
136	0.27489239729584\\
137	0.2748875364156\\
138	0.274883259769798\\
139	0.274879501088352\\
140	0.27487620157437\\
141	0.27487330904885\\
142	0.274870777195645\\
143	0.2748685648946\\
144	0.27486663563228\\
145	0.274864956980975\\
146	0.274863500137898\\
147	0.274862239517394\\
148	0.274861152389931\\
149	0.274860218562351\\
150	0.274859420094557\\
};
\end{axis}
\end{tikzpicture}%
%
 }  \hfil
   \subfigure{%
    % This file was created by matlab2tikz v0.4.7 running on MATLAB 8.4.
% Copyright (c) 2008--2014, Nico Schlömer <nico.schloemer@gmail.com>
% All rights reserved.
% Minimal pgfplots version: 1.3
% 
% The latest updates can be retrieved from
%   http://www.mathworks.com/matlabcentral/fileexchange/22022-matlab2tikz
% where you can also make suggestions and rate matlab2tikz.
% 
\begin{tikzpicture}

\begin{axis}[%
width=2.8cm,
height=1.4cm,
scale only axis,
every outer x axis line/.append style={white!15!black},
every x tick label/.append style={font=\color{white!15!black}},
xmin=1,
xmax=150,
xlabel={$t$},
every outer y axis line/.append style={white!15!black},
every y tick label/.append style={font=\color{white!15!black}},
ymin=0,
ymax=6,
ylabel={g'},
title style={font=\bfseries},
title={$\sigma{=}0.3$, $n{=}50$, $d{=}3$, $\rho{=}0.5$},
axis x line*=bottom,
axis y line*=left
]
\addplot [color=blue,solid,line width=1.0pt,forget plot]
  table[row sep=crcr]{%
1	0.287751524992962\\
2	0.287751524992962\\
3	0.287751524992962\\
4	0.287751524992962\\
5	0.287751524992962\\
6	0.287751524992962\\
7	0.287751524992962\\
8	0.287751524992962\\
9	0.287751524992962\\
10	0.287751524992962\\
11	0.287751524992962\\
12	0.287751524992962\\
13	0.287751524992962\\
14	0.287751524992962\\
15	0.287751524992962\\
16	0.287751524992962\\
17	0.287751524992962\\
18	0.287751524992962\\
19	0.287751524992962\\
20	0.287751524992962\\
21	0.287751524992962\\
22	0.287751524992962\\
23	0.287751524992962\\
24	0.287751524992962\\
25	0.287751524992962\\
26	0.287751524992962\\
27	0.287751524992962\\
28	0.287751524992962\\
29	0.287751524992962\\
30	0.287751524992962\\
31	0.287751524992962\\
32	0.287751524992962\\
33	0.287751524992962\\
34	0.287751524992962\\
35	0.287751524992962\\
36	0.287751524992962\\
37	0.287751524992962\\
38	0.287751524992962\\
39	0.287751524992962\\
40	0.287751524992962\\
41	0.287751524992962\\
42	0.287751524992962\\
43	0.287751524992962\\
44	0.287751524992962\\
45	0.287751524992962\\
46	0.287751524992962\\
47	0.287751524992962\\
48	0.287751524992962\\
49	0.287751524992962\\
50	0.287751524992962\\
51	0.287751524992962\\
52	0.287751524992962\\
53	0.287751524992962\\
54	0.287751524992962\\
55	0.287751524992962\\
56	0.287751524992962\\
57	0.287751524992962\\
58	0.287751524992962\\
59	0.287751524992962\\
60	0.287751524992962\\
61	0.287751524992962\\
62	0.287751524992962\\
63	0.287751524992962\\
64	0.287751524992962\\
65	0.287751524992962\\
66	0.287751524992962\\
67	0.287751524992962\\
68	0.287751524992962\\
69	0.287751524992962\\
70	0.287751524992962\\
71	0.287751524992962\\
72	0.287751524992962\\
73	0.287751524992962\\
74	0.287751524992962\\
75	0.287751524992962\\
76	0.287751524992962\\
77	0.287751524992962\\
78	0.287751524992962\\
79	0.287751524992962\\
80	0.287751524992962\\
81	0.287751524992962\\
82	0.287751524992962\\
83	0.287751524992962\\
84	0.287751524992962\\
85	0.287751524992962\\
86	0.287751524992962\\
87	0.287751524992962\\
88	0.287751524992962\\
89	0.287751524992962\\
90	0.287751524992962\\
91	0.287751524992962\\
92	0.287751524992962\\
93	0.287751524992962\\
94	0.287751524992962\\
95	0.287751524992962\\
96	0.287751524992962\\
97	0.287751524992962\\
98	0.287751524992962\\
99	0.287751524992962\\
100	0.287751524992962\\
101	0.287751524992962\\
102	0.287751524992962\\
103	0.287751524992962\\
104	0.287751524992962\\
105	0.287751524992962\\
106	0.287751524992962\\
107	0.287751524992962\\
108	0.287751524992962\\
109	0.287751524992962\\
110	0.287751524992962\\
111	0.287751524992962\\
112	0.287751524992962\\
113	0.287751524992962\\
114	0.287751524992962\\
115	0.287751524992962\\
116	0.287751524992962\\
117	0.287751524992962\\
118	0.287751524992962\\
119	0.287751524992962\\
120	0.287751524992962\\
121	0.287751524992962\\
122	0.287751524992962\\
123	0.287751524992962\\
124	0.287751524992962\\
125	0.287751524992962\\
126	0.287751524992962\\
127	0.287751524992962\\
128	0.287751524992962\\
129	0.287751524992962\\
130	0.287751524992962\\
131	0.287751524992962\\
132	0.287751524992962\\
133	0.287751524992962\\
134	0.287751524992962\\
135	0.287751524992962\\
136	0.287751524992962\\
137	0.287751524992962\\
138	0.287751524992962\\
139	0.287751524992962\\
140	0.287751524992962\\
141	0.287751524992962\\
142	0.287751524992962\\
143	0.287751524992962\\
144	0.287751524992962\\
145	0.287751524992962\\
146	0.287751524992962\\
147	0.287751524992962\\
148	0.287751524992962\\
149	0.287751524992962\\
150	0.287751524992962\\
};
\addplot [color=red,solid,line width=1.0pt,forget plot]
  table[row sep=crcr]{%
1	5.68842635506432\\
2	5.59028694353382\\
3	5.43758302531509\\
4	5.21717145094578\\
5	4.91320113667002\\
6	4.53287030195325\\
7	4.058605818379\\
8	3.5580411355134\\
9	3.08117529609842\\
10	2.61635792316603\\
11	2.21600992169918\\
12	1.87824650979157\\
13	1.59375172032277\\
14	1.34265631111993\\
15	1.13932573892435\\
16	0.989056340440984\\
17	0.858933092552339\\
18	0.750277328970379\\
19	0.674838495251883\\
20	0.600462445134837\\
21	0.53469304915301\\
22	0.492950592922609\\
23	0.453341182622389\\
24	0.419802346471457\\
25	0.396356465986755\\
26	0.375638901811273\\
27	0.361035826759243\\
28	0.348513962138556\\
29	0.340690271288017\\
30	0.334819396079789\\
31	0.328813781813036\\
32	0.325872808726686\\
33	0.323626580532924\\
34	0.321215578733288\\
35	0.31981398684308\\
36	0.318762025949887\\
37	0.317491058510502\\
38	0.316814501233813\\
39	0.316399387773092\\
40	0.31584017396251\\
41	0.315427413346887\\
42	0.314174986019924\\
43	0.313499667444413\\
44	0.312708717487024\\
45	0.311700355448673\\
46	0.311079646725623\\
47	0.309885489642419\\
48	0.308150644374479\\
49	0.306328285687395\\
50	0.304115551999065\\
51	0.301541032376092\\
52	0.299672763330045\\
53	0.296967784948133\\
54	0.295147919833712\\
55	0.292440824157268\\
56	0.292018135332177\\
57	0.291709024888839\\
58	0.291481715367765\\
59	0.290102740517531\\
60	0.289963923100462\\
61	0.289862740218841\\
62	0.289788481333871\\
63	0.289733833142938\\
64	0.289693610338074\\
65	0.289664051041773\\
66	0.289642387772767\\
67	0.289626568154447\\
68	0.289615064068343\\
69	0.289606736565822\\
70	0.289600737753289\\
71	0.28959643818741\\
72	0.289593372451806\\
73	0.289591198046151\\
74	0.289589664231996\\
75	0.289588588436328\\
76	0.289587838438455\\
77	0.289587318992089\\
78	0.289586961840882\\
79	0.289586718316455\\
80	0.289586553887949\\
81	0.289586444174756\\
82	0.289586372047619\\
83	0.289586325533109\\
84	0.289586296306895\\
85	0.289586278615689\\
86	0.289586268509448\\
87	0.289586263296919\\
88	0.28958626116116\\
89	0.289586260889144\\
90	0.289586261682398\\
91	0.289586263024942\\
92	0.289586264591639\\
93	0.289586266184903\\
94	0.289586267691247\\
95	0.289586269051665\\
96	0.289586270241597\\
97	0.289586271257523\\
98	0.289586272108105\\
99	0.289586272808425\\
100	0.289586273376336\\
101	0.289586273830213\\
102	0.289586274187638\\
103	0.289586274464697\\
104	0.289586274675662\\
105	0.289586274832913\\
106	0.289586274947008\\
107	0.28958627502682\\
108	0.289586275079722\\
109	0.289586275111772\\
110	0.289586275127905\\
111	0.289586275132101\\
112	0.289586275127541\\
113	0.289586275116745\\
114	0.289586275101689\\
115	0.289586275083906\\
116	0.289586275064567\\
117	0.289586275044558\\
118	0.289586275024533\\
119	0.289586275004966\\
120	0.289586274986187\\
121	0.289586274968414\\
122	0.289586274951783\\
123	0.289586274936363\\
124	0.289586274922176\\
125	0.289586274909209\\
126	0.289586274897422\\
127	0.289586274886762\\
128	0.289586274877161\\
129	0.289586274868546\\
130	0.289586274860843\\
131	0.289586274853976\\
132	0.28958627484787\\
133	0.289586274842455\\
134	0.289586274837663\\
135	0.289586274833431\\
136	0.289586274829701\\
137	0.289586274826418\\
138	0.289586274823534\\
139	0.289586274821005\\
140	0.289586274818788\\
141	0.28958627481685\\
142	0.289586274815155\\
143	0.289586274813676\\
144	0.289586274812387\\
145	0.289586274811263\\
146	0.289586274810285\\
147	0.289586274809435\\
148	0.289586274808696\\
149	0.289586274808055\\
150	0.289586274807499\\
};
\end{axis}
\end{tikzpicture}%
%
 }
    } \vspace{-0.4cm}
    \centerline{
    \subfigure{%
    \input{tikz/errorPlotDistributed-Zmatrix-orthogonal-k-30-d-3-sigma-0.3-graphDensity-0.1-nDraws-\nDraws.tikz}%
 } \hfil 
   \subfigure{%
    % This file was created by matlab2tikz v0.4.7 running on MATLAB 8.4.
% Copyright (c) 2008--2014, Nico Schlömer <nico.schloemer@gmail.com>
% All rights reserved.
% Minimal pgfplots version: 1.3
% 
% The latest updates can be retrieved from
%   http://www.mathworks.com/matlabcentral/fileexchange/22022-matlab2tikz
% where you can also make suggestions and rate matlab2tikz.
% 
\begin{tikzpicture}

\begin{axis}[%
width=2.8cm,
height=1.4cm,
scale only axis,
every outer x axis line/.append style={white!15!black},
every x tick label/.append style={font=\color{white!15!black}},
xmin=1,
xmax=150,
xlabel={$t$},
every outer y axis line/.append style={white!15!black},
every y tick label/.append style={font=\color{white!15!black}},
ymin=0,
ymax=6,
ylabel={g'},
title style={font=\bfseries},
title={$\sigma{=}0.3$, $n{=}30$, $d{=}3$, $\rho{=}0.3$},
axis x line*=bottom,
axis y line*=left
]
\addplot [color=blue,solid,line width=1.0pt,forget plot]
  table[row sep=crcr]{%
1	0.270174248244598\\
2	0.270174248244598\\
3	0.270174248244598\\
4	0.270174248244598\\
5	0.270174248244598\\
6	0.270174248244598\\
7	0.270174248244598\\
8	0.270174248244598\\
9	0.270174248244598\\
10	0.270174248244598\\
11	0.270174248244598\\
12	0.270174248244598\\
13	0.270174248244598\\
14	0.270174248244598\\
15	0.270174248244598\\
16	0.270174248244598\\
17	0.270174248244598\\
18	0.270174248244598\\
19	0.270174248244598\\
20	0.270174248244598\\
21	0.270174248244598\\
22	0.270174248244598\\
23	0.270174248244598\\
24	0.270174248244598\\
25	0.270174248244598\\
26	0.270174248244598\\
27	0.270174248244598\\
28	0.270174248244598\\
29	0.270174248244598\\
30	0.270174248244598\\
31	0.270174248244598\\
32	0.270174248244598\\
33	0.270174248244598\\
34	0.270174248244598\\
35	0.270174248244598\\
36	0.270174248244598\\
37	0.270174248244598\\
38	0.270174248244598\\
39	0.270174248244598\\
40	0.270174248244598\\
41	0.270174248244598\\
42	0.270174248244598\\
43	0.270174248244598\\
44	0.270174248244598\\
45	0.270174248244598\\
46	0.270174248244598\\
47	0.270174248244598\\
48	0.270174248244598\\
49	0.270174248244598\\
50	0.270174248244598\\
51	0.270174248244598\\
52	0.270174248244598\\
53	0.270174248244598\\
54	0.270174248244598\\
55	0.270174248244598\\
56	0.270174248244598\\
57	0.270174248244598\\
58	0.270174248244598\\
59	0.270174248244598\\
60	0.270174248244598\\
61	0.270174248244598\\
62	0.270174248244598\\
63	0.270174248244598\\
64	0.270174248244598\\
65	0.270174248244598\\
66	0.270174248244598\\
67	0.270174248244598\\
68	0.270174248244598\\
69	0.270174248244598\\
70	0.270174248244598\\
71	0.270174248244598\\
72	0.270174248244598\\
73	0.270174248244598\\
74	0.270174248244598\\
75	0.270174248244598\\
76	0.270174248244598\\
77	0.270174248244598\\
78	0.270174248244598\\
79	0.270174248244598\\
80	0.270174248244598\\
81	0.270174248244598\\
82	0.270174248244598\\
83	0.270174248244598\\
84	0.270174248244598\\
85	0.270174248244598\\
86	0.270174248244598\\
87	0.270174248244598\\
88	0.270174248244598\\
89	0.270174248244598\\
90	0.270174248244598\\
91	0.270174248244598\\
92	0.270174248244598\\
93	0.270174248244598\\
94	0.270174248244598\\
95	0.270174248244598\\
96	0.270174248244598\\
97	0.270174248244598\\
98	0.270174248244598\\
99	0.270174248244598\\
100	0.270174248244598\\
101	0.270174248244598\\
102	0.270174248244598\\
103	0.270174248244598\\
104	0.270174248244598\\
105	0.270174248244598\\
106	0.270174248244598\\
107	0.270174248244598\\
108	0.270174248244598\\
109	0.270174248244598\\
110	0.270174248244598\\
111	0.270174248244598\\
112	0.270174248244598\\
113	0.270174248244598\\
114	0.270174248244598\\
115	0.270174248244598\\
116	0.270174248244598\\
117	0.270174248244598\\
118	0.270174248244598\\
119	0.270174248244598\\
120	0.270174248244598\\
121	0.270174248244598\\
122	0.270174248244598\\
123	0.270174248244598\\
124	0.270174248244598\\
125	0.270174248244598\\
126	0.270174248244598\\
127	0.270174248244598\\
128	0.270174248244598\\
129	0.270174248244598\\
130	0.270174248244598\\
131	0.270174248244598\\
132	0.270174248244598\\
133	0.270174248244598\\
134	0.270174248244598\\
135	0.270174248244598\\
136	0.270174248244598\\
137	0.270174248244598\\
138	0.270174248244598\\
139	0.270174248244598\\
140	0.270174248244598\\
141	0.270174248244598\\
142	0.270174248244598\\
143	0.270174248244598\\
144	0.270174248244598\\
145	0.270174248244598\\
146	0.270174248244598\\
147	0.270174248244598\\
148	0.270174248244598\\
149	0.270174248244598\\
150	0.270174248244598\\
};
\addplot [color=red,solid,line width=1.0pt,forget plot]
  table[row sep=crcr]{%
1	5.38106195864361\\
2	5.33677403540601\\
3	5.28325233252759\\
4	5.21916747989999\\
5	5.14879144094243\\
6	5.06619401671537\\
7	4.97428094277714\\
8	4.87934143682242\\
9	4.78030805567107\\
10	4.67609803862094\\
11	4.56968499376717\\
12	4.45269362451682\\
13	4.32967455469005\\
14	4.21503880038793\\
15	4.09511000635873\\
16	3.96948658425453\\
17	3.83457426640438\\
18	3.69792618747749\\
19	3.56938211994184\\
20	3.43596261069282\\
21	3.30954643217349\\
22	3.17885325471874\\
23	3.06296108410112\\
24	2.9610650133403\\
25	2.83295594525296\\
26	2.71580990480064\\
27	2.60867502578715\\
28	2.515496687711\\
29	2.42213824592022\\
30	2.31797344793244\\
31	2.24178621061688\\
32	2.15148483206595\\
33	2.06521000549572\\
34	1.98126399174266\\
35	1.89777813223717\\
36	1.8219075626629\\
37	1.75350480104535\\
38	1.68923347315641\\
39	1.63052678074482\\
40	1.56612638087982\\
41	1.50917821363666\\
42	1.44774441804611\\
43	1.39397433920624\\
44	1.33962456080565\\
45	1.28820641978592\\
46	1.24089067345445\\
47	1.20162310783929\\
48	1.15607349033398\\
49	1.10592346927181\\
50	1.06614700790069\\
51	1.02570234538126\\
52	0.99068488418499\\
53	0.954543280859768\\
54	0.920301435168707\\
55	0.885644525417739\\
56	0.854981133999205\\
57	0.832260571285296\\
58	0.807714679076578\\
59	0.779680114319711\\
60	0.75167631506763\\
61	0.727769031132774\\
62	0.70596042652128\\
63	0.681904137620394\\
64	0.661024320411424\\
65	0.63866899210019\\
66	0.616384978674369\\
67	0.598160745700214\\
68	0.581311135050225\\
69	0.56590621168156\\
70	0.555134049811644\\
71	0.542733769085435\\
72	0.527387810168452\\
73	0.507831790305676\\
74	0.496743629900073\\
75	0.487008111787831\\
76	0.479106232380739\\
77	0.470747419315342\\
78	0.459621539187853\\
79	0.452680651735739\\
80	0.443365759293574\\
81	0.437522709167873\\
82	0.423174305583647\\
83	0.417394104877944\\
84	0.409713680849248\\
85	0.404350676443003\\
86	0.400619338064282\\
87	0.394458177609765\\
88	0.390863659689779\\
89	0.385896195622773\\
90	0.38053513913892\\
91	0.376080244697859\\
92	0.372358618127015\\
93	0.368100826675711\\
94	0.361251156114462\\
95	0.357737513193606\\
96	0.353805765311438\\
97	0.351119204839039\\
98	0.348608106686309\\
99	0.346263786822143\\
100	0.344077218110893\\
101	0.34124304040203\\
102	0.33931577186215\\
103	0.335733103508897\\
104	0.332693906491039\\
105	0.331064283290379\\
106	0.329546964053244\\
107	0.328133718342426\\
108	0.326816769715731\\
109	0.324165758644278\\
110	0.321288689207471\\
111	0.320108712314519\\
112	0.319011455481933\\
113	0.317990191835279\\
114	0.317038677930704\\
115	0.316151149127437\\
116	0.314209493163485\\
117	0.313433520289645\\
118	0.312706998842284\\
119	0.312025852525436\\
120	0.311386358223692\\
121	0.310785116845403\\
122	0.310219025839889\\
123	0.309685253820149\\
124	0.309181217463028\\
125	0.308704560678706\\
126	0.30809116406231\\
127	0.307637821892488\\
128	0.306971251411883\\
129	0.3065540468742\\
130	0.306154408173248\\
131	0.305123277437196\\
132	0.304740419095328\\
133	0.303014948404661\\
134	0.302633323339783\\
135	0.302263688875688\\
136	0.301905464006406\\
137	0.301558164903504\\
138	0.30122139597062\\
139	0.300894840951345\\
140	0.300578254103034\\
141	0.300271451471914\\
142	0.299196390451222\\
143	0.297340789980786\\
144	0.297016326984962\\
145	0.296700964024305\\
146	0.296394864302921\\
147	0.296098211380327\\
148	0.295811201232635\\
149	0.295534034860071\\
150	0.295266911545916\\
};
\end{axis}
\end{tikzpicture}%
%
 }\hfil  
   \subfigure{%
    % This file was created by matlab2tikz v0.4.7 running on MATLAB 8.4.
% Copyright (c) 2008--2014, Nico Schlömer <nico.schloemer@gmail.com>
% All rights reserved.
% Minimal pgfplots version: 1.3
% 
% The latest updates can be retrieved from
%   http://www.mathworks.com/matlabcentral/fileexchange/22022-matlab2tikz
% where you can also make suggestions and rate matlab2tikz.
% 
\begin{tikzpicture}

\begin{axis}[%
width=2.8cm,
height=1.4cm,
scale only axis,
every outer x axis line/.append style={white!15!black},
every x tick label/.append style={font=\color{white!15!black}},
xmin=1,
xmax=150,
xlabel={$t$},
every outer y axis line/.append style={white!15!black},
every y tick label/.append style={font=\color{white!15!black}},
ymin=0,
ymax=6,
ylabel={g'},
title style={font=\bfseries},
title={$\sigma{=}0.3$, $n{=}30$, $d{=}3$, $\rho{=}0.8$},
axis x line*=bottom,
axis y line*=left
]
\addplot [color=blue,solid,line width=1.0pt,forget plot]
  table[row sep=crcr]{%
1	0.288251488996135\\
2	0.288251488996135\\
3	0.288251488996135\\
4	0.288251488996135\\
5	0.288251488996135\\
6	0.288251488996135\\
7	0.288251488996135\\
8	0.288251488996135\\
9	0.288251488996135\\
10	0.288251488996135\\
11	0.288251488996135\\
12	0.288251488996135\\
13	0.288251488996135\\
14	0.288251488996135\\
15	0.288251488996135\\
16	0.288251488996135\\
17	0.288251488996135\\
18	0.288251488996135\\
19	0.288251488996135\\
20	0.288251488996135\\
21	0.288251488996135\\
22	0.288251488996135\\
23	0.288251488996135\\
24	0.288251488996135\\
25	0.288251488996135\\
26	0.288251488996135\\
27	0.288251488996135\\
28	0.288251488996135\\
29	0.288251488996135\\
30	0.288251488996135\\
31	0.288251488996135\\
32	0.288251488996135\\
33	0.288251488996135\\
34	0.288251488996135\\
35	0.288251488996135\\
36	0.288251488996135\\
37	0.288251488996135\\
38	0.288251488996135\\
39	0.288251488996135\\
40	0.288251488996135\\
41	0.288251488996135\\
42	0.288251488996135\\
43	0.288251488996135\\
44	0.288251488996135\\
45	0.288251488996135\\
46	0.288251488996135\\
47	0.288251488996135\\
48	0.288251488996135\\
49	0.288251488996135\\
50	0.288251488996135\\
51	0.288251488996135\\
52	0.288251488996135\\
53	0.288251488996135\\
54	0.288251488996135\\
55	0.288251488996135\\
56	0.288251488996135\\
57	0.288251488996135\\
58	0.288251488996135\\
59	0.288251488996135\\
60	0.288251488996135\\
61	0.288251488996135\\
62	0.288251488996135\\
63	0.288251488996135\\
64	0.288251488996135\\
65	0.288251488996135\\
66	0.288251488996135\\
67	0.288251488996135\\
68	0.288251488996135\\
69	0.288251488996135\\
70	0.288251488996135\\
71	0.288251488996135\\
72	0.288251488996135\\
73	0.288251488996135\\
74	0.288251488996135\\
75	0.288251488996135\\
76	0.288251488996135\\
77	0.288251488996135\\
78	0.288251488996135\\
79	0.288251488996135\\
80	0.288251488996135\\
81	0.288251488996135\\
82	0.288251488996135\\
83	0.288251488996135\\
84	0.288251488996135\\
85	0.288251488996135\\
86	0.288251488996135\\
87	0.288251488996135\\
88	0.288251488996135\\
89	0.288251488996135\\
90	0.288251488996135\\
91	0.288251488996135\\
92	0.288251488996135\\
93	0.288251488996135\\
94	0.288251488996135\\
95	0.288251488996135\\
96	0.288251488996135\\
97	0.288251488996135\\
98	0.288251488996135\\
99	0.288251488996135\\
100	0.288251488996135\\
101	0.288251488996135\\
102	0.288251488996135\\
103	0.288251488996135\\
104	0.288251488996135\\
105	0.288251488996135\\
106	0.288251488996135\\
107	0.288251488996135\\
108	0.288251488996135\\
109	0.288251488996135\\
110	0.288251488996135\\
111	0.288251488996135\\
112	0.288251488996135\\
113	0.288251488996135\\
114	0.288251488996135\\
115	0.288251488996135\\
116	0.288251488996135\\
117	0.288251488996135\\
118	0.288251488996135\\
119	0.288251488996135\\
120	0.288251488996135\\
121	0.288251488996135\\
122	0.288251488996135\\
123	0.288251488996135\\
124	0.288251488996135\\
125	0.288251488996135\\
126	0.288251488996135\\
127	0.288251488996135\\
128	0.288251488996135\\
129	0.288251488996135\\
130	0.288251488996135\\
131	0.288251488996135\\
132	0.288251488996135\\
133	0.288251488996135\\
134	0.288251488996135\\
135	0.288251488996135\\
136	0.288251488996135\\
137	0.288251488996135\\
138	0.288251488996135\\
139	0.288251488996135\\
140	0.288251488996135\\
141	0.288251488996135\\
142	0.288251488996135\\
143	0.288251488996135\\
144	0.288251488996135\\
145	0.288251488996135\\
146	0.288251488996135\\
147	0.288251488996135\\
148	0.288251488996135\\
149	0.288251488996135\\
150	0.288251488996135\\
};
\addplot [color=red,solid,line width=1.0pt,forget plot]
  table[row sep=crcr]{%
1	5.66536002816156\\
2	5.53601513551905\\
3	5.36089853239269\\
4	5.11964667400376\\
5	4.78260878846668\\
6	4.3546692193112\\
7	3.88212811232705\\
8	3.34638274465837\\
9	2.86588691339996\\
10	2.43486544339228\\
11	2.0585674607167\\
12	1.75193253895995\\
13	1.4490679703755\\
14	1.22262080538557\\
15	1.03084524738137\\
16	0.888894377157814\\
17	0.77251593097118\\
18	0.689771783280086\\
19	0.606969515765285\\
20	0.533229156870088\\
21	0.482427890868659\\
22	0.446869081112844\\
23	0.415972890081584\\
24	0.390944576411874\\
25	0.375949047465195\\
26	0.361007875733084\\
27	0.349066638011329\\
28	0.338909950321669\\
29	0.333913820883163\\
30	0.324887794617485\\
31	0.323165076400475\\
32	0.318376130940299\\
33	0.316557062905263\\
34	0.313913199041344\\
35	0.310046453189181\\
36	0.304115084513313\\
37	0.300570646120344\\
38	0.298925992272191\\
39	0.295874247938139\\
40	0.294864512684819\\
41	0.294085827859744\\
42	0.292070063280934\\
43	0.291428567306887\\
44	0.290975431815592\\
45	0.290673537027392\\
46	0.290479828343637\\
47	0.290357936392144\\
48	0.290281846712342\\
49	0.290234455605788\\
50	0.290204947090996\\
51	0.290186583557208\\
52	0.290175181600698\\
53	0.290168138384641\\
54	0.29016382797035\\
55	0.290161230365995\\
56	0.290159703432636\\
57	0.290158841954195\\
58	0.290158390010682\\
59	0.29015818614018\\
60	0.290158128815062\\
61	0.290158154595149\\
62	0.290158224245219\\
63	0.290158313885968\\
64	0.29015840934225\\
65	0.290158502530623\\
66	0.290158589151808\\
67	0.290158667219966\\
68	0.290158736129307\\
69	0.290158796065765\\
70	0.290158847640037\\
71	0.290158891662264\\
72	0.290158929006935\\
73	0.29015896053489\\
74	0.290158987051056\\
75	0.290159009284233\\
76	0.290159027880162\\
77	0.29015904340232\\
78	0.29015905633695\\
79	0.290159067100141\\
80	0.290159076045679\\
81	0.290159083472872\\
82	0.290159089633933\\
83	0.290159094740717\\
84	0.290159098970716\\
85	0.290159102472316\\
86	0.290159105369347\\
87	0.290159107764994\\
88	0.290159109745129\\
89	0.290159111381134\\
90	0.290159112732293\\
91	0.290159113847802\\
92	0.29015911476845\\
93	0.290159115528037\\
94	0.290159116154554\\
95	0.290159116671169\\
96	0.290159117097048\\
97	0.29015911744804\\
98	0.290159117737243\\
99	0.290159117975481\\
100	0.290159118171692\\
101	0.290159118333256\\
102	0.290159118466265\\
103	0.290159118575745\\
104	0.290159118665841\\
105	0.290159118739973\\
106	0.290159118800958\\
107	0.29015911885112\\
108	0.290159118892374\\
109	0.290159118926295\\
110	0.290159118954184\\
111	0.29015911897711\\
112	0.290159118995953\\
113	0.290159119011439\\
114	0.290159119024165\\
115	0.29015911903462\\
116	0.290159119043209\\
117	0.290159119050264\\
118	0.290159119056059\\
119	0.290159119060818\\
120	0.290159119064726\\
121	0.290159119067935\\
122	0.29015911907057\\
123	0.290159119072733\\
124	0.290159119074508\\
125	0.290159119075966\\
126	0.290159119077162\\
127	0.290159119078143\\
128	0.290159119078949\\
129	0.29015911907961\\
130	0.290159119080152\\
131	0.290159119080597\\
132	0.290159119080962\\
133	0.290159119081262\\
134	0.290159119081507\\
135	0.290159119081709\\
136	0.290159119081874\\
137	0.29015911908201\\
138	0.290159119082121\\
139	0.290159119082212\\
140	0.290159119082287\\
141	0.290159119082348\\
142	0.290159119082399\\
143	0.29015911908244\\
144	0.290159119082474\\
145	0.290159119082501\\
146	0.290159119082524\\
147	0.290159119082543\\
148	0.290159119082558\\
149	0.290159119082571\\
150	0.290159119082581\\
};
\end{axis}
\end{tikzpicture}%
%
 }  \hfil
   \subfigure{%
    % This file was created by matlab2tikz v0.4.7 running on MATLAB 8.4.
% Copyright (c) 2008--2014, Nico Schlömer <nico.schloemer@gmail.com>
% All rights reserved.
% Minimal pgfplots version: 1.3
% 
% The latest updates can be retrieved from
%   http://www.mathworks.com/matlabcentral/fileexchange/22022-matlab2tikz
% where you can also make suggestions and rate matlab2tikz.
% 
\begin{tikzpicture}

\begin{axis}[%
width=2.8cm,
height=1.4cm,
scale only axis,
every outer x axis line/.append style={white!15!black},
every x tick label/.append style={font=\color{white!15!black}},
xmin=1,
xmax=150,
xlabel={$t$},
every outer y axis line/.append style={white!15!black},
every y tick label/.append style={font=\color{white!15!black}},
ymin=0,
ymax=6,
ylabel={g'},
title style={font=\bfseries},
title={$\sigma{=}0.3$, $n{=}30$, $d{=}3$, $\rho{=}1$},
axis x line*=bottom,
axis y line*=left
]
\addplot [color=blue,solid,line width=1.0pt,forget plot]
  table[row sep=crcr]{%
1	0.288326318530242\\
2	0.288326318530242\\
3	0.288326318530242\\
4	0.288326318530242\\
5	0.288326318530242\\
6	0.288326318530242\\
7	0.288326318530242\\
8	0.288326318530242\\
9	0.288326318530242\\
10	0.288326318530242\\
11	0.288326318530242\\
12	0.288326318530242\\
13	0.288326318530242\\
14	0.288326318530242\\
15	0.288326318530242\\
16	0.288326318530242\\
17	0.288326318530242\\
18	0.288326318530242\\
19	0.288326318530242\\
20	0.288326318530242\\
21	0.288326318530242\\
22	0.288326318530242\\
23	0.288326318530242\\
24	0.288326318530242\\
25	0.288326318530242\\
26	0.288326318530242\\
27	0.288326318530242\\
28	0.288326318530242\\
29	0.288326318530242\\
30	0.288326318530242\\
31	0.288326318530242\\
32	0.288326318530242\\
33	0.288326318530242\\
34	0.288326318530242\\
35	0.288326318530242\\
36	0.288326318530242\\
37	0.288326318530242\\
38	0.288326318530242\\
39	0.288326318530242\\
40	0.288326318530242\\
41	0.288326318530242\\
42	0.288326318530242\\
43	0.288326318530242\\
44	0.288326318530242\\
45	0.288326318530242\\
46	0.288326318530242\\
47	0.288326318530242\\
48	0.288326318530242\\
49	0.288326318530242\\
50	0.288326318530242\\
51	0.288326318530242\\
52	0.288326318530242\\
53	0.288326318530242\\
54	0.288326318530242\\
55	0.288326318530242\\
56	0.288326318530242\\
57	0.288326318530242\\
58	0.288326318530242\\
59	0.288326318530242\\
60	0.288326318530242\\
61	0.288326318530242\\
62	0.288326318530242\\
63	0.288326318530242\\
64	0.288326318530242\\
65	0.288326318530242\\
66	0.288326318530242\\
67	0.288326318530242\\
68	0.288326318530242\\
69	0.288326318530242\\
70	0.288326318530242\\
71	0.288326318530242\\
72	0.288326318530242\\
73	0.288326318530242\\
74	0.288326318530242\\
75	0.288326318530242\\
76	0.288326318530242\\
77	0.288326318530242\\
78	0.288326318530242\\
79	0.288326318530242\\
80	0.288326318530242\\
81	0.288326318530242\\
82	0.288326318530242\\
83	0.288326318530242\\
84	0.288326318530242\\
85	0.288326318530242\\
86	0.288326318530242\\
87	0.288326318530242\\
88	0.288326318530242\\
89	0.288326318530242\\
90	0.288326318530242\\
91	0.288326318530242\\
92	0.288326318530242\\
93	0.288326318530242\\
94	0.288326318530242\\
95	0.288326318530242\\
96	0.288326318530242\\
97	0.288326318530242\\
98	0.288326318530242\\
99	0.288326318530242\\
100	0.288326318530242\\
101	0.288326318530242\\
102	0.288326318530242\\
103	0.288326318530242\\
104	0.288326318530242\\
105	0.288326318530242\\
106	0.288326318530242\\
107	0.288326318530242\\
108	0.288326318530242\\
109	0.288326318530242\\
110	0.288326318530242\\
111	0.288326318530242\\
112	0.288326318530242\\
113	0.288326318530242\\
114	0.288326318530242\\
115	0.288326318530242\\
116	0.288326318530242\\
117	0.288326318530242\\
118	0.288326318530242\\
119	0.288326318530242\\
120	0.288326318530242\\
121	0.288326318530242\\
122	0.288326318530242\\
123	0.288326318530242\\
124	0.288326318530242\\
125	0.288326318530242\\
126	0.288326318530242\\
127	0.288326318530242\\
128	0.288326318530242\\
129	0.288326318530242\\
130	0.288326318530242\\
131	0.288326318530242\\
132	0.288326318530242\\
133	0.288326318530242\\
134	0.288326318530242\\
135	0.288326318530242\\
136	0.288326318530242\\
137	0.288326318530242\\
138	0.288326318530242\\
139	0.288326318530242\\
140	0.288326318530242\\
141	0.288326318530242\\
142	0.288326318530242\\
143	0.288326318530242\\
144	0.288326318530242\\
145	0.288326318530242\\
146	0.288326318530242\\
147	0.288326318530242\\
148	0.288326318530242\\
149	0.288326318530242\\
150	0.288326318530242\\
};
\addplot [color=red,solid,line width=1.0pt,forget plot]
  table[row sep=crcr]{%
1	5.69062358180133\\
2	5.50028223974668\\
3	5.18846003651223\\
4	4.73278993520499\\
5	4.13509356927212\\
6	3.37484022431805\\
7	2.6555081622708\\
8	2.07070180740914\\
9	1.62271577805569\\
10	1.27968507156132\\
11	1.03777976007103\\
12	0.824355704364511\\
13	0.660736477053972\\
14	0.533493568135404\\
15	0.433781043096146\\
16	0.372877136466629\\
17	0.33601176466071\\
18	0.315778045955694\\
19	0.302404605538254\\
20	0.296214202209603\\
21	0.293037612851714\\
22	0.291426578320648\\
23	0.290610798134594\\
24	0.290196766968948\\
25	0.289985926925836\\
26	0.289878189949841\\
27	0.289822961183825\\
28	0.289794567018856\\
29	0.28977993041222\\
30	0.289772367453981\\
31	0.289768451276908\\
32	0.289766419944952\\
33	0.289765365149909\\
34	0.289764817440811\\
35	0.289764533567614\\
36	0.289764387166366\\
37	0.289764312427841\\
38	0.28976427499826\\
39	0.289764256909393\\
40	0.289764248753173\\
41	0.28976424560633\\
42	0.289764244901273\\
43	0.289764245310004\\
44	0.289764246157035\\
45	0.289764247109541\\
46	0.289764248013403\\
47	0.289764248806317\\
48	0.289764249471803\\
49	0.289764250014996\\
50	0.289764250450056\\
51	0.289764250793799\\
52	0.289764251062625\\
53	0.289764251271185\\
54	0.289764251431949\\
55	0.289764251555207\\
56	0.289764251649283\\
57	0.289764251720805\\
58	0.289764251774994\\
59	0.289764251815926\\
60	0.289764251846759\\
61	0.289764251869928\\
62	0.289764251887299\\
63	0.289764251900295\\
64	0.28976425191\\
65	0.289764251917234\\
66	0.289764251922618\\
67	0.289764251926618\\
68	0.289764251929586\\
69	0.289764251931786\\
70	0.289764251933413\\
71	0.289764251934616\\
72	0.289764251935504\\
73	0.289764251936159\\
74	0.289764251936641\\
75	0.289764251936996\\
76	0.289764251937258\\
77	0.28976425193745\\
78	0.28976425193759\\
79	0.289764251937693\\
80	0.289764251937769\\
81	0.289764251937824\\
82	0.289764251937865\\
83	0.289764251937895\\
84	0.289764251937916\\
85	0.289764251937932\\
86	0.289764251937944\\
87	0.289764251937952\\
88	0.289764251937958\\
89	0.289764251937963\\
90	0.289764251937966\\
91	0.289764251937968\\
92	0.28976425193797\\
93	0.289764251937972\\
94	0.289764251937972\\
95	0.289764251937973\\
96	0.289764251937973\\
97	0.289764251937974\\
98	0.289764251937974\\
99	0.289764251937974\\
100	0.289764251937974\\
101	0.289764251937974\\
102	0.289764251937974\\
103	0.289764251937974\\
104	0.289764251937975\\
105	0.289764251937975\\
106	0.289764251937975\\
107	0.289764251937975\\
108	0.289764251937975\\
109	0.289764251937975\\
110	0.289764251937975\\
111	0.289764251937975\\
112	0.289764251937975\\
113	0.289764251937975\\
114	0.289764251937975\\
115	0.289764251937975\\
116	0.289764251937975\\
117	0.289764251937975\\
118	0.289764251937975\\
119	0.289764251937975\\
120	0.289764251937975\\
121	0.289764251937975\\
122	0.289764251937975\\
123	0.289764251937975\\
124	0.289764251937975\\
125	0.289764251937975\\
126	0.289764251937975\\
127	0.289764251937975\\
128	0.289764251937975\\
129	0.289764251937975\\
130	0.289764251937975\\
131	0.289764251937975\\
132	0.289764251937975\\
133	0.289764251937975\\
134	0.289764251937975\\
135	0.289764251937975\\
136	0.289764251937975\\
137	0.289764251937975\\
138	0.289764251937975\\
139	0.289764251937975\\
140	0.289764251937975\\
141	0.289764251937975\\
142	0.289764251937975\\
143	0.289764251937975\\
144	0.289764251937975\\
145	0.289764251937975\\
146	0.289764251937975\\
147	0.289764251937975\\
148	0.289764251937975\\
149	0.289764251937975\\
150	0.289764251937975\\
};
\end{axis}
\end{tikzpicture}%
%
 }
    } \vspace{-0.4cm}
\caption{Normalised error in \eqref{gPrimeFcn} on the vertical axis for the Z-matrix method (blue) and its distributed version (red) when considering transformations in $O(d)$. The horizontal axis shows the number of steps, where the step size has been chosen as $\epsilon = 0.01$ in each sub-figure. 
}
    \label{distributedZOrthogonal} 
\end{figure*}

\begin{figure*}
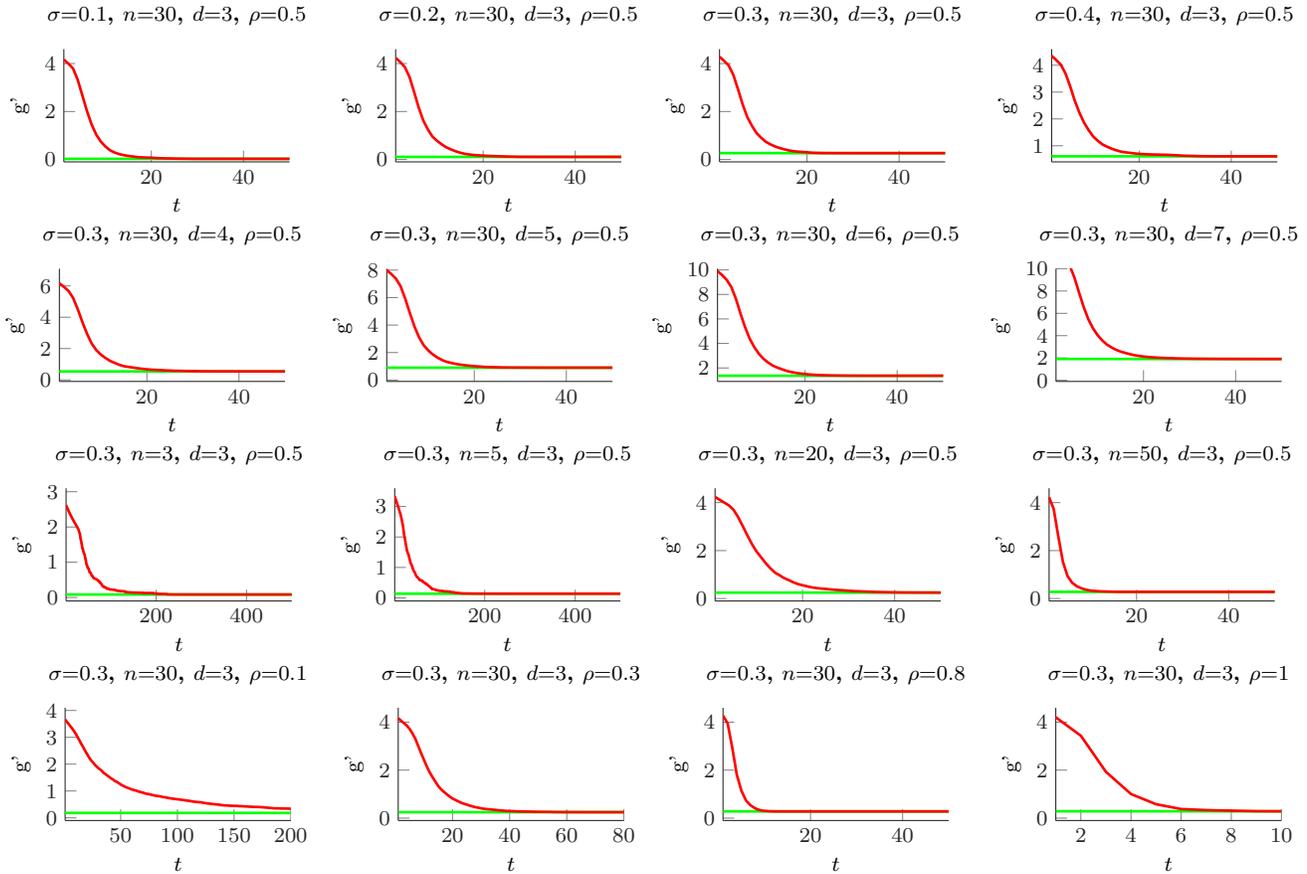

  \centerline{
    \subfigure{%
    % This file was created by matlab2tikz v0.4.7 running on MATLAB 8.4.
% Copyright (c) 2008--2014, Nico Schlömer <nico.schloemer@gmail.com>
% All rights reserved.
% Minimal pgfplots version: 1.3
% 
% The latest updates can be retrieved from
%   http://www.mathworks.com/matlabcentral/fileexchange/22022-matlab2tikz
% where you can also make suggestions and rate matlab2tikz.
% 
\begin{tikzpicture}

\begin{axis}[%
width=3cm,
height=1.5cm,
scale only axis,
every outer x axis line/.append style={white!15!black},
every x tick label/.append style={font=\color{white!15!black}},
xmin=1,
xmax=50,
xlabel={$t$},
every outer y axis line/.append style={white!15!black},
every y tick label/.append style={font=\color{white!15!black}},
ymin=-0.1,
ymax=4.6,
ylabel={g'},
title style={font=\bfseries},
title={$\sigma{=}0.1$, $n{=}30$, $d{=}3$, $\rho{=}0.5$},
axis x line*=bottom,
axis y line*=left
]
\addplot [color=green,solid,line width=1.0pt,forget plot]
  table[row sep=crcr]{%
1	0.0264665133218367\\
2	0.0264665133218367\\
3	0.0264665133218367\\
4	0.0264665133218367\\
5	0.0264665133218367\\
6	0.0264665133218367\\
7	0.0264665133218367\\
8	0.0264665133218367\\
9	0.0264665133218367\\
10	0.0264665133218367\\
11	0.0264665133218367\\
12	0.0264665133218367\\
13	0.0264665133218367\\
14	0.0264665133218367\\
15	0.0264665133218367\\
16	0.0264665133218367\\
17	0.0264665133218367\\
18	0.0264665133218367\\
19	0.0264665133218367\\
20	0.0264665133218367\\
21	0.0264665133218367\\
22	0.0264665133218367\\
23	0.0264665133218367\\
24	0.0264665133218367\\
25	0.0264665133218367\\
26	0.0264665133218367\\
27	0.0264665133218367\\
28	0.0264665133218367\\
29	0.0264665133218367\\
30	0.0264665133218367\\
31	0.0264665133218367\\
32	0.0264665133218367\\
33	0.0264665133218367\\
34	0.0264665133218367\\
35	0.0264665133218367\\
36	0.0264665133218367\\
37	0.0264665133218367\\
38	0.0264665133218367\\
39	0.0264665133218367\\
40	0.0264665133218367\\
41	0.0264665133218367\\
42	0.0264665133218367\\
43	0.0264665133218367\\
44	0.0264665133218367\\
45	0.0264665133218367\\
46	0.0264665133218367\\
47	0.0264665133218367\\
48	0.0264665133218367\\
49	0.0264665133218367\\
50	0.0264665133218367\\
};
\addplot [color=red,solid,line width=1.0pt,forget plot]
  table[row sep=crcr]{%
1	4.1648518295303\\
2	3.99096000057046\\
3	3.79102827525802\\
4	3.32387173537767\\
5	2.67831757190174\\
6	2.02719007477894\\
7	1.47100058163748\\
8	1.04641185713822\\
9	0.751698000116342\\
10	0.543920966098112\\
11	0.386303027369221\\
12	0.29147998208791\\
13	0.223703121013246\\
14	0.181690960063155\\
15	0.147726897105473\\
16	0.124664848217693\\
17	0.0978491386288122\\
18	0.0874676156773447\\
19	0.0791620048297727\\
20	0.0734456937566333\\
21	0.061135459428297\\
22	0.0534327190700933\\
23	0.0455209846462102\\
24	0.0422311141473069\\
25	0.0374650506641064\\
26	0.033094238025453\\
27	0.0320784776313627\\
28	0.0296679666824578\\
29	0.027149390391609\\
30	0.026870702524935\\
31	0.0267084074864977\\
32	0.0266135175219466\\
33	0.0265578282531307\\
34	0.0265250127523152\\
35	0.0265055922733769\\
36	0.0264940473830982\\
37	0.0264871527382754\\
38	0.0264830160500391\\
39	0.0264805224033218\\
40	0.026479012026587\\
41	0.0264780927595593\\
42	0.0264775304727077\\
43	0.0264771847651464\\
44	0.0264769710694122\\
45	0.0264768382221312\\
46	0.0264767551314608\\
47	0.0264767028179913\\
48	0.0264766696432366\\
49	0.0264766484370122\\
50	0.0264766347607332\\
};
\end{axis}
\end{tikzpicture}%
%
 } \hfil 
   \subfigure{%
    % This file was created by matlab2tikz v0.4.7 running on MATLAB 8.4.
% Copyright (c) 2008--2014, Nico Schlömer <nico.schloemer@gmail.com>
% All rights reserved.
% Minimal pgfplots version: 1.3
% 
% The latest updates can be retrieved from
%   http://www.mathworks.com/matlabcentral/fileexchange/22022-matlab2tikz
% where you can also make suggestions and rate matlab2tikz.
% 
\begin{tikzpicture}

\begin{axis}[%
width=3cm,
height=1.5cm,
scale only axis,
every outer x axis line/.append style={white!15!black},
every x tick label/.append style={font=\color{white!15!black}},
xmin=1,
xmax=50,
xlabel={$t$},
every outer y axis line/.append style={white!15!black},
every y tick label/.append style={font=\color{white!15!black}},
ymin=-0.1,
ymax=4.6,
ylabel={g'},
title style={font=\bfseries},
title={$\sigma{=}0.2$, $n{=}30$, $d{=}3$, $\rho{=}0.5$},
axis x line*=bottom,
axis y line*=left
]
\addplot [color=green,solid,line width=1.0pt,forget plot]
  table[row sep=crcr]{%
1	0.107672341078717\\
2	0.107672341078717\\
3	0.107672341078717\\
4	0.107672341078717\\
5	0.107672341078717\\
6	0.107672341078717\\
7	0.107672341078717\\
8	0.107672341078717\\
9	0.107672341078717\\
10	0.107672341078717\\
11	0.107672341078717\\
12	0.107672341078717\\
13	0.107672341078717\\
14	0.107672341078717\\
15	0.107672341078717\\
16	0.107672341078717\\
17	0.107672341078717\\
18	0.107672341078717\\
19	0.107672341078717\\
20	0.107672341078717\\
21	0.107672341078717\\
22	0.107672341078717\\
23	0.107672341078717\\
24	0.107672341078717\\
25	0.107672341078717\\
26	0.107672341078717\\
27	0.107672341078717\\
28	0.107672341078717\\
29	0.107672341078717\\
30	0.107672341078717\\
31	0.107672341078717\\
32	0.107672341078717\\
33	0.107672341078717\\
34	0.107672341078717\\
35	0.107672341078717\\
36	0.107672341078717\\
37	0.107672341078717\\
38	0.107672341078717\\
39	0.107672341078717\\
40	0.107672341078717\\
41	0.107672341078717\\
42	0.107672341078717\\
43	0.107672341078717\\
44	0.107672341078717\\
45	0.107672341078717\\
46	0.107672341078717\\
47	0.107672341078717\\
48	0.107672341078717\\
49	0.107672341078717\\
50	0.107672341078717\\
};
\addplot [color=red,solid,line width=1.0pt,forget plot]
  table[row sep=crcr]{%
1	4.23724956241858\\
2	4.05700875562371\\
3	3.85079764004175\\
4	3.43990253832013\\
5	2.81553695816291\\
6	2.16464024192603\\
7	1.62402530389292\\
8	1.24811390634684\\
9	0.952225177938504\\
10	0.77487671822672\\
11	0.630268385809698\\
12	0.506119964173984\\
13	0.420182057407215\\
14	0.343416024080157\\
15	0.282451938951071\\
16	0.236014688303009\\
17	0.208148372171074\\
18	0.186075950118026\\
19	0.171031779911492\\
20	0.159098161932436\\
21	0.150285783977082\\
22	0.143647596297795\\
23	0.135619818352058\\
24	0.129463334509765\\
25	0.123680079937022\\
26	0.117671349292212\\
27	0.113088539673038\\
28	0.109878861643778\\
29	0.109163578231998\\
30	0.108672043514036\\
31	0.108340978588226\\
32	0.10812172931088\\
33	0.107978341248193\\
34	0.107885347060372\\
35	0.107825337370295\\
36	0.107786712338662\\
37	0.107761872677881\\
38	0.107745892338612\\
39	0.107735598296716\\
40	0.10772895395721\\
41	0.107724654350074\\
42	0.107721863625854\\
43	0.107720046099946\\
44	0.107718857997199\\
45	0.107718078266315\\
46	0.107717564415873\\
47	0.107717224326841\\
48	0.107716998249941\\
49	0.107716847292345\\
50	0.107716746040422\\
};
\end{axis}
\end{tikzpicture}%
%
 }\hfil  
   \subfigure{%
    % This file was created by matlab2tikz v0.4.7 running on MATLAB 8.4.
% Copyright (c) 2008--2014, Nico Schlömer <nico.schloemer@gmail.com>
% All rights reserved.
% Minimal pgfplots version: 1.3
% 
% The latest updates can be retrieved from
%   http://www.mathworks.com/matlabcentral/fileexchange/22022-matlab2tikz
% where you can also make suggestions and rate matlab2tikz.
% 
\begin{tikzpicture}

\begin{axis}[%
width=3cm,
height=1.5cm,
scale only axis,
every outer x axis line/.append style={white!15!black},
every x tick label/.append style={font=\color{white!15!black}},
xmin=1,
xmax=50,
xlabel={$t$},
every outer y axis line/.append style={white!15!black},
every y tick label/.append style={font=\color{white!15!black}},
ymin=-0.1,
ymax=4.6,
ylabel={g'},
title style={font=\bfseries},
title={$\sigma{=}0.3$, $n{=}30$, $d{=}3$, $\rho{=}0.5$},
axis x line*=bottom,
axis y line*=left
]
\addplot [color=green,solid,line width=1.0pt,forget plot]
  table[row sep=crcr]{%
1	0.264709587397121\\
2	0.264709587397121\\
3	0.264709587397121\\
4	0.264709587397121\\
5	0.264709587397121\\
6	0.264709587397121\\
7	0.264709587397121\\
8	0.264709587397121\\
9	0.264709587397121\\
10	0.264709587397121\\
11	0.264709587397121\\
12	0.264709587397121\\
13	0.264709587397121\\
14	0.264709587397121\\
15	0.264709587397121\\
16	0.264709587397121\\
17	0.264709587397121\\
18	0.264709587397121\\
19	0.264709587397121\\
20	0.264709587397121\\
21	0.264709587397121\\
22	0.264709587397121\\
23	0.264709587397121\\
24	0.264709587397121\\
25	0.264709587397121\\
26	0.264709587397121\\
27	0.264709587397121\\
28	0.264709587397121\\
29	0.264709587397121\\
30	0.264709587397121\\
31	0.264709587397121\\
32	0.264709587397121\\
33	0.264709587397121\\
34	0.264709587397121\\
35	0.264709587397121\\
36	0.264709587397121\\
37	0.264709587397121\\
38	0.264709587397121\\
39	0.264709587397121\\
40	0.264709587397121\\
41	0.264709587397121\\
42	0.264709587397121\\
43	0.264709587397121\\
44	0.264709587397121\\
45	0.264709587397121\\
46	0.264709587397121\\
47	0.264709587397121\\
48	0.264709587397121\\
49	0.264709587397121\\
50	0.264709587397121\\
};
\addplot [color=red,solid,line width=1.0pt,forget plot]
  table[row sep=crcr]{%
1	4.28314885338754\\
2	4.10976009389633\\
3	3.91527202574918\\
4	3.5231986405821\\
5	2.93454767727169\\
6	2.30880424199581\\
7	1.77818783642306\\
8	1.4157761910269\\
9	1.10891224642569\\
10	0.912152498615187\\
11	0.74616441214554\\
12	0.63412914755926\\
13	0.54782964093325\\
14	0.48098930297048\\
15	0.426794418084136\\
16	0.379741249916324\\
17	0.345087061093633\\
18	0.326153209795169\\
19	0.313682499315902\\
20	0.298303539391091\\
21	0.281792273076201\\
22	0.27463181475994\\
23	0.270040542059557\\
24	0.268149229516274\\
25	0.266928274453751\\
26	0.266149418128028\\
27	0.265655857161369\\
28	0.265344191050473\\
29	0.265147704891089\\
30	0.265023876799581\\
31	0.264945796770226\\
32	0.264896504564659\\
33	0.264865333987464\\
34	0.264845582635503\\
35	0.264833038222778\\
36	0.264825051120257\\
37	0.26481995228375\\
38	0.264816688431688\\
39	0.264814593446808\\
40	0.264813245059309\\
41	0.264812374894482\\
42	0.264811811924992\\
43	0.264811446854208\\
44	0.264811209633852\\
45	0.264811055238469\\
46	0.264810954641326\\
47	0.264810889074542\\
48	0.264810846368019\\
49	0.264810818607609\\
50	0.264810800632739\\
};
\end{axis}
\end{tikzpicture}%
%
 }  \hfil
   \subfigure{%
    % This file was created by matlab2tikz v0.4.7 running on MATLAB 8.4.
% Copyright (c) 2008--2014, Nico Schlömer <nico.schloemer@gmail.com>
% All rights reserved.
% Minimal pgfplots version: 1.3
% 
% The latest updates can be retrieved from
%   http://www.mathworks.com/matlabcentral/fileexchange/22022-matlab2tikz
% where you can also make suggestions and rate matlab2tikz.
% 
\begin{tikzpicture}

\begin{axis}[%
width=3cm,
height=1.5cm,
scale only axis,
every outer x axis line/.append style={white!15!black},
every x tick label/.append style={font=\color{white!15!black}},
xmin=1,
xmax=50,
xlabel={$t$},
every outer y axis line/.append style={white!15!black},
every y tick label/.append style={font=\color{white!15!black}},
ymin=0.4,
ymax=4.6,
ylabel={g'},
title style={font=\bfseries},
title={$\sigma{=}0.4$, $n{=}30$, $d{=}3$, $\rho{=}0.5$},
axis x line*=bottom,
axis y line*=left
]
\addplot [color=green,solid,line width=1.0pt,forget plot]
  table[row sep=crcr]{%
1	0.606992550043298\\
2	0.606992550043298\\
3	0.606992550043298\\
4	0.606992550043298\\
5	0.606992550043298\\
6	0.606992550043298\\
7	0.606992550043298\\
8	0.606992550043298\\
9	0.606992550043298\\
10	0.606992550043298\\
11	0.606992550043298\\
12	0.606992550043298\\
13	0.606992550043298\\
14	0.606992550043298\\
15	0.606992550043298\\
16	0.606992550043298\\
17	0.606992550043298\\
18	0.606992550043298\\
19	0.606992550043298\\
20	0.606992550043298\\
21	0.606992550043298\\
22	0.606992550043298\\
23	0.606992550043298\\
24	0.606992550043298\\
25	0.606992550043298\\
26	0.606992550043298\\
27	0.606992550043298\\
28	0.606992550043298\\
29	0.606992550043298\\
30	0.606992550043298\\
31	0.606992550043298\\
32	0.606992550043298\\
33	0.606992550043298\\
34	0.606992550043298\\
35	0.606992550043298\\
36	0.606992550043298\\
37	0.606992550043298\\
38	0.606992550043298\\
39	0.606992550043298\\
40	0.606992550043298\\
41	0.606992550043298\\
42	0.606992550043298\\
43	0.606992550043298\\
44	0.606992550043298\\
45	0.606992550043298\\
46	0.606992550043298\\
47	0.606992550043298\\
48	0.606992550043298\\
49	0.606992550043298\\
50	0.606992550043298\\
};
\addplot [color=red,solid,line width=1.0pt,forget plot]
  table[row sep=crcr]{%
1	4.35239106032275\\
2	4.19294834227796\\
3	4.02030084644957\\
4	3.69613891146228\\
5	3.2261489004752\\
6	2.68758531336449\\
7	2.24373984793017\\
8	1.87205344915119\\
9	1.59129264163477\\
10	1.34949407813891\\
11	1.1886068034273\\
12	1.05403626897585\\
13	0.967639192592711\\
14	0.892640302195333\\
15	0.833098653024433\\
16	0.783700481482091\\
17	0.755451495477486\\
18	0.735949700231112\\
19	0.714524493962141\\
20	0.698486500467229\\
21	0.688013583403919\\
22	0.681553833224228\\
23	0.673949392872552\\
24	0.672233290111295\\
25	0.66528413271868\\
26	0.664003259944955\\
27	0.655560094835767\\
28	0.647092491697657\\
29	0.636523511147928\\
30	0.627066734088926\\
31	0.624703719176872\\
32	0.618459512443653\\
33	0.614415782135353\\
34	0.611889377042122\\
35	0.611116142710582\\
36	0.608709460589607\\
37	0.608271242786968\\
38	0.607955907697141\\
39	0.607731277222336\\
40	0.607572327309936\\
41	0.607460295094699\\
42	0.607381496983426\\
43	0.607326136272334\\
44	0.607287273662237\\
45	0.607260018314554\\
46	0.607240928017302\\
47	0.607227579289947\\
48	0.607218264585321\\
49	0.607211780315773\\
50	0.607207278378067\\
};
\end{axis}
\end{tikzpicture}%
%
 }
    } \vspace{-0.4cm}
 	\centerline{
    \subfigure{%
    % This file was created by matlab2tikz v0.4.7 running on MATLAB 8.4.
% Copyright (c) 2008--2014, Nico Schlömer <nico.schloemer@gmail.com>
% All rights reserved.
% Minimal pgfplots version: 1.3
% 
% The latest updates can be retrieved from
%   http://www.mathworks.com/matlabcentral/fileexchange/22022-matlab2tikz
% where you can also make suggestions and rate matlab2tikz.
% 
\begin{tikzpicture}

\begin{axis}[%
width=3cm,
height=1.5cm,
scale only axis,
every outer x axis line/.append style={white!15!black},
every x tick label/.append style={font=\color{white!15!black}},
xmin=1,
xmax=50,
xlabel={$t$},
every outer y axis line/.append style={white!15!black},
every y tick label/.append style={font=\color{white!15!black}},
ymin=-0.1,
ymax=7.1,
ylabel={g'},
title style={font=\bfseries},
title={$\sigma{=}0.3$, $n{=}30$, $d{=}4$, $\rho{=}0.5$},
axis x line*=bottom,
axis y line*=left
]
\addplot [color=green,solid,line width=1.0pt,forget plot]
  table[row sep=crcr]{%
1	0.539075283787872\\
2	0.539075283787872\\
3	0.539075283787872\\
4	0.539075283787872\\
5	0.539075283787872\\
6	0.539075283787872\\
7	0.539075283787872\\
8	0.539075283787872\\
9	0.539075283787872\\
10	0.539075283787872\\
11	0.539075283787872\\
12	0.539075283787872\\
13	0.539075283787872\\
14	0.539075283787872\\
15	0.539075283787872\\
16	0.539075283787872\\
17	0.539075283787872\\
18	0.539075283787872\\
19	0.539075283787872\\
20	0.539075283787872\\
21	0.539075283787872\\
22	0.539075283787872\\
23	0.539075283787872\\
24	0.539075283787872\\
25	0.539075283787872\\
26	0.539075283787872\\
27	0.539075283787872\\
28	0.539075283787872\\
29	0.539075283787872\\
30	0.539075283787872\\
31	0.539075283787872\\
32	0.539075283787872\\
33	0.539075283787872\\
34	0.539075283787872\\
35	0.539075283787872\\
36	0.539075283787872\\
37	0.539075283787872\\
38	0.539075283787872\\
39	0.539075283787872\\
40	0.539075283787872\\
41	0.539075283787872\\
42	0.539075283787872\\
43	0.539075283787872\\
44	0.539075283787872\\
45	0.539075283787872\\
46	0.539075283787872\\
47	0.539075283787872\\
48	0.539075283787872\\
49	0.539075283787872\\
50	0.539075283787872\\
};
\addplot [color=red,solid,line width=1.0pt,forget plot]
  table[row sep=crcr]{%
1	6.16763102102221\\
2	5.93453974038169\\
3	5.67221644987925\\
4	5.21730257255899\\
5	4.50087777083865\\
6	3.68132606853918\\
7	2.94147596373776\\
8	2.32348155104677\\
9	1.90149089684367\\
10	1.61041620116761\\
11	1.39588491437454\\
12	1.21108769943726\\
13	1.08297102662663\\
14	0.961072134046987\\
15	0.871791908185517\\
16	0.827318924994808\\
17	0.78264069235527\\
18	0.740294101197827\\
19	0.704307207269796\\
20	0.6730511109535\\
21	0.649263935996696\\
22	0.63614378868886\\
23	0.621737915709254\\
24	0.604414012455548\\
25	0.594219444891788\\
26	0.58047005325917\\
27	0.576749648756518\\
28	0.569594094875303\\
29	0.562436610743966\\
30	0.556629917697513\\
31	0.551815212669128\\
32	0.546210406622749\\
33	0.544920786722051\\
34	0.543905978644848\\
35	0.543118449450569\\
36	0.541212208110059\\
37	0.540631445728513\\
38	0.540219295761073\\
39	0.539927936704987\\
40	0.539722984409956\\
41	0.53958004214639\\
42	0.539481418390633\\
43	0.539414116434346\\
44	0.539368637896629\\
45	0.539338154503536\\
46	0.539317852598022\\
47	0.539304398997024\\
48	0.539295518958691\\
49	0.539289677122662\\
50	0.539285845528232\\
};
\end{axis}
\end{tikzpicture}%
%
 } \hfil 
   \subfigure{%
    % This file was created by matlab2tikz v0.4.7 running on MATLAB 8.4.
% Copyright (c) 2008--2014, Nico Schlömer <nico.schloemer@gmail.com>
% All rights reserved.
% Minimal pgfplots version: 1.3
% 
% The latest updates can be retrieved from
%   http://www.mathworks.com/matlabcentral/fileexchange/22022-matlab2tikz
% where you can also make suggestions and rate matlab2tikz.
% 
\begin{tikzpicture}

\begin{axis}[%
width=3cm,
height=1.5cm,
scale only axis,
every outer x axis line/.append style={white!15!black},
every x tick label/.append style={font=\color{white!15!black}},
xmin=1,
xmax=50,
xlabel={$t$},
every outer y axis line/.append style={white!15!black},
every y tick label/.append style={font=\color{white!15!black}},
ymin=-0.1,
ymax=8.1,
ylabel={g'},
title style={font=\bfseries},
title={$\sigma{=}0.3$, $n{=}30$, $d{=}5$, $\rho{=}0.5$},
axis x line*=bottom,
axis y line*=left
]
\addplot [color=green,solid,line width=1.0pt,forget plot]
  table[row sep=crcr]{%
1	0.899109615871829\\
2	0.899109615871829\\
3	0.899109615871829\\
4	0.899109615871829\\
5	0.899109615871829\\
6	0.899109615871829\\
7	0.899109615871829\\
8	0.899109615871829\\
9	0.899109615871829\\
10	0.899109615871829\\
11	0.899109615871829\\
12	0.899109615871829\\
13	0.899109615871829\\
14	0.899109615871829\\
15	0.899109615871829\\
16	0.899109615871829\\
17	0.899109615871829\\
18	0.899109615871829\\
19	0.899109615871829\\
20	0.899109615871829\\
21	0.899109615871829\\
22	0.899109615871829\\
23	0.899109615871829\\
24	0.899109615871829\\
25	0.899109615871829\\
26	0.899109615871829\\
27	0.899109615871829\\
28	0.899109615871829\\
29	0.899109615871829\\
30	0.899109615871829\\
31	0.899109615871829\\
32	0.899109615871829\\
33	0.899109615871829\\
34	0.899109615871829\\
35	0.899109615871829\\
36	0.899109615871829\\
37	0.899109615871829\\
38	0.899109615871829\\
39	0.899109615871829\\
40	0.899109615871829\\
41	0.899109615871829\\
42	0.899109615871829\\
43	0.899109615871829\\
44	0.899109615871829\\
45	0.899109615871829\\
46	0.899109615871829\\
47	0.899109615871829\\
48	0.899109615871829\\
49	0.899109615871829\\
50	0.899109615871829\\
};
\addplot [color=red,solid,line width=1.0pt,forget plot]
  table[row sep=crcr]{%
1	7.99648746349688\\
2	7.71123019815341\\
3	7.37316317341311\\
4	6.82588845021207\\
5	5.92619242165445\\
6	4.89274600227118\\
7	3.91093890812117\\
8	3.11862629406964\\
9	2.54475939354803\\
10	2.15590382326316\\
11	1.86755160129359\\
12	1.63574403138596\\
13	1.45594382067725\\
14	1.32746121963852\\
15	1.24438522824186\\
16	1.17493016542747\\
17	1.12034368949\\
18	1.07748133148425\\
19	1.04232170909279\\
20	1.0113519587254\\
21	0.988985093106483\\
22	0.963221020318925\\
23	0.947663250037316\\
24	0.935549174750984\\
25	0.92612759739377\\
26	0.91730144844556\\
27	0.91335669618191\\
28	0.909121103146252\\
29	0.906780938324117\\
30	0.90507290917951\\
31	0.903855371634483\\
32	0.90300596152523\\
33	0.902423656458995\\
34	0.900494863863649\\
35	0.90014680901865\\
36	0.899920416318783\\
37	0.899773090462819\\
38	0.8996769383012\\
39	0.899613961229533\\
40	0.899572567616796\\
41	0.899545271479192\\
42	0.899527218236023\\
43	0.899515246202935\\
44	0.899507287913834\\
45	0.899501986435557\\
46	0.899498448195105\\
47	0.899496082961497\\
48	0.899494499783948\\
49	0.899493439044768\\
50	0.899492727931443\\
};
\end{axis}
\end{tikzpicture}%
%
 }\hfil  
   \subfigure{%
    % This file was created by matlab2tikz v0.4.7 running on MATLAB 8.4.
% Copyright (c) 2008--2014, Nico Schlömer <nico.schloemer@gmail.com>
% All rights reserved.
% Minimal pgfplots version: 1.3
% 
% The latest updates can be retrieved from
%   http://www.mathworks.com/matlabcentral/fileexchange/22022-matlab2tikz
% where you can also make suggestions and rate matlab2tikz.
% 
\begin{tikzpicture}

\begin{axis}[%
width=3cm,
height=1.5cm,
scale only axis,
every outer x axis line/.append style={white!15!black},
every x tick label/.append style={font=\color{white!15!black}},
xmin=1,
xmax=50,
xlabel={$t$},
every outer y axis line/.append style={white!15!black},
every y tick label/.append style={font=\color{white!15!black}},
ymin=0.9,
ymax=10.1,
ylabel={g'},
title style={font=\bfseries},
title={$\sigma{=}0.3$, $n{=}30$, $d{=}6$, $\rho{=}0.5$},
axis x line*=bottom,
axis y line*=left
]
\addplot [color=green,solid,line width=1.0pt,forget plot]
  table[row sep=crcr]{%
1	1.36242635569454\\
2	1.36242635569454\\
3	1.36242635569454\\
4	1.36242635569454\\
5	1.36242635569454\\
6	1.36242635569454\\
7	1.36242635569454\\
8	1.36242635569454\\
9	1.36242635569454\\
10	1.36242635569454\\
11	1.36242635569454\\
12	1.36242635569454\\
13	1.36242635569454\\
14	1.36242635569454\\
15	1.36242635569454\\
16	1.36242635569454\\
17	1.36242635569454\\
18	1.36242635569454\\
19	1.36242635569454\\
20	1.36242635569454\\
21	1.36242635569454\\
22	1.36242635569454\\
23	1.36242635569454\\
24	1.36242635569454\\
25	1.36242635569454\\
26	1.36242635569454\\
27	1.36242635569454\\
28	1.36242635569454\\
29	1.36242635569454\\
30	1.36242635569454\\
31	1.36242635569454\\
32	1.36242635569454\\
33	1.36242635569454\\
34	1.36242635569454\\
35	1.36242635569454\\
36	1.36242635569454\\
37	1.36242635569454\\
38	1.36242635569454\\
39	1.36242635569454\\
40	1.36242635569454\\
41	1.36242635569454\\
42	1.36242635569454\\
43	1.36242635569454\\
44	1.36242635569454\\
45	1.36242635569454\\
46	1.36242635569454\\
47	1.36242635569454\\
48	1.36242635569454\\
49	1.36242635569454\\
50	1.36242635569454\\
};
\addplot [color=red,solid,line width=1.0pt,forget plot]
  table[row sep=crcr]{%
1	9.90374639041393\\
2	9.57729713813464\\
3	9.18420518907162\\
4	8.59959932202767\\
5	7.63195439433428\\
6	6.42154348146885\\
7	5.30482308802894\\
8	4.37601360854144\\
9	3.68532388181084\\
10	3.14562733861119\\
11	2.73662383025714\\
12	2.43789971231008\\
13	2.19040363200046\\
14	2.02411075370333\\
15	1.89039008120283\\
16	1.75920318192144\\
17	1.66604280294703\\
18	1.59905145164628\\
19	1.54653816788635\\
20	1.50516904487264\\
21	1.46536145496381\\
22	1.44084498602901\\
23	1.41812128728385\\
24	1.40423287264917\\
25	1.39337336672392\\
26	1.38542773227659\\
27	1.3785065629712\\
28	1.3750479874658\\
29	1.37105936335637\\
30	1.36795453169349\\
31	1.36645898085862\\
32	1.36540170671707\\
33	1.36466285247235\\
34	1.36415075881136\\
35	1.36379782798082\\
36	1.36355549237975\\
37	1.36338947367153\\
38	1.36327587313116\\
39	1.36319816950022\\
40	1.36314500763323\\
41	1.36310861253825\\
42	1.36308367253924\\
43	1.36306656279159\\
44	1.36305481018912\\
45	1.36304672679497\\
46	1.36304115965182\\
47	1.36303732037841\\
48	1.36303466921222\\
49	1.36303283610631\\
50	1.36303156701068\\
};
\end{axis}
\end{tikzpicture}%
%
 }  \hfil
   \subfigure{%
    % This file was created by matlab2tikz v0.4.7 running on MATLAB 8.4.
% Copyright (c) 2008--2014, Nico Schlömer <nico.schloemer@gmail.com>
% All rights reserved.
% Minimal pgfplots version: 1.3
% 
% The latest updates can be retrieved from
%   http://www.mathworks.com/matlabcentral/fileexchange/22022-matlab2tikz
% where you can also make suggestions and rate matlab2tikz.
% 
\begin{tikzpicture}

\begin{axis}[%
width=3cm,
height=1.5cm,
scale only axis,
every outer x axis line/.append style={white!15!black},
every x tick label/.append style={font=\color{white!15!black}},
xmin=1,
xmax=50,
xlabel={$t$},
every outer y axis line/.append style={white!15!black},
every y tick label/.append style={font=\color{white!15!black}},
ymin=-0.1,
ymax=10,
ylabel={g'},
title style={font=\bfseries},
title={$\sigma{=}0.3$, $n{=}30$, $d{=}7$, $\rho{=}0.5$},
axis x line*=bottom,
axis y line*=left
]
\addplot [color=green,solid,line width=1.0pt,forget plot]
  table[row sep=crcr]{%
1	1.91106237085954\\
2	1.91106237085954\\
3	1.91106237085954\\
4	1.91106237085954\\
5	1.91106237085954\\
6	1.91106237085954\\
7	1.91106237085954\\
8	1.91106237085954\\
9	1.91106237085954\\
10	1.91106237085954\\
11	1.91106237085954\\
12	1.91106237085954\\
13	1.91106237085954\\
14	1.91106237085954\\
15	1.91106237085954\\
16	1.91106237085954\\
17	1.91106237085954\\
18	1.91106237085954\\
19	1.91106237085954\\
20	1.91106237085954\\
21	1.91106237085954\\
22	1.91106237085954\\
23	1.91106237085954\\
24	1.91106237085954\\
25	1.91106237085954\\
26	1.91106237085954\\
27	1.91106237085954\\
28	1.91106237085954\\
29	1.91106237085954\\
30	1.91106237085954\\
31	1.91106237085954\\
32	1.91106237085954\\
33	1.91106237085954\\
34	1.91106237085954\\
35	1.91106237085954\\
36	1.91106237085954\\
37	1.91106237085954\\
38	1.91106237085954\\
39	1.91106237085954\\
40	1.91106237085954\\
41	1.91106237085954\\
42	1.91106237085954\\
43	1.91106237085954\\
44	1.91106237085954\\
45	1.91106237085954\\
46	1.91106237085954\\
47	1.91106237085954\\
48	1.91106237085954\\
49	1.91106237085954\\
50	1.91106237085954\\
};
\addplot [color=red,solid,line width=1.0pt,forget plot]
  table[row sep=crcr]{%
1	11.7481924161392\\
2	11.3722686004206\\
3	10.9150696373376\\
4	10.2457736629152\\
5	9.19644833398169\\
6	7.86558727360288\\
7	6.5781279628608\\
8	5.48737988619235\\
9	4.67726997479646\\
10	4.06279423183971\\
11	3.58721038488231\\
12	3.21649696803674\\
13	2.95623033767625\\
14	2.74120537564895\\
15	2.57436198437157\\
16	2.45810500035146\\
17	2.34503693992969\\
18	2.25156389614094\\
19	2.18403549011642\\
20	2.12955985486062\\
21	2.08701487183268\\
22	2.05713702185886\\
23	2.03102418908485\\
24	2.00822090085067\\
25	1.99414765156911\\
26	1.98205810162286\\
27	1.96511929896281\\
28	1.95392214717721\\
29	1.94464201661004\\
30	1.93983784488529\\
31	1.93478260795138\\
32	1.93210034495019\\
33	1.92724309232422\\
34	1.92430187800633\\
35	1.92274233983957\\
36	1.91968155922941\\
37	1.91851246765703\\
38	1.91751065032495\\
39	1.91665576594524\\
40	1.91593627723186\\
41	1.91446096813171\\
42	1.91385078309797\\
43	1.91336030832276\\
44	1.91298080372896\\
45	1.91269723443946\\
46	1.91249126829503\\
47	1.91234474553815\\
48	1.91224195896987\\
49	1.91217047612422\\
50	1.91212100221106\\
};
\end{axis}
\end{tikzpicture}%
%
 }
    } \vspace{-0.4cm}
    \centerline{
    \subfigure{%
    \input{tikz/errorPlotDistributed-Hmatrix-orthogonal-k-3-d-3-sigma-0.3-graphDensity-0.5-nDraws-\nDraws.tikz}%
 } \hfil 
   \subfigure{%
    \input{tikz/errorPlotDistributed-Hmatrix-orthogonal-k-5-d-3-sigma-0.3-graphDensity-0.5-nDraws-\nDraws.tikz}%
 }\hfil  
   \subfigure{%
    % This file was created by matlab2tikz v0.4.7 running on MATLAB 8.4.
% Copyright (c) 2008--2014, Nico Schlömer <nico.schloemer@gmail.com>
% All rights reserved.
% Minimal pgfplots version: 1.3
% 
% The latest updates can be retrieved from
%   http://www.mathworks.com/matlabcentral/fileexchange/22022-matlab2tikz
% where you can also make suggestions and rate matlab2tikz.
% 
\begin{tikzpicture}

\begin{axis}[%
width=3cm,
height=1.5cm,
scale only axis,
every outer x axis line/.append style={white!15!black},
every x tick label/.append style={font=\color{white!15!black}},
xmin=1,
xmax=50,
xlabel={$t$},
every outer y axis line/.append style={white!15!black},
every y tick label/.append style={font=\color{white!15!black}},
ymin=-0.1,
ymax=4.6,
ylabel={g'},
title style={font=\bfseries},
title={$\sigma{=}0.3$, $n{=}20$, $d{=}3$, $\rho{=}0.5$},
axis x line*=bottom,
axis y line*=left
]
\addplot [color=green,solid,line width=1.0pt,forget plot]
  table[row sep=crcr]{%
1	0.248187730552994\\
2	0.248187730552994\\
3	0.248187730552994\\
4	0.248187730552994\\
5	0.248187730552994\\
6	0.248187730552994\\
7	0.248187730552994\\
8	0.248187730552994\\
9	0.248187730552994\\
10	0.248187730552994\\
11	0.248187730552994\\
12	0.248187730552994\\
13	0.248187730552994\\
14	0.248187730552994\\
15	0.248187730552994\\
16	0.248187730552994\\
17	0.248187730552994\\
18	0.248187730552994\\
19	0.248187730552994\\
20	0.248187730552994\\
21	0.248187730552994\\
22	0.248187730552994\\
23	0.248187730552994\\
24	0.248187730552994\\
25	0.248187730552994\\
26	0.248187730552994\\
27	0.248187730552994\\
28	0.248187730552994\\
29	0.248187730552994\\
30	0.248187730552994\\
31	0.248187730552994\\
32	0.248187730552994\\
33	0.248187730552994\\
34	0.248187730552994\\
35	0.248187730552994\\
36	0.248187730552994\\
37	0.248187730552994\\
38	0.248187730552994\\
39	0.248187730552994\\
40	0.248187730552994\\
41	0.248187730552994\\
42	0.248187730552994\\
43	0.248187730552994\\
44	0.248187730552994\\
45	0.248187730552994\\
46	0.248187730552994\\
47	0.248187730552994\\
48	0.248187730552994\\
49	0.248187730552994\\
50	0.248187730552994\\
};
\addplot [color=red,solid,line width=1.0pt,forget plot]
  table[row sep=crcr]{%
1	4.23757206672648\\
2	4.11734714786887\\
3	4.00170501416502\\
4	3.88532090577722\\
5	3.70235275218087\\
6	3.39673344109661\\
7	3.02732758197908\\
8	2.63808015281842\\
9	2.27225039275487\\
10	1.96458281066073\\
11	1.72150433545926\\
12	1.47506523099756\\
13	1.26271864276884\\
14	1.08328063012201\\
15	0.973672268945812\\
16	0.853037572050432\\
17	0.769134022653555\\
18	0.678147650660757\\
19	0.604606619851223\\
20	0.554921000674161\\
21	0.509908071849586\\
22	0.466134697732602\\
23	0.438192393332252\\
24	0.41577040322002\\
25	0.402191524370575\\
26	0.383752538305146\\
27	0.366946037261441\\
28	0.352360348856203\\
29	0.342417240644055\\
30	0.325904752171544\\
31	0.312451147366622\\
32	0.305921463569654\\
33	0.296917476891912\\
34	0.290593854974671\\
35	0.281889577421513\\
36	0.271124086700937\\
37	0.269226818034739\\
38	0.265289726953887\\
39	0.261076980871885\\
40	0.259871832078868\\
41	0.255305735050741\\
42	0.254262605808499\\
43	0.253377376650744\\
44	0.252625955457072\\
45	0.251987466499139\\
46	0.251444094045687\\
47	0.250980858720909\\
48	0.250585307863391\\
49	0.250247155006474\\
50	0.249957914638229\\
};
\end{axis}
\end{tikzpicture}%
%
 }  \hfil
   \subfigure{%
    % This file was created by matlab2tikz v0.4.7 running on MATLAB 8.4.
% Copyright (c) 2008--2014, Nico Schlömer <nico.schloemer@gmail.com>
% All rights reserved.
% Minimal pgfplots version: 1.3
% 
% The latest updates can be retrieved from
%   http://www.mathworks.com/matlabcentral/fileexchange/22022-matlab2tikz
% where you can also make suggestions and rate matlab2tikz.
% 
\begin{tikzpicture}

\begin{axis}[%
width=3cm,
height=1.5cm,
scale only axis,
every outer x axis line/.append style={white!15!black},
every x tick label/.append style={font=\color{white!15!black}},
xmin=1,
xmax=50,
xlabel={$t$},
every outer y axis line/.append style={white!15!black},
every y tick label/.append style={font=\color{white!15!black}},
ymin=-0.1,
ymax=4.6,
ylabel={g'},
title style={font=\bfseries},
title={$\sigma{=}0.3$, $n{=}50$, $d{=}3$, $\rho{=}0.5$},
axis x line*=bottom,
axis y line*=left
]
\addplot [color=green,solid,line width=1.0pt,forget plot]
  table[row sep=crcr]{%
1	0.279222853781557\\
2	0.279222853781557\\
3	0.279222853781557\\
4	0.279222853781557\\
5	0.279222853781557\\
6	0.279222853781557\\
7	0.279222853781557\\
8	0.279222853781557\\
9	0.279222853781557\\
10	0.279222853781557\\
11	0.279222853781557\\
12	0.279222853781557\\
13	0.279222853781557\\
14	0.279222853781557\\
15	0.279222853781557\\
16	0.279222853781557\\
17	0.279222853781557\\
18	0.279222853781557\\
19	0.279222853781557\\
20	0.279222853781557\\
21	0.279222853781557\\
22	0.279222853781557\\
23	0.279222853781557\\
24	0.279222853781557\\
25	0.279222853781557\\
26	0.279222853781557\\
27	0.279222853781557\\
28	0.279222853781557\\
29	0.279222853781557\\
30	0.279222853781557\\
31	0.279222853781557\\
32	0.279222853781557\\
33	0.279222853781557\\
34	0.279222853781557\\
35	0.279222853781557\\
36	0.279222853781557\\
37	0.279222853781557\\
38	0.279222853781557\\
39	0.279222853781557\\
40	0.279222853781557\\
41	0.279222853781557\\
42	0.279222853781557\\
43	0.279222853781557\\
44	0.279222853781557\\
45	0.279222853781557\\
46	0.279222853781557\\
47	0.279222853781557\\
48	0.279222853781557\\
49	0.279222853781557\\
50	0.279222853781557\\
};
\addplot [color=red,solid,line width=1.0pt,forget plot]
  table[row sep=crcr]{%
1	4.22172434052588\\
2	3.73387294242835\\
3	2.59261960821607\\
4	1.53326787387648\\
5	0.940477631245098\\
6	0.660579872054199\\
7	0.515866906659587\\
8	0.429333998883337\\
9	0.370115725982792\\
10	0.33235308905455\\
11	0.315857823828331\\
12	0.302378382042081\\
13	0.295967734500503\\
14	0.289543427056502\\
15	0.284380523484301\\
16	0.280641878985767\\
17	0.279803197742517\\
18	0.279466268430191\\
19	0.279339255135489\\
20	0.279292154930131\\
21	0.279274732762452\\
22	0.279268267839536\\
23	0.27926585184002\\
24	0.279264939774488\\
25	0.279264591192762\\
26	0.279264456172053\\
27	0.279264403198553\\
28	0.27926438221471\\
29	0.279264373885485\\
30	0.279264370622873\\
31	0.279264369401214\\
32	0.279264368996974\\
33	0.279264368910737\\
34	0.279264368938649\\
35	0.279264368998941\\
36	0.279264369060502\\
37	0.279264369113048\\
38	0.279264369154578\\
39	0.279264369186099\\
40	0.279264369209464\\
41	0.279264369226527\\
42	0.279264369238867\\
43	0.279264369247733\\
44	0.279264369254073\\
45	0.279264369258592\\
46	0.279264369261806\\
47	0.279264369264087\\
48	0.279264369265705\\
49	0.27926436926685\\
50	0.279264369267662\\
};
\end{axis}
\end{tikzpicture}%
%
 } 
    } \vspace{-0.4cm}
    \centerline{
    \subfigure{%
    \input{tikz/errorPlotDistributed-Hmatrix-orthogonal-k-30-d-3-sigma-0.3-graphDensity-0.1-nDraws-\nDraws.tikz}%
 } \hfil 
   \subfigure{%
    % This file was created by matlab2tikz v0.4.7 running on MATLAB 8.4.
% Copyright (c) 2008--2014, Nico Schlömer <nico.schloemer@gmail.com>
% All rights reserved.
% Minimal pgfplots version: 1.3
% 
% The latest updates can be retrieved from
%   http://www.mathworks.com/matlabcentral/fileexchange/22022-matlab2tikz
% where you can also make suggestions and rate matlab2tikz.
% 
\begin{tikzpicture}

\begin{axis}[%
width=3cm,
height=1.5cm,
scale only axis,
every outer x axis line/.append style={white!15!black},
every x tick label/.append style={font=\color{white!15!black}},
xmin=1,
xmax=80,
xlabel={$t$},
every outer y axis line/.append style={white!15!black},
every y tick label/.append style={font=\color{white!15!black}},
ymin=-0.1,
ymax=4.6,
ylabel={g'},
title style={font=\bfseries},
title={$\sigma{=}0.3$, $n{=}30$, $d{=}3$, $\rho{=}0.3$},
axis x line*=bottom,
axis y line*=left
]
\addplot [color=green,solid,line width=1.0pt,forget plot]
  table[row sep=crcr]{%
1	0.247060139275121\\
2	0.247060139275121\\
3	0.247060139275121\\
4	0.247060139275121\\
5	0.247060139275121\\
6	0.247060139275121\\
7	0.247060139275121\\
8	0.247060139275121\\
9	0.247060139275121\\
10	0.247060139275121\\
11	0.247060139275121\\
12	0.247060139275121\\
13	0.247060139275121\\
14	0.247060139275121\\
15	0.247060139275121\\
16	0.247060139275121\\
17	0.247060139275121\\
18	0.247060139275121\\
19	0.247060139275121\\
20	0.247060139275121\\
21	0.247060139275121\\
22	0.247060139275121\\
23	0.247060139275121\\
24	0.247060139275121\\
25	0.247060139275121\\
26	0.247060139275121\\
27	0.247060139275121\\
28	0.247060139275121\\
29	0.247060139275121\\
30	0.247060139275121\\
31	0.247060139275121\\
32	0.247060139275121\\
33	0.247060139275121\\
34	0.247060139275121\\
35	0.247060139275121\\
36	0.247060139275121\\
37	0.247060139275121\\
38	0.247060139275121\\
39	0.247060139275121\\
40	0.247060139275121\\
41	0.247060139275121\\
42	0.247060139275121\\
43	0.247060139275121\\
44	0.247060139275121\\
45	0.247060139275121\\
46	0.247060139275121\\
47	0.247060139275121\\
48	0.247060139275121\\
49	0.247060139275121\\
50	0.247060139275121\\
51	0.247060139275121\\
52	0.247060139275121\\
53	0.247060139275121\\
54	0.247060139275121\\
55	0.247060139275121\\
56	0.247060139275121\\
57	0.247060139275121\\
58	0.247060139275121\\
59	0.247060139275121\\
60	0.247060139275121\\
61	0.247060139275121\\
62	0.247060139275121\\
63	0.247060139275121\\
64	0.247060139275121\\
65	0.247060139275121\\
66	0.247060139275121\\
67	0.247060139275121\\
68	0.247060139275121\\
69	0.247060139275121\\
70	0.247060139275121\\
71	0.247060139275121\\
72	0.247060139275121\\
73	0.247060139275121\\
74	0.247060139275121\\
75	0.247060139275121\\
76	0.247060139275121\\
77	0.247060139275121\\
78	0.247060139275121\\
79	0.247060139275121\\
80	0.247060139275121\\
};
\addplot [color=red,solid,line width=1.0pt,forget plot]
  table[row sep=crcr]{%
1	4.16632365745065\\
2	4.05963975464426\\
3	3.95878256638006\\
4	3.85814732823836\\
5	3.71708707605187\\
6	3.52478342778604\\
7	3.30645333676705\\
8	3.00747678101272\\
9	2.71362439700127\\
10	2.43220634404794\\
11	2.15387966687091\\
12	1.91109529460002\\
13	1.70870617819191\\
14	1.52195450781341\\
15	1.34136783539602\\
16	1.20100128844768\\
17	1.09125899099795\\
18	0.994980499131521\\
19	0.901440602740595\\
20	0.816022443588578\\
21	0.758017940732724\\
22	0.699960691895872\\
23	0.641789948829591\\
24	0.598558782202034\\
25	0.563387424754511\\
26	0.521391813822216\\
27	0.48826672558345\\
28	0.456360521529441\\
29	0.436763257532051\\
30	0.40781444568546\\
31	0.389457262944792\\
32	0.376310616715409\\
33	0.362811582825451\\
34	0.353723493819958\\
35	0.341565589278609\\
36	0.328400901063694\\
37	0.318066733079768\\
38	0.308045610207429\\
39	0.298896674570206\\
40	0.29421868642138\\
41	0.290294550833362\\
42	0.283942237273913\\
43	0.279200267139424\\
44	0.275706018031416\\
45	0.273400187985357\\
46	0.272106532494069\\
47	0.270967698427658\\
48	0.267655098429168\\
49	0.265240057518459\\
50	0.261818921371666\\
51	0.260953528689844\\
52	0.260182282049819\\
53	0.259496576403891\\
54	0.258014212182786\\
55	0.254787332151191\\
56	0.252782572090547\\
57	0.252226811222807\\
58	0.250327240802681\\
59	0.249912152277642\\
60	0.249560956178187\\
61	0.249264080723801\\
62	0.249013207622471\\
63	0.248801174341804\\
64	0.248621857283422\\
65	0.24847005194229\\
66	0.248341357154818\\
67	0.248232066713339\\
68	0.248139070089801\\
69	0.248059763229727\\
70	0.247991969778909\\
71	0.247933872584128\\
72	0.247883954896707\\
73	0.247840950433058\\
74	0.247803801305589\\
75	0.247771622803608\\
76	0.247743674043242\\
77	0.247719333587389\\
78	0.247698079238426\\
79	0.247679471312473\\
80	0.247663138805544\\
};
\end{axis}
\end{tikzpicture}%
%
 }\hfil  
   \subfigure{%
    % This file was created by matlab2tikz v0.4.7 running on MATLAB 8.4.
% Copyright (c) 2008--2014, Nico Schlömer <nico.schloemer@gmail.com>
% All rights reserved.
% Minimal pgfplots version: 1.3
% 
% The latest updates can be retrieved from
%   http://www.mathworks.com/matlabcentral/fileexchange/22022-matlab2tikz
% where you can also make suggestions and rate matlab2tikz.
% 
\begin{tikzpicture}

\begin{axis}[%
width=3cm,
height=1.5cm,
scale only axis,
every outer x axis line/.append style={white!15!black},
every x tick label/.append style={font=\color{white!15!black}},
xmin=1,
xmax=50,
xlabel={$t$},
every outer y axis line/.append style={white!15!black},
every y tick label/.append style={font=\color{white!15!black}},
ymin=-0.1,
ymax=4.6,
ylabel={g'},
title style={font=\bfseries},
title={$\sigma{=}0.3$, $n{=}30$, $d{=}3$, $\rho{=}0.8$},
axis x line*=bottom,
axis y line*=left
]
\addplot [color=green,solid,line width=1.0pt,forget plot]
  table[row sep=crcr]{%
1	0.27968651131728\\
2	0.27968651131728\\
3	0.27968651131728\\
4	0.27968651131728\\
5	0.27968651131728\\
6	0.27968651131728\\
7	0.27968651131728\\
8	0.27968651131728\\
9	0.27968651131728\\
10	0.27968651131728\\
11	0.27968651131728\\
12	0.27968651131728\\
13	0.27968651131728\\
14	0.27968651131728\\
15	0.27968651131728\\
16	0.27968651131728\\
17	0.27968651131728\\
18	0.27968651131728\\
19	0.27968651131728\\
20	0.27968651131728\\
21	0.27968651131728\\
22	0.27968651131728\\
23	0.27968651131728\\
24	0.27968651131728\\
25	0.27968651131728\\
26	0.27968651131728\\
27	0.27968651131728\\
28	0.27968651131728\\
29	0.27968651131728\\
30	0.27968651131728\\
31	0.27968651131728\\
32	0.27968651131728\\
33	0.27968651131728\\
34	0.27968651131728\\
35	0.27968651131728\\
36	0.27968651131728\\
37	0.27968651131728\\
38	0.27968651131728\\
39	0.27968651131728\\
40	0.27968651131728\\
41	0.27968651131728\\
42	0.27968651131728\\
43	0.27968651131728\\
44	0.27968651131728\\
45	0.27968651131728\\
46	0.27968651131728\\
47	0.27968651131728\\
48	0.27968651131728\\
49	0.27968651131728\\
50	0.27968651131728\\
};
\addplot [color=red,solid,line width=1.0pt,forget plot]
  table[row sep=crcr]{%
1	4.26977112643794\\
2	3.95075477999274\\
3	2.94188556625338\\
4	1.79955961051301\\
5	1.13660001193642\\
6	0.730649094118858\\
7	0.530704033641796\\
8	0.405393526078842\\
9	0.332130920220465\\
10	0.301328369618038\\
11	0.286467195700904\\
12	0.281234604988188\\
13	0.280220395127807\\
14	0.279877765040979\\
15	0.279762347878087\\
16	0.279723316533558\\
17	0.279710011179896\\
18	0.279705423011021\\
19	0.279703815388665\\
20	0.279703239304172\\
21	0.279703026181779\\
22	0.279702943772107\\
23	0.279702910006857\\
24	0.279702895183448\\
25	0.279702888182309\\
26	0.279702884643164\\
27	0.27970288275127\\
28	0.279702881697049\\
29	0.279702881092521\\
30	0.279702880739258\\
31	0.279702880530306\\
32	0.279702880405751\\
33	0.279702880331132\\
34	0.27970288028628\\
35	0.279702880259259\\
36	0.279702880242953\\
37	0.279702880233101\\
38	0.279702880227141\\
39	0.279702880223533\\
40	0.279702880221347\\
41	0.279702880220022\\
42	0.279702880219218\\
43	0.27970288021873\\
44	0.279702880218433\\
45	0.279702880218253\\
46	0.279702880218143\\
47	0.279702880218077\\
48	0.279702880218036\\
49	0.279702880218011\\
50	0.279702880217996\\
};
\end{axis}
\end{tikzpicture}%
%
 }  \hfil
   \subfigure{%
    % This file was created by matlab2tikz v0.4.7 running on MATLAB 8.4.
% Copyright (c) 2008--2014, Nico Schlömer <nico.schloemer@gmail.com>
% All rights reserved.
% Minimal pgfplots version: 1.3
% 
% The latest updates can be retrieved from
%   http://www.mathworks.com/matlabcentral/fileexchange/22022-matlab2tikz
% where you can also make suggestions and rate matlab2tikz.
% 
\begin{tikzpicture}

\begin{axis}[%
width=3cm,
height=1.5cm,
scale only axis,
every outer x axis line/.append style={white!15!black},
every x tick label/.append style={font=\color{white!15!black}},
xmin=1,
xmax=10,
xlabel={$t$},
every outer y axis line/.append style={white!15!black},
every y tick label/.append style={font=\color{white!15!black}},
ymin=-0.1,
ymax=4.6,
ylabel={g'},
title style={font=\bfseries},
title={$\sigma{=}0.3$, $n{=}30$, $d{=}3$, $\rho{=}1$},
axis x line*=bottom,
axis y line*=left
]
\addplot [color=green,solid,line width=1.0pt,forget plot]
  table[row sep=crcr]{%
1	0.283922253366453\\
2	0.283922253366453\\
3	0.283922253366453\\
4	0.283922253366453\\
5	0.283922253366453\\
6	0.283922253366453\\
7	0.283922253366453\\
8	0.283922253366453\\
9	0.283922253366453\\
10	0.283922253366453\\
};
\addplot [color=red,solid,line width=1.0pt,forget plot]
  table[row sep=crcr]{%
1	4.20339629566153\\
2	3.43452885553825\\
3	1.92952515476771\\
4	1.00353735817344\\
5	0.578804846342979\\
6	0.378147755098654\\
7	0.333216967024571\\
8	0.309969892173674\\
9	0.287642731084674\\
10	0.284749194030023\\
};
\end{axis}
\end{tikzpicture}%
%
 }
    } \vspace{-0.4cm}
\caption{Normalised error in \eqref{gPrimeFcn} on the vertical axis for the H-matrix method (green) and its distributed version (red) when considering transformations in $O(d)$. The horizontal axis shows the number of steps, where the step size has been chosen as $\epsilon = 0.01$ in each sub-figure. }
    \label{distributedHOrthogonal} 
\end{figure*}

Results of the distributed Z-matrix method are shown in Fig.~\ref{distributedZOrthogonal} and results for the distributed H-matrix method are shown in Fig.~\ref{distributedHOrthogonal}.

\section*{Conclusions}
This worked addressed transitive consistency
of linear inverible transformations between 
Euclidean coordinate systems. Given a set of linear
invertible transformations (or matrices) -- that are not transitively consistent -- the proposed methods synchronize the transformations. This means that they provide transformations that are both transitively consistent and close to the original non-synchronized transformations. 
First two different direct or centralized approaches were proposed. In the first
approach -- the $Z$-matrix approach -- linear algebraic conditions were formulated that 
must hold for transitively consistent transformations. 
Then the sought transformations are obtained from the  
 solution of a least squares problem.
In the second approach -- the $H$-matrix approach -- optimization problems were formulated directly, without
taking a detour via linear algebraic constraints. The
sought transformations are obtained from 
the solution of the optimization problems. 

A Gauss-Newton iterative method was also proposed where the solution from the $H$-matrix method was used as initialization. This method was 
later adapted to the case of affine and Euclidean transformations.
It was shown in numerical simulations that for the case of affine and Euclidean transformations, this approach outperforms the $H$-matrix approach and the 
$Z$-matrix approach. However, for orthogonal transformations 
no improvement is possible over the $H$-matrix method.

In a later part of the paper, for orthogonal matrices, two distributed algorithms were presented. These algorithms share similarities with linear consensus algorithms for distributed averaging. It was shown that these simple consensus-like protocols can be used to provide a 
solution to our problem that is very close to the global 
optimum -- even for noise large in magnitude. 
The proposed methods -- both the direct/centralized and the iterative/distributed -- have been verified to work in numerical experiments for a wide range of parameter 
settings. 

\bibliographystyle{unsrt}       
\bibliography{refExtr}
\end{document}